\documentclass[11pt]{amsart}
\usepackage[utf8]{inputenc}
\usepackage[T1]{fontenc}
\usepackage[english]{babel}
\usepackage{amsfonts,amssymb,amsmath,amsthm,mathtools}
\usepackage{enumitem}
\usepackage{booktabs}
\usepackage{float}
\usepackage{tikz}
\usetikzlibrary{decorations.pathreplacing,calc,positioning}
\usepackage{pgfplots}
\pgfplotsset{compat=1.18}
\usepackage{url,hyperref}
\usepackage[top=3cm,left=3cm,right=3cm,bottom=3cm]{geometry}

\theoremstyle{plain}
\newtheorem{theorem}{Theorem}[section]
\newtheorem{lemma}[theorem]{Lemma}
\newtheorem{proposition}[theorem]{Proposition}
\newtheorem{corollary}[theorem]{Corollary}
\newtheorem{conjecture}[theorem]{Conjecture}

\newtheorem{property}[theorem]{Property}
\theoremstyle{definition}
\newtheorem{definition}[theorem]{Definition}
\newtheorem{example}[theorem]{Example}
\theoremstyle{remark}
\newtheorem{remark}[theorem]{Remark}

\newcommand{\lmid}{\hat{\lambda}}
\newcommand{\Nal}{N_\alpha}
\newcommand{\Fo}{F_2^{\circ}}
\newcommand{\At}{\widetilde{A}_\alpha}
\newcommand{\npos}{\nu_{\le 0}}
\newcommand{\Diag}{\operatorname{Diagonalize}}

\tikzset{
  leaf/.style={circle,draw,fill=white,inner sep=1.7pt},
  dtwo/.style={circle,draw,fill=gray!55,inner sep=1.7pt},
  dtwob/.style={circle,draw,fill=gray!15,inner sep=1.7pt},
  hub/.style={circle,draw,fill=black,inner sep=1.9pt},
  gone/.style={circle,draw,dashed,fill=white,inner sep=1.7pt},
}

\begin{document}

\title[The midrange $A_\alpha$-eigenvalue of a tree]{The midrange $A_\alpha$-eigenvalue of a tree: a sharp threshold beyond $\alpha=\frac{1}{2}$}

\author[E. P. Bartur\'{e}n]{Elvia P. Bartur\'{e}n}
\address{Escuela Profesional de Matem\'{a}tica, Universidad Nacional Mayor de San Marcos, Ciudad Universitaria, Av.~Venezuela s/n, cuadra 34, Lima 1, Per\'{u}}
\email{\tt eperezb@unmsm.edu.pe}
\urladdr{https://orcid.org/0000-0002-0637-3014}

\author[E. R. Oliveira]{Elismar R. Oliveira}
\address{Instituto de Matem\'atica e Estat\'{\i}stica, UFRGS, Porto Alegre, Brazil}
\email{\tt elismar.oliveira@ufrgs.br}
\urladdr{https://orcid.org/0000-0003-2611-0489}

\author[V. Trevisan]{Vilmar Trevisan}
\address{Instituto de Matem\'atica e Estat\'{\i}stica, UFRGS, Porto Alegre, Brazil}
\email{\tt trevisan@mat.ufrgs.br}


\begin{abstract}
For a tree $T$ with $n$ vertices, let $A_\alpha(T)=(1-\alpha)A(T)+\alpha D(T)$ be its Nikiforov matrix and $d_\alpha=\alpha(2-2/n)$ the average of its eigenvalues. It is known that the midrange eigenvalue $\lmid_\alpha(T)$, the $\lceil n/2\rceil$-th smallest one, satisfies $\lmid_\alpha(T)\le d_\alpha$ for $0\le\alpha\le\frac12$. We study the range $\alpha>\frac12$. We prove that $\lmid_\alpha(T)\le d_\alpha$ for all $\alpha\in[0,1]$ if and only if $T$ has at least $\lceil n/2\rceil$ leaves, and that every other tree violates this inequality for all $\alpha\in[\frac{n}{n+1},1]$. This threshold is asymptotically attained by an explicit spider: the largest value $\alpha_{\max}(n)$ at which a tree with fewer leaves still satisfies the inequality obeys $\frac{1}{n+1}<1-\alpha_{\max}(n)\le\frac{1}{n+1}+O(n^{-3})$, and $\alpha_{\max}(n)>1-\frac1n$ for $n=16$ and all $n\ge18$. We then show that the number $p$ of leaves controls when the midrange eigenvalue overtakes $d_\alpha$: if $p\le\lceil n/2\rceil-1$, then $\lmid_\alpha(T)>d_\alpha$ for all $\alpha>\alpha_n(p)$, where $\frac{1}{1-\alpha_n(p)}\le1+n\kappa_n(p)+\min\{3/\kappa_n(p),\sqrt{3n}\}$ for an explicit $\kappa_n(p)\in(0,1)$. Explicit subdivided trees show that this bound is sharp up to an additive constant whenever $p\ge\lceil(\lceil n/2\rceil+3)/3\rceil$, and it yields the asymptotic profile of the latest crossing as a function of the leaf density. We also develop a first-order theory of the crossing based on the forest of degree-two vertices, and we transfer all results to the matrices $B_\beta$ and to the deformed Laplacian.
\end{abstract}

\maketitle

\tableofcontents

\section{Introduction}\label{sec:intro}

The distribution of the eigenvalues of graph matrices inside prescribed intervals of the real line is a classical and prominent subject of spectral graph theory. A striking instance of this problem is the result of Jacobs, Oliveira and Trevisan~\cite{JOT}, who proved that at least half of the Laplacian eigenvalues of a tree are smaller than the average degree $2-\frac{2}{n}$, settling a conjecture of Trevisan, Carvalho, Del Vecchio and Vinagre~\cite{TCDV}; the same result was obtained independently by Sin~\cite{Sin}.

The matrix $A_\alpha(G)=(1-\alpha)A(G)+\alpha D(G)$, $\alpha\in[0,1]$, was introduced by Nikiforov~\cite{Nik} with the purpose of merging the spectral theories of the adjacency matrix $A(G)$ and of the signless Laplacian matrix $Q(G)=D(G)+A(G)$, where $D(G)$ is the diagonal matrix of vertex degrees. We recall that $A_0=A$, that $2A_{1/2}=Q$ and that $A_1=D$. For trees, the Laplacian and the signless Laplacian have the same spectrum, so the value $\alpha=\frac12$ recovers the Laplacian setting of~\cite{JOT}. In~\cite{BO} two of the present authors extended the result of~\cite{JOT} to the $A_\alpha$-matrix: for every tree $T$ with $n$ vertices and every $0\le\alpha\le\frac12$, at least $\lceil n/2\rceil$ of the $A_\alpha$-eigenvalues of $T$ are less than or equal to their average
\[
d_\alpha=\frac{1}{n}\operatorname{tr}A_\alpha(T)=\alpha\Bigl(2-\frac{2}{n}\Bigr).
\]
It was also shown in~\cite{BO} that the same statement may fail for $\alpha>\frac12$, and that there are counterexamples accumulating both at $\alpha=\frac12$ and at $\alpha=1$. The interval $\frac12<\alpha\le1$ therefore exhibits an anomalous behaviour, and the purpose of the present paper is to understand it.

The property studied in~\cite{BO} is conveniently expressed through a single eigenvalue. If $\lambda_1\le\cdots\le\lambda_n$ are the $A_\alpha$-eigenvalues of $T$, the \emph{midrange eigenvalue} is $\lmid_\alpha(T)=\lambda_{\lceil n/2\rceil}$, and the property holds if and only if $\lmid_\alpha(T)\le d_\alpha$ (see Section~\ref{sec:notation}). For each tree we consider the largest value $\alpha_0(T)\in[\frac12,1]$ at which this inequality is still true. In other words, $\alpha_0(T)$ is the moment after which the midrange eigenvalue of $T$ stays strictly above the average $d_\alpha$.

Our main results can be summarized as follows. First, the number $p$ of leaves decides whether the midrange eigenvalue ever overtakes the average: the inequality $\lmid_\alpha(T)\le d_\alpha$ holds for every $\alpha\in[0,1]$ if and only if $p\ge\lceil n/2\rceil$ (Theorem~\ref{thm:dichotomy}); such exceptional trees form an exponentially small proportion of all trees (Remark~\ref{rem:rare}). Every other tree violates the inequality on the whole interval $[\frac{n}{n+1},1]$, independently of its structure (Theorem~\ref{thm:threshold}), and this universal threshold is asymptotically attained: an explicit spider $S_n^\ast$ satisfies $\alpha_0(S_n^\ast)\ge1-\frac{1}{1+n\cos(\pi/\lfloor n/2\rfloor)}$, so that the largest crossing $\alpha_{\max}(n)$ among the non-exceptional trees satisfies $1-\alpha_{\max}(n)=\frac{1}{n+1}+O(n^{-3})$, and $\alpha_{\max}(n)>1-\frac1n$ for $n=16$ and for all $n\ge18$ (Theorem~\ref{thm:spider} and Corollaries~\ref{cor:asymp} and~\ref{cor:false}). Near $\alpha=1$, the position of the crossing is governed, to first order, by the spectrum and the matching number of the forest $F_2$ induced by the vertices of degree two, compared with the \emph{leaf deficit} $j=\lceil n/2\rceil-p$ (Theorem~\ref{thm:firstorder}). Our main quantitative result shows that the number of leaves alone controls the moment of the crossing: if $T$ has $p\le\lceil n/2\rceil-1$ leaves, then $\lmid_\alpha(T)>d_\alpha$ for every $\alpha$ larger than an explicit number $\alpha_n(p)$ with
\[
\frac{1}{1-\alpha_n(p)}\le1+n\,\kappa_n(p)+\min\Bigl\{\frac{3}{\kappa_n(p)},\sqrt{3n}\Bigr\},
\]
where $\kappa_n(p)\in(0,1)$ depends only on $n$ and $p$ (Theorem~\ref{thm:leaves}). This bound is attained up to an additive constant by explicit subdivided trees whenever $p\ge\lceil(\lceil n/2\rceil+3)/3\rceil$ (Theorem~\ref{thm:leaves-sharp} and Corollary~\ref{cor:twosided}), and it gives the asymptotic profile of $n(1-\alpha_{\max}(n,p))$ as a function of the leaf density $p/n$ (Corollary~\ref{cor:profile}). Finally, since the matrices $B_\beta$ and the deformed Laplacian $M_T(s)$ are affine reparametrisations of $A_\alpha$ with positive slope, all these results transfer to those families (Remark~\ref{rem:families}).

The proofs rely on Cauchy interlacing, on the eigenvalue location algorithm of Jacobs and Trevisan~\cite{JT} (see also~\cite{HJT,BO}), on elementary facts about tridiagonal matrices, Chebyshev polynomials and rational recurrences~\cite{OT}, and on a simple combinatorial estimate for the small eigenvalues of a union of paths.

\bigskip
The paper is organized as follows. In Section~\ref{sec:notation} we fix the notation, introduce the centred Nikiforov matrix, the forest of degree-two vertices, and recall the matrix tools that are used throughout the paper, in particular Cauchy interlacing and the algorithm $\Diag$. In Section~\ref{sec:background} we recall the result of~\cite{BO}, prove that every $A_\alpha$-eigenvalue of a tree is strictly increasing in $\alpha$, analyse the extreme case $\alpha=1$ and define the quantities $\alpha_0(T)$ and $\alpha_{\max}(n)$. Section~\ref{sec:exceptional} is devoted to the dichotomy for trees with many leaves. In Section~\ref{sec:threshold} we prove the universal threshold $\frac{n}{n+1}$, and in Section~\ref{sec:spider} we show that it is asymptotically sharp by means of the spider $S_n^\ast$. Section~\ref{sec:firstorder} contains the first-order theory in terms of the degree-two forest. In Section~\ref{sec:leaves} we answer the question of how the crossing depends on the number of leaves. Finally, in Section~\ref{sec:families} we explain how the results transfer to the matrices $B_\beta$ and to the deformed Laplacian.

\section{Notation and preliminaries}\label{sec:notation}

In this section we fix the notation and recall the tools that will be used in the rest of the paper. For the benefit of the reader we include short proofs of the facts that are less standard.

\subsection{Asymptotic notation}\label{sec:asymptotic}

We use the standard Landau notation. Let $f$ and $g$ be real functions of a variable that tends to a limit specified by the context; most of the time this variable is the number of vertices $n\to\infty$, and occasionally it is $\varepsilon=1-\alpha\to0$.
\begin{itemize}
\item $f=O(g)$ means that there are constants $C>0$ and $n_0$ such that $|f(n)|\le C|g(n)|$ for every $n\ge n_0$ (respectively, for all sufficiently small $\varepsilon$). In particular, $f=O(1)$ means that $f$ is bounded. Unless stated otherwise, the constant $C$ is absolute; when it depends on a fixed parameter (such as $\delta$ or $\beta$) this is indicated, and statements that hold \emph{uniformly} over a family of trees mean that the same constants work for all trees in the family.
\item $f=o(g)$ means that $f(n)/g(n)\to0$; in particular, $f=o(1)$ means that $f\to0$.
\item $f\asymp g$ means that $f=O(g)$ and $g=O(f)$, that is, $f$ and $g$ have the same order of magnitude.
\end{itemize}
An expression such as $f=g+O(h)$ means $f-g=O(h)$, and similarly for $o$; the equality sign is read from left to right. An expression such as $a^{\,n+o(n)}$ denotes $a^{\,n+h(n)}$ for some function $h$ with $h(n)=o(n)$, so that $\frac1n\log\bigl(a^{\,n+o(n)}\bigr)\to\log a$. The symbols $\approx$ and $\gtrsim$ are used only informally, in heuristic discussions, to indicate approximate equality and approximate inequality, respectively; they never appear in the statement or in the proof of a result. Finally, a sequence of random variables $X_n$ \emph{converges in probability} to a constant $c$ if $\Pr(|X_n-c|>\delta)\to0$ for every $\delta>0$.

\subsection{Graph matrices}

All graphs are finite, simple and undirected. For a graph $G$ with vertex set $V(G)=\{v_1,\dots,v_n\}$ and edge set $E(G)$ we denote by $A=A(G)$ its adjacency matrix, by $D=D(G)$ the diagonal matrix whose $(i,i)$ entry is the degree $d_{v_i}$ of $v_i$, by $L=D-A$ its Laplacian matrix and by $Q=D+A$ its signless Laplacian matrix. We write $uv$ for the edge joining $u$ and $v$, and $u\sim v$ when $u$ and $v$ are adjacent.

The Laplacian matrix is positive semidefinite. Indeed, for $x\in\mathbb{R}^n$ (indexed by the vertices),
\begin{equation}\label{eq:Lform}
x^\top Lx=\sum_{v}d_vx_v^2-2\sum_{uv\in E(G)}x_ux_v=\sum_{uv\in E(G)}(x_u-x_v)^2\ \ge 0,
\end{equation}
because each $x_v^2$ appears in the last sum once for each edge incident with $v$, that is, $d_v$ times. In particular $x^\top Lx=0$ if and only if $x$ is constant on each connected component of $G$.

The \emph{Nikiforov matrix} (or $A_\alpha$-matrix) of $G$ is
\[
A_\alpha=A_\alpha(G)=(1-\alpha)A+\alpha D,\qquad \alpha\in[0,1].
\]
Since $(1-\alpha)A+\alpha D=A+\alpha(D-A)=D-(1-\alpha)(D-A)$, it can be written in three equivalent ways:
\begin{equation}\label{eq:threeforms}
A_\alpha=(1-\alpha)A+\alpha D=A+\alpha L=D-\varepsilon L,\qquad \varepsilon:=1-\alpha .
\end{equation}
The first form is the definition, the second one shows that $A_\alpha$ increases with $\alpha$ in the sense of quadratic forms (because $L$ is positive semidefinite), and the third one is the natural one near $\alpha=1$, where $A_\alpha$ is a small perturbation of the degree matrix. Throughout the paper we use the abbreviation $\varepsilon=1-\alpha$. The matrix $A_\alpha$ is real and symmetric, hence all its eigenvalues are real. We always list them in \emph{nondecreasing} order,
\[
\lambda_1(A_\alpha)\le\lambda_2(A_\alpha)\le\cdots\le\lambda_n(A_\alpha),
\]
counted with multiplicity (note that in~\cite{BO} the opposite order is used). More generally, for a real symmetric matrix $X$ of order $n$ we write $\lambda_1(X)\le\cdots\le\lambda_n(X)$ for its eigenvalues.

\subsection{\texorpdfstring{The average, the midrange eigenvalue and the property $\mathcal P_\alpha$}{The average, the midrange eigenvalue and the property P-alpha}}

Throughout, $T$ denotes a tree with $n$ vertices. Since a tree with $n$ vertices has exactly $n-1$ edges, and the sum of the degrees counts every edge twice,
\[
\operatorname{tr}A_\alpha(T)=(1-\alpha)\operatorname{tr}A+\alpha\operatorname{tr}D=\alpha\sum_{v}d_v=2\alpha(n-1),
\]
where we used that $A$ has zero diagonal. Therefore the average of the $A_\alpha$-eigenvalues of $T$ (the trace divided by $n$) is
\[
d_\alpha=\alpha\Bigl(2-\frac{2}{n}\Bigr),
\]
that is, $\alpha$ times the average degree $2-\frac2n$ of $T$. We put
\[
k_n:=\left\lceil \frac n2\right\rceil ,
\]
so that $k_n=k$ if $n=2k$ and $k_n=k+1$ if $n=2k+1$.

Since the number of leaves of a tree is an integer, comparisons between it and $\frac n2$ will always be written in terms of $k_n$. We record the elementary translation, paying attention to the parity of $n$: for an integer $p$,
\begin{equation}\label{eq:halfn}
p<\frac n2\iff p\le k_n-1,\qquad\qquad p\ge\frac n2\iff p\ge k_n .
\end{equation}
Indeed, if $n=2k$, then $p<k$ means $p\le k-1=k_n-1$, and $p\ge k=k_n$; note that in this case $p=\frac n2=k_n$ is a possible number of leaves, and it falls into the second alternative. If $n=2k+1$, then $p<k+\frac12$ means $p\le k=k_n-1$, and $p\ge k+\frac12$ means $p\ge k+1=k_n$. Equivalently, $k_n-1=\lceil\frac n2\rceil-1=\lfloor\frac{n-1}{2}\rfloor$ is the largest integer strictly smaller than $\frac n2$, and $k_n$ is the smallest integer larger than or equal to $\frac n2$. In particular, ``at most $k_n-1$ leaves'' and ``at least $k_n$ leaves'' are complementary conditions.

\begin{definition}\label{def:midrange}
Let $T$ be a tree with $n$ vertices and $A_\alpha$-eigenvalues $\lambda_1\le\dots\le\lambda_n$. The \emph{midrange eigenvalue} of $T$ is
\[
\lmid_\alpha(T)=\lambda_{k_n}=\begin{cases}\lambda_k,& n=2k,\\ \lambda_{k+1},& n=2k+1.\end{cases}
\]
We also put
\[
\Nal(T)=\#\{i:\ \lambda_i\le d_\alpha\},
\]
which, in the notation of~\cite{BO}, is the number $m_T(\alpha)(-\infty,d_\alpha]$ of \emph{small} $A_\alpha$-eigenvalues of $T$. We say that $T$ has \emph{property $\mathcal P_\alpha$}, and write $\mathcal P_\alpha(T)$, if
\[
\lmid_\alpha(T)\le d_\alpha .
\]
\end{definition}

\begin{lemma}\label{lem:equivP}
For every tree $T$ with $n$ vertices and every $\alpha$, the following are equivalent: (a) $\lmid_\alpha(T)\le d_\alpha$; (b) $\Nal(T)\ge k_n$; (c) at least $\lceil n/2\rceil$ of the $A_\alpha$-eigenvalues of $T$ are less than or equal to $d_\alpha$.
\end{lemma}

\begin{proof}
(b) and (c) are the same statement. If (a) holds, then $\lambda_1\le\cdots\le\lambda_{k_n}=\lmid_\alpha(T)\le d_\alpha$, so at least $k_n$ eigenvalues are at most $d_\alpha$, which is (b). Conversely, if (b) holds, the set $\{i:\lambda_i\le d_\alpha\}$ has at least $k_n$ elements; since the eigenvalues are ordered, this set is an initial segment $\{1,\dots,\Nal(T)\}$ of the indices, and it contains $k_n$. Hence $\lambda_{k_n}\le d_\alpha$, which is (a).
\end{proof}

Thus the main result of~\cite{BO} says precisely that $\mathcal P_\alpha(T)$ holds for every tree $T$ and every $0\le\alpha\le\frac12$ (Property~\ref{property} below), and the present paper studies the validity of $\mathcal P_\alpha(T)$ for $\frac12<\alpha\le1$.

\subsection{The centred Nikiforov matrix}

It is convenient to shift the matrix by its average.

\begin{definition}
The \emph{centred Nikiforov matrix} of $T$ is
\[
\At=\At(T):=A_\alpha(T)-d_\alpha I .
\]
For a real symmetric matrix $X$ we denote by $\npos(X)$ the number of nonpositive eigenvalues of $X$, counted with multiplicity, and for a set $U$ of indices we denote by $X[U]$ the principal submatrix of $X$ with rows and columns indexed by $U$.
\end{definition}

The eigenvalues of $\At$ are $\lambda_i-d_\alpha$, $1\le i\le n$. Hence $\lambda_i\le d_\alpha$ if and only if the corresponding eigenvalue of $\At$ is nonpositive, and
\begin{equation}\label{eq:Nnpos}
\Nal(T)=\npos(\At(T)).
\end{equation}
By the definition of $A_\alpha$, the diagonal entry of $\At$ at a vertex $v$ is $\alpha d_v-d_\alpha=\alpha\bigl(d_v-2+\frac{2}{n}\bigr)$, and the off-diagonal entry at a pair $u\neq v$ is $(1-\alpha)$ if $uv$ is an edge and $0$ otherwise. Separating the vertices according to their degree we obtain
\begin{equation}\label{eq:entries}
(\At)_{vv}=
\begin{cases}
-\alpha\,\dfrac{n-2}{n}, & d_v=1,\\[2mm]
\dfrac{2\alpha}{n}, & d_v=2,\\[2mm]
\alpha\bigl(d_v-2+\frac2n\bigr), & d_v\ge3,
\end{cases}
\qquad (\At)_{uv}=\begin{cases}\varepsilon, & uv\in E(T),\\ 0, & \text{otherwise},\end{cases}
\end{equation}
for $u\ne v$. We shall often use the abbreviations
\begin{equation}\label{eq:x1c}
x_1=x_1(\alpha):=-\alpha\,\frac{n-2}{n},\qquad c=c(\alpha):=\frac{2\alpha}{n},
\end{equation}
for the diagonal entries at leaves and at vertices of degree $2$, respectively. The notation $x_1$ agrees with~\cite{BO}, where $x_1(\alpha)=\alpha-d_\alpha$ is the first value produced by the algorithm $\Diag$ at a leaf (see Section~\ref{sec:tools}).

\begin{remark}[Three scales]\label{rem:scales}
When $\alpha$ is close to $1$, the diagonal of $\At$ lives on three different scales: the entries at leaves are close to $-1$, the entries at vertices of degree at least $3$ are at least $\alpha(1+\frac{2}{n})$, hence close to or larger than $1$, while the entries at vertices of degree $2$ are only $\frac{2\alpha}{n}$, which is small when $n$ is large. The off-diagonal entries are equal to $\varepsilon=1-\alpha$, which is also small. Hence, near $\alpha=1$, every leaf contributes a negative eigenvalue and every vertex of degree at least $3$ contributes a positive eigenvalue of $\At$, while the sign of the eigenvalues associated with the vertices of degree $2$ is decided by the competition between $\frac{2\alpha}{n}$ and $\varepsilon$. This simple observation drives the whole analysis of the paper.
\end{remark}

\subsection{Degree classes and the forest of degree-two vertices}\label{sec:runs}

For a tree $T$ we denote by $p=p(T)$ the number of leaves (vertices of degree $1$), by $n_2$ the number of vertices of degree $2$ and by $n_{\ge3}$ the number of vertices of degree at least $3$, which we call \emph{hubs}. The corresponding sets of vertices are $V_1$, $V_2$ and $V_{\ge3}$, so that $n=p+n_2+n_{\ge3}$. For a set $U$ of vertices, $T[U]$ denotes the subgraph induced by $U$, and $T-U$ the subgraph induced by the complement of $U$.

Let $F_2=T[V_2]$ be the subgraph induced by the vertices of degree $2$. It is a forest, because it is a subgraph of a tree, and each of its vertices has at most two neighbours in $F_2$, because it has exactly two neighbours in $T$. A connected acyclic graph in which every vertex has degree at most $2$ is a path (starting at a vertex of degree at most one and walking along the unique available edges we visit all its vertices). Therefore $F_2$ is a disjoint union of paths.

\begin{definition}\label{def:runs}
The connected components of $F_2$ are called the \emph{runs} of $T$; their number is denoted by $r=r(T)$, and the number of vertices of a run is its \emph{length}. Let $V_2^\circ\subseteq V_2$ be the set of vertices of degree $2$ having no leaf as a neighbour, and let $\Fo=T[V_2^\circ]$, which is again a disjoint union of paths. Finally we set
\[
j(T):=k_n-p .
\]
\end{definition}

By Lemma~\ref{lem:degree} below, when $T$ has at most $k_n-1$ leaves, $j(T)$ is a positive integer, and it is the number of eigenvalues that the degree-$2$ part of the tree must contribute below $d_\alpha$ in order for $\mathcal P_\alpha$ to hold (the leaves contribute $p$ of them and the hubs essentially none, see Remark~\ref{rem:scales}).

A \emph{matching} of a graph is a set of pairwise disjoint edges, and the \emph{matching number} $\mu(F)$ of a graph $F$ is the maximum size of a matching. For a forest $F$ we write $\theta_1(F)\le\theta_2(F)\le\cdots$ for the eigenvalues of its adjacency matrix $A(F)$.

\begin{figure}[ht]
\centering
\begin{tikzpicture}[scale=1.05]
\node[hub,label=above:$a$] (a) at (0,0) {};
\node[leaf] (l1) at (-1,0.8) {};
\node[leaf] (l2) at (-1,-0.8) {};
\node[dtwo] (r1) at (1,0) {};
\node[dtwo] (r2) at (2,0) {};
\node[dtwo] (r3) at (3,0) {};
\node[hub,label=below:$b$] (b) at (4,0) {};
\node[leaf] (l3) at (4,1.1) {};
\node[dtwo] (s1) at (5,0.6) {};
\node[dtwob] (s2) at (6,0.6) {};
\node[leaf] (l4) at (7,0.6) {};
\node[dtwob] (t1) at (5,-0.7) {};
\node[leaf] (l5) at (6,-0.7) {};
\draw (l1)--(a)--(l2); \draw (a)--(r1)--(r2)--(r3)--(b)--(l3);
\draw (b)--(s1)--(s2)--(l4); \draw (b)--(t1)--(l5);
\begin{scope}[shift={(8.6,0.9)}]
\node[hub] at (0,0) {}; \node[right] at (0.15,0) {\small hub ($d_v\ge3$)};
\node[dtwo] at (0,-0.5) {}; \node[right] at (0.15,-0.5) {\small $V_2^\circ$};
\node[dtwob] at (0,-1) {}; \node[right] at (0.15,-1) {\small $V_2\setminus V_2^\circ$};
\node[leaf] at (0,-1.5) {}; \node[right] at (0.15,-1.5) {\small leaf};
\end{scope}
\end{tikzpicture}
\caption{A tree with $n=13$ vertices, $p=5$ leaves and $k_{13}=7$, so that $j(T)=2$. There are $r=3$ runs: $P_3$ (between $a$ and $b$), $P_2$ and $P_1$. Hence $n_2=6$ and $\mu(F_2)=1+1+0=2$. The spectrum of $F_2$ is $\{-\sqrt2,-1,0,0,1,\sqrt2\}$, so $\theta_2(F_2)=-1$. The vertices of $V_2^\circ$ induce $\Fo=P_3\cup P_1$, and $\mu(\Fo)=1$. The first-order prediction of Section~\ref{sec:firstorder} for this tree is $1-\frac{2}{2-n\theta_j(F_2)}=\frac{13}{15}$.}
\label{fig:example}
\end{figure}
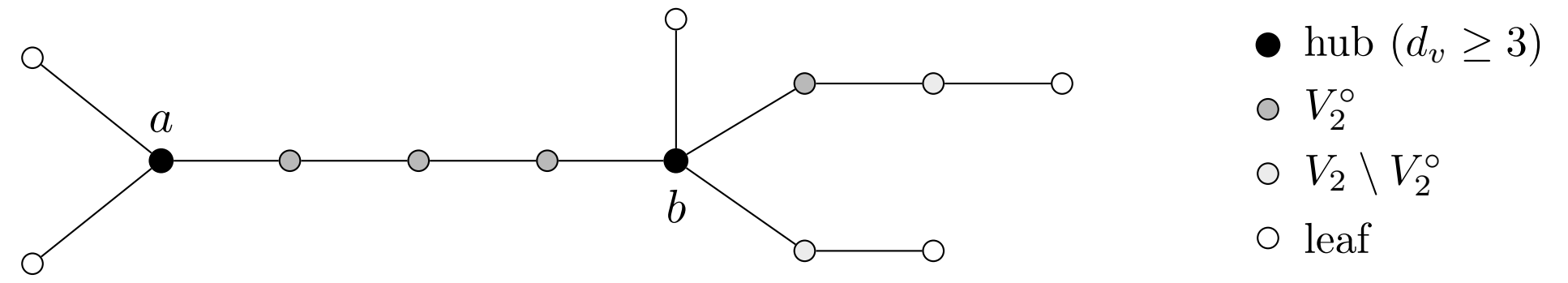

\subsection{Spectra of paths}

We denote by $P_\ell$ the path with $\ell$ vertices $1,2,\dots,\ell$ (in this order).

\begin{lemma}\label{lem:pathspec}
The adjacency eigenvalues of $P_\ell$ are
\begin{equation}\label{eq:pathspec}
2\cos\frac{k\pi}{\ell+1},\qquad k=1,\dots,\ell,
\end{equation}
with eigenvectors $x^{(k)}$ given by $x^{(k)}_i=\sin\frac{ik\pi}{\ell+1}$, $1\le i\le\ell$. Exactly $\lfloor \ell/2\rfloor=\mu(P_\ell)$ of these eigenvalues are negative. Consequently, if $F$ is a disjoint union of paths, then
\begin{equation}\label{eq:matching}
\#\{i:\theta_i(F)<0\}=\mu(F),\qquad\text{hence}\qquad \theta_j(F)<0\iff j\le\mu(F).
\end{equation}
\end{lemma}

\begin{proof}
Put $t=\frac{k\pi}{\ell+1}$ and $x_i=\sin(it)$ for $0\le i\le \ell+1$; note that $x_0=0$ and $x_{\ell+1}=\sin(k\pi)=0$. For $1\le i\le\ell$, the $i$-th entry of $A(P_\ell)x^{(k)}$ is $x_{i-1}+x_{i+1}$ (the terms $x_0$ and $x_{\ell+1}$ vanish, which accounts for the end vertices). By the identity $\sin((i-1)t)+\sin((i+1)t)=2\sin(it)\cos t$ we get $A(P_\ell)x^{(k)}=2\cos t\ x^{(k)}$. The vector $x^{(k)}$ is nonzero because $x_1=\sin t\neq0$. The $\ell$ numbers $2\cos\frac{k\pi}{\ell+1}$, $1\le k\le\ell$, are pairwise distinct because $\cos$ is injective on $(0,\pi)$; hence they are all the eigenvalues of $P_\ell$.

The eigenvalue $2\cos\frac{k\pi}{\ell+1}$ is negative if and only if $\frac{k}{\ell+1}>\frac12$, that is, $k>\frac{\ell+1}{2}$. If $\ell=2s$ this means $k\in\{s+1,\dots,2s\}$, and if $\ell=2s+1$ it means $k\in\{s+2,\dots,2s+1\}$; in both cases there are exactly $s=\lfloor\ell/2\rfloor$ such values of $k$. On the other hand, the edges of a matching of $P_\ell$ cover $2$ vertices each and are disjoint, so a matching has at most $\lfloor\ell/2\rfloor$ edges, and the edges $\{1,2\},\{3,4\},\dots$ form a matching with exactly $\lfloor\ell/2\rfloor$ edges. Thus $\mu(P_\ell)=\lfloor\ell/2\rfloor$.

If $F$ is a disjoint union of paths, the adjacency matrix of $F$ is block diagonal, so its spectrum is the union (with multiplicities) of the spectra of its components. Also, a matching of $F$ is the union of matchings of its components, so $\mu(F)$ is the sum of the matching numbers of the components. Adding the counts above over the components gives the first identity in~\eqref{eq:matching}. The equivalence follows because the eigenvalues $\theta_i(F)$ are ordered.
\end{proof}

\subsection{Matrix tools}\label{sec:tools}

\paragraph{The Loewner order.} For real symmetric matrices $X,Y$ of the same order we write $X\succeq Y$ if $X-Y$ is positive semidefinite, $X\succ0$ if $X$ is positive definite, and we denote by $\|X\|$ the spectral norm of $X$ (the largest absolute value of an eigenvalue of $X$).

\begin{lemma}[Weyl's monotonicity]\label{lem:weyl}
Let $X,Y$ be real symmetric matrices of order $n$.
\begin{enumerate}[label=(W\arabic*)]
\item\label{W1} If $X\succeq Y$, then $\lambda_i(X)\ge\lambda_i(Y)$ for every $i$; in particular $\npos(X)\le\npos(Y)$.
\item\label{W2} $|\lambda_i(X)-\lambda_i(Y)|\le\|X-Y\|$ for every $i$. In particular, the eigenvalues of a symmetric matrix depend continuously on its entries.
\end{enumerate}
\end{lemma}

\begin{proof}
By the Courant--Fischer theorem (see~\cite[Theorem 4.2.6]{HJ}),
\[
\lambda_i(X)=\min_{\dim S=i}\ \max_{x\in S,\ \|x\|=1}x^\top Xx,
\]
where $S$ runs over the subspaces of $\mathbb R^n$ of dimension $i$. If $X\succeq Y$, then $x^\top Xx\ge x^\top Yx$ for every $x$, so for every $S$ the inner maximum for $X$ is at least the inner maximum for $Y$, and the minima satisfy the same inequality. This proves the first claim of~\ref{W1}. If $\lambda_i(X)\le0$, then $\lambda_i(Y)\le\lambda_i(X)\le0$; hence the number of nonpositive eigenvalues of $Y$ is at least that of $X$.
For~\ref{W2}, put $\delta=\|X-Y\|$. The eigenvalues of $X-Y$ lie in $[-\delta,\delta]$, so $Y+\delta I\succeq X\succeq Y-\delta I$, and~\ref{W1} gives $\lambda_i(Y)+\delta\ge\lambda_i(X)\ge\lambda_i(Y)-\delta$. Since $\|X-Y\|$ tends to zero when the entries of $X$ tend to those of $Y$, the eigenvalues are continuous.
\end{proof}

\paragraph{Cauchy interlacing.} If $X$ is a real symmetric matrix of order $n$ and $Y=X[U]$ is a principal submatrix of order $n'$, then (see~\cite[Theorem 4.3.28]{HJ})
\begin{equation}\label{eq:cauchy}
\lambda_i(X)\le\lambda_i(Y)\le\lambda_{i+n-n'}(X),\qquad 1\le i\le n' .
\end{equation}
We shall use it through the following counting consequences.

\begin{lemma}\label{lem:interlacing}
Let $Y=X[U]$ be a principal submatrix of order $n'$ of a real symmetric matrix $X$ of order $n$. Then
\begin{enumerate}[label=(I\arabic*)]
\item\label{I1} $\npos(X)\ge\npos(Y)$;
\item\label{I2} $\npos(X)\le(n-n')+\npos(Y)$. In particular, if $Y$ is positive definite, then $\npos(X)\le n-n'$.
\end{enumerate}
\end{lemma}

\begin{proof}
\ref{I1} Let $q=\npos(Y)$. If $q=0$ there is nothing to prove. Otherwise $\lambda_q(Y)\le0$, and the left inequality in~\eqref{eq:cauchy} gives $\lambda_q(X)\le\lambda_q(Y)\le0$. Since the eigenvalues of $X$ are ordered, $\lambda_1(X),\dots,\lambda_q(X)$ are all nonpositive, so $\npos(X)\ge q$.

\ref{I2} Let $q=\npos(X)$. If $q\le n-n'$ the inequality is trivial. Suppose $q>n-n'$ and put $i=q-(n-n')$, so that $1\le i\le n'$. The right inequality in~\eqref{eq:cauchy} gives $\lambda_i(Y)\le\lambda_{i+n-n'}(X)=\lambda_{q}(X)\le0$, so $Y$ has at least $i=q-(n-n')$ nonpositive eigenvalues, that is, $q\le(n-n')+\npos(Y)$. If $Y\succ0$, then $\npos(Y)=0$.
\end{proof}

\paragraph{Sylvester's law of inertia.} Two real symmetric matrices $X$ and $Y$ are \emph{congruent} if $Y=E^\top XE$ for some invertible matrix $E$. Sylvester's law of inertia (see~\cite[Theorem 4.5.8]{HJ}) states that congruent matrices have the same numbers of positive, negative and zero eigenvalues.

\paragraph{\texorpdfstring{The algorithm $\Diag$.}{The algorithm Diagonalize.}} Our main tool for eigenvalue location is the algorithm $\Diag$ of Jacobs and Trevisan~\cite{JT}, whose versions for other matrices are discussed in~\cite{HJT,BO}. It applies to any real symmetric matrix $M=(m_{uv})$ whose underlying graph is a tree $T$, that is, $m_{uv}\neq0$ for $u\neq v$ only if $uv\in E(T)$. We root $T$ at some vertex and order the vertices $v_1,\dots,v_n$ bottom-up, so that every child precedes its parent. The algorithm is reproduced in Figure~\ref{fig:diagonalize}.

\begin{figure}[ht]
\centering
{\small
{\tt
\begin{tabbing}
aaa\=aaa\=aaa\=aaa\=aaa\=aaa\=aaa\=aaa\= \kill
\> Input: a symmetric matrix $M=(m_{uv})$ with underlying tree $T$\\
\> Input: a bottom-up ordering $v_1,\ldots,v_n$ of $V(T)$\\
\> Input: a real number $x$ \\
\> Output: a diagonal matrix $\operatorname{diag}(a(v_1),\ldots,a(v_n))$ congruent to $M+xI$ \\
\> \\
\> Algorithm $\Diag(M,x)$ \\
\> \> initialize $a(v):=m_{vv}+x$, for all vertices $v$ \\
\> \> {\bf for } $k=1$ to $n$ \\
\> \> \> {\bf if} $v_k$ is a leaf of the rooted tree {\bf then} continue \\
\> \> \> {\bf else if} $a(c)\neq0$ for all children $c$ of $v_k$ {\bf then} \\
\> \> \> \> $a(v_k):=a(v_k)-\sum_c \frac{m_{v_kc}^2}{a(c)}$, summing over all children of $v_k$ \\
\> \> \> {\bf else } \\
\> \> \> \> select one child $c_0$ of $v_k$ for which $a(c_0)=0$ \\
\> \> \> \> $a(v_k):=-\frac{m_{v_kc_0}^2}{2}$ \\
\> \> \> \> $a(c_0):=2$ \\
\> \> \> \> if $v_k$ has a parent $v_\ell$, remove the edge $v_kv_\ell$ \\
\> \> {\bf end loop}
\end{tabbing}
}
}
\caption{The algorithm $\Diag(M,x)$, which produces a diagonal matrix congruent to $M+xI$.}
\label{fig:diagonalize}
\end{figure}

\begin{lemma}[{\cite{JT}}]\label{lem:diag}
The diagonal matrix produced by $\Diag(M,x)$ is congruent to $M+xI$. Consequently, the numbers of negative, zero and positive values $a(v)$ are, respectively, the numbers of eigenvalues of $M$ that are smaller than, equal to and larger than $-x$. Moreover, if the third branch of the algorithm is never used, then the values satisfy
\begin{equation}\label{eq:diag}
a(v)=m_{vv}+x-\sum_{c\ \text{child of}\ v}\frac{m_{vc}^2}{a(c)},
\end{equation}
and $\det(M+xI)=\prod_va(v)$.
\end{lemma}

\begin{proof}
We prove the statement when the third branch is never used, which is the only case needed in this paper; the general case is proved in~\cite{JT} (see also~\cite{HJT}). Put $N=M+xI$; its underlying graph is still $T$. Let $c$ be a leaf of the rooted tree (a vertex without children), with parent $v$, and assume $a(c)=N_{cc}\neq0$. Let $E$ be the identity matrix with the additional entry $-N_{vc}/a(c)$ in position $(c,v)$. The matrix $N'=E^\top NE$ is obtained from $N$ by subtracting $\frac{N_{vc}}{a(c)}$ times row $c$ from row $v$ and the same multiple of column $c$ from column $v$. Since the only nonzero off-diagonal entry in row $c$ is $N_{cv}$ (the vertex $c$ has no other neighbour), this operation annihilates the entries $(v,c)$ and $(c,v)$, replaces $N_{vv}$ by $N_{vv}-N_{vc}^2/a(c)$ and does not change any other entry. Hence $N'$ is the direct sum of the $1\times1$ matrix $(a(c))$ and a symmetric matrix whose underlying graph is the tree $T-c$, with the modified diagonal entry at $v$. Repeating this procedure in the bottom-up order, and subtracting at each vertex the contributions of all its children, we obtain exactly the values~\eqref{eq:diag} and a diagonal matrix $E_{\rm tot}^\top NE_{\rm tot}=\operatorname{diag}(a(v))$, where $E_{\rm tot}$ is a product of matrices of the form $E$, each of determinant $1$. The statement on the signs follows from Sylvester's law of inertia, since the eigenvalues of $N$ are $\lambda_i(M)+x$, and $\det N=\det\operatorname{diag}(a(v))=\prod_v a(v)$ because $\det E_{\rm tot}=1$.
\end{proof}

We will apply Lemma~\ref{lem:diag} to $M=A_\alpha(T)$ and $x=-d_\alpha$, that is, to the matrix $\At(T)=A_\alpha(T)-d_\alpha I$, whose off-diagonal entries on edges are all equal to $\varepsilon$. In this case~\eqref{eq:diag} reads
\[
a(v)=(\At)_{vv}-\varepsilon^2\sum_{c\ \text{child of}\ v}\frac{1}{a(c)},
\]
and, by~\eqref{eq:Nnpos} and Lemma~\ref{lem:diag}, $\Nal(T)$ is the number of nonpositive values $a(v)$. At a leaf $v$ of $T$ (which is a leaf of the rooted tree whenever it is not the root) the algorithm produces $a(v)=\alpha-d_\alpha=x_1(\alpha)$, see~\eqref{eq:x1c}.

\section{Background and first properties}\label{sec:background}

We start by recalling the main result of~\cite{BO}, which is the starting point of this work.

\begin{property}[{\cite[Theorem 1.1]{BO}}]\label{property}
For every tree $T$ with $n$ vertices and every $0\le\alpha\le\frac12$, at least $k_n=\lceil\frac n2\rceil$ of the $A_\alpha$-eigenvalues of $T$ are less than or equal to $d_\alpha$. Equivalently (Lemma~\ref{lem:equivP}), $\lmid_\alpha(T)\le d_\alpha$, that is, $\mathcal P_\alpha(T)$ holds.
\end{property}

The next proposition shows that every $A_\alpha$-eigenvalue of a tree moves to the right when $\alpha$ increases. This is natural in view of the second form in~\eqref{eq:threeforms}, but the strictness requires an argument.

\begin{proposition}[Monotonicity in $\alpha$]\label{prop:monotone}
Let $G$ be a graph with $n$ vertices and let $\lambda_1(\alpha)\le\dots\le\lambda_n(\alpha)$ be the eigenvalues of $A_\alpha(G)=A+\alpha L$, defined for every real $\alpha$. Then each function $\alpha\mapsto\lambda_i(\alpha)$ is nondecreasing on $\mathbb R$. If no connected component of $G$ is regular, then each $\lambda_i(\alpha)$ is strictly increasing on $\mathbb R$. In particular, if $T$ is a tree with $n\ge3$ vertices, every $A_\alpha$-eigenvalue of $T$, and hence $\lmid_\alpha(T)$, is strictly increasing in $\alpha$.
\end{proposition}

\begin{proof}
Let $\alpha<\alpha'$ and put $t=\alpha'-\alpha>0$, $X=A_\alpha$ and $Y=A_{\alpha'}$. By~\eqref{eq:threeforms}, $Y=X+tL$, and $tL\succeq0$ by~\eqref{eq:Lform}. Hence $Y\succeq X$, and~\ref{W1} of Lemma~\ref{lem:weyl} gives $\lambda_i(\alpha')\ge\lambda_i(\alpha)$. This proves monotonicity.

Assume now that no component of $G$ is regular, and suppose, by contradiction, that $\lambda_i(\alpha)=\lambda_i(\alpha')=:\lambda$ for some $i$. Let $u_1,\dots,u_n$ be an orthonormal basis of eigenvectors of $X$, with $Xu_k=\lambda_k(\alpha)u_k$, and let $w_1,\dots,w_n$ be an orthonormal basis of eigenvectors of $Y$, with $Yw_k=\lambda_k(\alpha')w_k$. Consider the subspaces
\[
U'=\operatorname{span}\{w_1,\dots,w_i\},\qquad W'=\operatorname{span}\{u_i,\dots,u_n\}.
\]
Since $\dim U'+\dim W'=i+(n-i+1)=n+1>n$, there is a unit vector $z\in U'\cap W'$. As $z\in U'$, we have $z^\top Yz\le\lambda_i(\alpha')=\lambda$; as $z\in W'$, we have $z^\top Xz\ge\lambda_i(\alpha)=\lambda$. Therefore
\[
\lambda\ \ge\ z^\top Yz\ =\ z^\top Xz+t\,z^\top Lz\ \ge\ \lambda+t\,z^\top Lz .
\]
Since $t>0$ and $z^\top Lz\ge0$, we conclude that $z^\top Lz=0$ and $z^\top Xz=\lambda$. Write $z=\sum_{k\ge i}\gamma_ku_k$. Then $z^\top Xz=\sum_{k\ge i}\gamma_k^2\lambda_k(\alpha)=\lambda=\lambda\sum_{k\ge i}\gamma_k^2$, that is, $\sum_{k\ge i}\gamma_k^2\bigl(\lambda_k(\alpha)-\lambda\bigr)=0$. All the terms of this sum are nonnegative, so $\gamma_k=0$ whenever $\lambda_k(\alpha)>\lambda$, and therefore $Xz=\lambda z$. On the other hand, $z^\top Lz=0$ means, by~\eqref{eq:Lform}, that $z$ is constant on each connected component of $G$. Choose a component $C$ on which $z$ takes a constant value $\gamma\neq0$ (it exists because $z\neq0$). For every vertex $v$ of $C$, all the neighbours of $v$ belong to $C$, so
\[
(Xz)_v=(1-\alpha)\sum_{u\sim v}z_u+\alpha d_vz_v=(1-\alpha)d_v\gamma+\alpha d_v\gamma=d_v\gamma .
\]
Since $(Xz)_v=\lambda\gamma$, we get $d_v=\lambda$ for every vertex $v$ of $C$, so $C$ is regular, a contradiction. Hence $\lambda_i(\alpha)<\lambda_i(\alpha')$.

Finally, a tree with $n\ge3$ vertices is connected and has vertices of degree $1$ and vertices of degree at least $2$, so it is not regular.
\end{proof}

\begin{remark}
The regularity hypothesis cannot be dropped: if $G$ is $r$-regular, then the all-ones vector is an eigenvector of $A_\alpha$ with eigenvalue $(1-\alpha)r+\alpha r=r$ for every $\alpha$. For instance, for $T=K_2$ the $A_\alpha$-eigenvalues are $2\alpha-1$ and $1$, and the second one is constant. We also stress that the monotonicity of $\lmid_\alpha(T)$ does \emph{not} make the difference $\alpha\mapsto\lmid_\alpha(T)-d_\alpha$ monotone, because $d_\alpha$ increases as well. This is why the set of values of $\alpha$ for which $\mathcal P_\alpha$ holds is not obviously an interval. We do not know whether $S(T)$, defined below, is always an interval.
\end{remark}

The opposite extreme of Property~\ref{property} is the value $\alpha=1$, where everything can be computed explicitly.

\begin{proposition}[The case $\alpha=1$]\label{prop:alpha1}
Let $T$ be a tree with $n\ge2$ vertices and $p$ leaves. Then
\[
\lmid_1(T)>d_1\iff p\le k_n-1,\qquad\text{equivalently}\qquad\lmid_1(T)\le d_1\iff p\ge k_n .
\]
\end{proposition}

\begin{proof}
Since $A_1=D$ is diagonal, the $A_1$-eigenvalues of $T$ are its degrees. Let $d_{(1)}\le\dots\le d_{(n)}$ be the degrees in nondecreasing order, so that $\lmid_1(T)=d_{(k_n)}$ is the $k_n$-th smallest degree. As $n\ge2$, we have $1\le d_1=2-\frac2n<2$. Since degrees are integers, $\lmid_1(T)>d_1$ holds if and only if $\lmid_1(T)\ge2$. A tree with $n\ge2$ vertices has no vertex of degree $0$, so the vertices of degree less than $2$ are exactly the $p$ leaves, which occupy the first $p$ positions of the ordered degree sequence. Therefore $d_{(k_n)}\ge2$ if and only if $k_n>p$, that is, $p\le k_n-1$. The second equivalence is the negation of the first one.
\end{proof}

We now introduce the quantities that measure how far $\mathcal P_\alpha$ survives beyond $\alpha=\frac12$.

\begin{definition}\label{def:alpha0}
For a tree $T$ with $n\ge2$ vertices let
\[
g_T(\alpha):=\lmid_\alpha(T)-d_\alpha,
\]
\[
S(T)=\Bigl\{\alpha\in\bigl[\tfrac12,1\bigr]:\ \mathcal P_\alpha(T)\text{ holds}\Bigr\}=\Bigl\{\alpha\in\bigl[\tfrac12,1\bigr]:\ g_T(\alpha)\le0\Bigr\},
\]
and $\alpha_0(T)=\max S(T)$. For $n\ge5$ we also define
\[
\alpha_{\max}(n)=\max\bigl\{\alpha_0(T)\ :\ |T|=n,\ p(T)\le k_n-1\bigr\},
\]
and, for $2\le p\le k_n-1$,
\[
\alpha_{\max}(n,p)=\max\bigl\{\alpha_0(T)\ :\ |T|=n,\ p(T)=p\bigr\}.
\]
\end{definition}

Let us check that these maxima exist. The function $g_T$ is continuous on $[\frac12,1]$, because $\alpha\mapsto A_\alpha$ is affine and eigenvalues depend continuously on the entries of a matrix (Lemma~\ref{lem:weyl}~\ref{W2}). Hence $S(T)$, the set of points of $[\frac12,1]$ where $g_T\le0$, is closed, and it is nonempty because $\frac12\in S(T)$ by Property~\ref{property}. A nonempty closed subset of the compact interval $[\frac12,1]$ has a maximum, so $\alpha_0(T)$ is well defined. By Proposition~\ref{prop:alpha1}, $\alpha_0(T)=1$ if $p\ge k_n$ and $\alpha_0(T)<1$ if $p\le k_n-1$. Every tree with $n\ge2$ vertices has $p\ge2$ leaves, and the path $P_n$ has $p=2$, so trees with $p\le k_n-1$ exist exactly when $k_n\ge3$, that is, when $n\ge5$. There are finitely many trees with $n$ vertices, so $\alpha_{\max}(n)$ and $\alpha_{\max}(n,p)$ are maxima over finite nonempty sets (for $\alpha_{\max}(n,p)$ note that, for $2\le p\le n-1$, there is a tree with $n$ vertices and exactly $p$ leaves, for instance a spider). Thus $\alpha_{\max}(n)$ is the maximum of $\alpha_0$ over the trees with $n$ vertices and at most $k_n-1$ leaves, and $\alpha_{\max}(n)=\max_{2\le p\le k_n-1}\alpha_{\max}(n,p)$. In words, $\alpha_0(T)$ is the moment after which the midrange eigenvalue of $T$ stays strictly above the average $d_\alpha$, and $\alpha_{\max}(n,p)$ is the latest such moment among the trees with $n$ vertices and $p$ leaves.

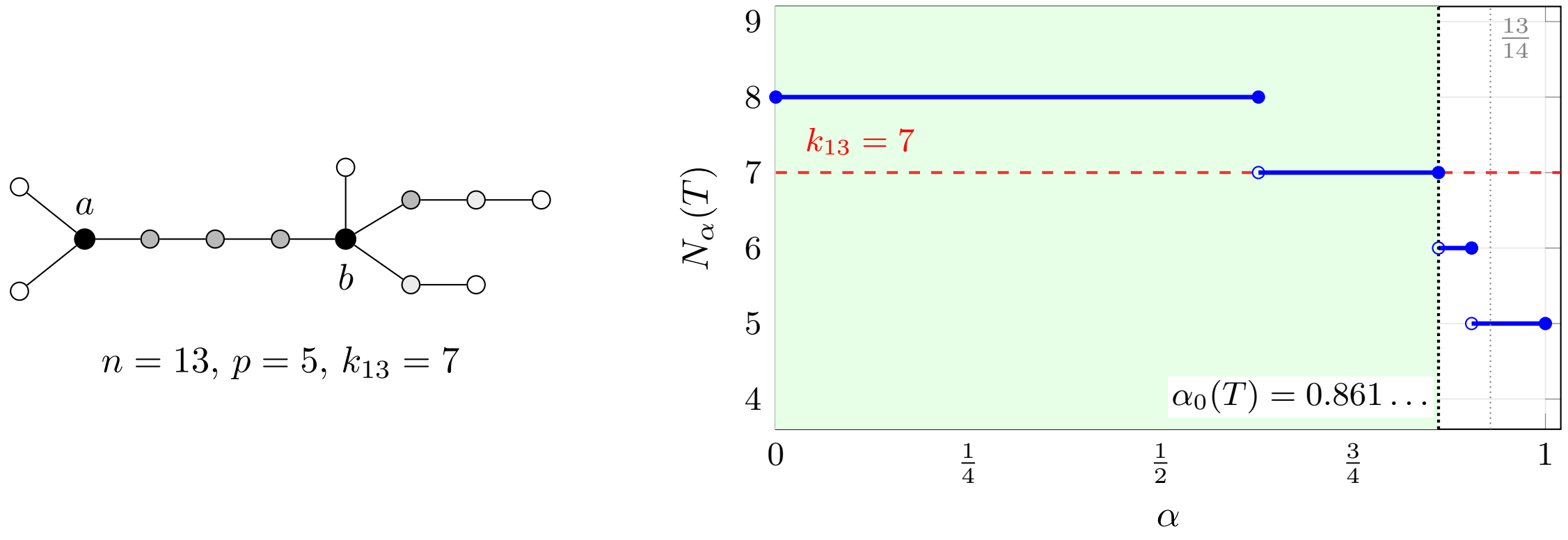
\begin{figure}[ht]
\centering
\begin{minipage}{0.40\textwidth}\centering
\begin{tikzpicture}[scale=0.62]
\node[hub,label=above:{\small $a$}] (a) at (0,0) {};
\node[leaf] (l1) at (-1,0.8) {};
\node[leaf] (l2) at (-1,-0.8) {};
\node[dtwo] (r1) at (1,0) {};
\node[dtwo] (r2) at (2,0) {};
\node[dtwo] (r3) at (3,0) {};
\node[hub,label=below:{\small $b$}] (b) at (4,0) {};
\node[leaf] (l3) at (4,1.1) {};
\node[dtwo] (s1) at (5,0.6) {};
\node[dtwob] (s2) at (6,0.6) {};
\node[leaf] (l4) at (7,0.6) {};
\node[dtwob] (t1) at (5,-0.7) {};
\node[leaf] (l5) at (6,-0.7) {};
\draw (l1)--(a)--(l2); \draw (a)--(r1)--(r2)--(r3)--(b)--(l3);
\draw (b)--(s1)--(s2)--(l4); \draw (b)--(t1)--(l5);
\end{tikzpicture}\\[2mm]
{\small $n=13$, $p=5$, $k_{13}=7$}
\end{minipage}\hfill
\begin{minipage}{0.58\textwidth}\centering
\begin{tikzpicture}
\begin{axis}[width=\textwidth,height=5.6cm,xmin=0,xmax=1.02,ymin=3.6,ymax=9.2,
  xlabel={$\alpha$},ylabel={$N_\alpha(T)$},
  xtick={0,0.25,0.5,0.75,1},xticklabels={$0$,$\frac14$,$\frac12$,$\frac34$,$1$},
  ytick={4,5,6,7,8,9},grid=major,grid style={gray!15},
  tick label style={font=\footnotesize},label style={font=\small}]
\fill[green!10] (axis cs:0,3.6) rectangle (axis cs:0.86112,9.2);
\addplot[dashed,red,thick,opacity=0.8] coordinates {(0,7) (1.02,7)};
\addplot[jump mark left,very thick,blue,mark=none] coordinates {(0,8) (0.62713,7) (0.86112,6) (0.90412,5) (1,5)};
\addplot[only marks,mark=*,mark size=1.6pt,blue] coordinates {(0,8) (0.62713,8) (0.86112,7) (0.90412,6) (1,5)};
\addplot[only marks,mark=o,mark size=1.6pt,blue,fill=white] coordinates {(0.62713,7) (0.86112,6) (0.90412,5)};
\node[red,anchor=south west,font=\footnotesize] at (axis cs:0.02,7.05) {$k_{13}=7$};
\addplot[densely dotted,thick,black] coordinates {(0.86112,3.6) (0.86112,9.2)};
\node[anchor=south east,font=\footnotesize,fill=white,inner sep=1pt] at (axis cs:0.855,3.75) {$\alpha_0(T)=0.861\ldots$};
\addplot[densely dotted,gray] coordinates {(0.92857,3.6) (0.92857,9.2)};
\node[gray,anchor=north west,font=\footnotesize,inner sep=1pt] at (axis cs:0.932,9.1) {$\frac{13}{14}$};
\end{axis}
\end{tikzpicture}
\end{minipage}
\caption{Left: the tree $T$ of Figure~\ref{fig:example}. Right: the step function $\alpha\mapsto\Nal(T)$, the number of $A_\alpha$-eigenvalues of $T$ that are less than or equal to $d_\alpha$, for $0\le\alpha\le1$. Filled dots indicate the value of $\Nal(T)$ at the jumps. The property $\mathcal P_\alpha(T)$, that is, $\Nal(T)\ge k_{13}=7$ (dashed line), holds exactly on the shaded region $[0,\alpha_0(T)]$, with $\alpha_0(T)=0.86112\ldots$, and fails beyond it; the dotted grey line marks the universal threshold $\frac{n}{n+1}=\frac{13}{14}$ of Theorem~\ref{thm:threshold}.}
\label{fig:Nalpha}
\end{figure}

Figure~\ref{fig:Nalpha} illustrates these definitions for the tree $T$ of Figure~\ref{fig:example}, which has $n=13$ vertices, $p=5$ leaves and $k_{13}=7$. In general, the function $\alpha\mapsto\Nal(T)$ is piecewise constant, and it can only change at the roots of the polynomial $\Phi_T(\alpha)=\det\bigl(A_\alpha(T)-d_\alpha I\bigr)$ of Lemma~\ref{lem:firstcrossing} below. Indeed, on an interval containing no root of $\Phi_T$, no eigenvalue of $A_\alpha(T)$ equals $d_\alpha$; since the eigenvalues $\lambda_i(\alpha)$ and $d_\alpha$ depend continuously on $\alpha$ (Lemma~\ref{lem:weyl}~\ref{W2}), each difference $\lambda_i(\alpha)-d_\alpha$ keeps its sign there, by the intermediate value theorem, and so $\Nal(T)$ is constant. For the tree of Figure~\ref{fig:example}, the roots of $\Phi_T$ in $(0,1)$ are $0.62713\ldots$, $0.86112\ldots$ and $0.90411\ldots$, and
\[
\Nal(T)=\begin{cases}
8, & 0\le\alpha\le0.62713\ldots,\\
7, & 0.62713\ldots<\alpha\le0.86112\ldots,\\
6, & 0.86112\ldots<\alpha\le0.90411\ldots,\\
5, & 0.90411\ldots<\alpha\le1 .
\end{cases}
\]
Hence $S(T)=[\frac12,\alpha_0(T)]$ with $\alpha_0(T)=0.86112\ldots$. Note that $\Nal(T)$ may exceed $k_n$: for $\alpha\le0.62713\ldots$ the property $\mathcal P_\alpha(T)$ holds with one small eigenvalue to spare. At $\alpha=1$ we recover $N_1(T)=p=5$, in agreement with the proof of Proposition~\ref{prop:alpha1}. The value $\alpha_0(T)$ lies below the universal threshold $\frac{13}{14}=0.9285\ldots$ of Theorem~\ref{thm:threshold}, and it is close to the first-order prediction $\frac{13}{15}=0.8666\ldots$ of Section~\ref{sec:firstorder} (see the caption of Figure~\ref{fig:example}).

\begin{figure}[ht]
\centering
\begin{minipage}{0.24\textwidth}\centering
\begin{tikzpicture}[scale=0.55]
\foreach \i in {0,...,6} {\ifnum\i=0 \node[leaf] (p\i) at (\i,0) {}; \else \ifnum\i=6 \node[leaf] (p\i) at (\i,0) {}; \else \node[dtwo] (p\i) at (\i,0) {}; \fi\fi}
\foreach \i [count=\j] in {0,...,5} \draw (p\i)--(p\j);
\end{tikzpicture}\\[1mm]
{\small $P_7$: $\alpha_0=0.7305$}
\end{minipage}
\begin{minipage}{0.24\textwidth}\centering
\begin{tikzpicture}[scale=0.55]
\node[hub] (u) at (0,0) {};
\node[leaf] (a) at (-0.9,0.8) {}; \node[leaf] (b) at (-0.9,-0.8) {};
\node[dtwo] (w1) at (1,0) {}; \node[dtwo] (w2) at (2,0) {}; \node[dtwo] (w3) at (3,0) {}; \node[leaf] (z) at (4,0) {};
\draw (a)--(u)--(b); \draw (u)--(w1)--(w2)--(w3)--(z);
\end{tikzpicture}\\[1mm]
{\small $S_7^\ast$: $\alpha_0=0.8300$}
\end{minipage}
\begin{minipage}{0.24\textwidth}\centering
\begin{tikzpicture}[scale=0.55]
\node[hub] (u) at (0,0) {};
\node[leaf] (a) at (0,1) {};
\node[dtwo] (b1) at (-1,-0.6) {}; \node[leaf] (b2) at (-2,-1.2) {};
\node[dtwo] (c1) at (1,-0.6) {}; \node[dtwo] (c2) at (2,-1.2) {}; \node[leaf] (c3) at (3,-1.8) {};
\draw (a)--(u)--(b1)--(b2); \draw (u)--(c1)--(c2)--(c3);
\end{tikzpicture}\\[1mm]
{\small legs $1,2,3$: $\alpha_0=0.7725$}
\end{minipage}
\begin{minipage}{0.24\textwidth}\centering
\begin{tikzpicture}[scale=0.55]
\node[hub] (u) at (0,0) {};
\foreach \ang/\nm in {90/a,210/b,330/c} {\node[dtwo] (\nm1) at (\ang:1) {}; \node[leaf] (\nm2) at (\ang:2) {}; \draw (u)--(\nm1)--(\nm2);}
\end{tikzpicture}\\[1mm]
{\small legs $2,2,2$: $\alpha_0=0.6438$}
\end{minipage}
\caption{The four trees with $7$ vertices and at most $k_7-1=3$ leaves, and their values of $\alpha_0$. The colours are as in Figure~\ref{fig:example}.}
\label{fig:n7}
\end{figure}
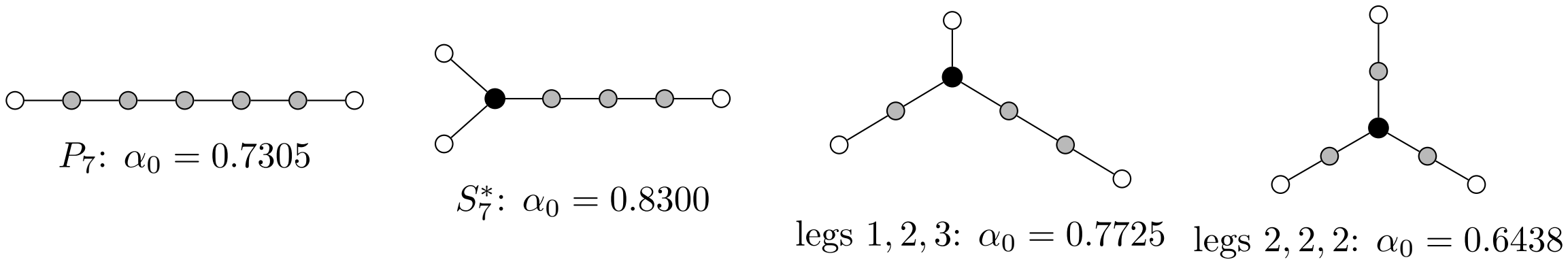

\begin{example}\label{ex:n7}
For $n=7$ there are $11$ trees, and $4$ of them have at most $k_7-1=3$ leaves: the path $P_7$ and three spiders (trees with a single hub). They are displayed in Figure~\ref{fig:n7} together with their values of $\alpha_0$. The largest one, $\alpha_{\max}(7)=0.83001\ldots$, is attained by the spider with one long leg, which is the tree $S_7^\ast$ of Section~\ref{sec:spider}. Note that the four values of $\alpha_0$ are well separated, and that the path has the second smallest one; this is consistent with the analysis of Sections~\ref{sec:firstorder} and~\ref{sec:leaves}, where long runs of degree-$2$ vertices attached to many leaves are seen to delay the crossing.
\end{example}

The next lemma explains how $\alpha_0(T)$ can be characterized algebraically.

\begin{lemma}[$\alpha_0$ is the first crossing from the right]\label{lem:firstcrossing}
Let $T$ be a tree with $n\ge3$ vertices and $p\le k_n-1$ leaves. Then $\mathcal P_\alpha(T)$ fails for every $\alpha\in(\alpha_0(T),1]$, $g_T(\alpha_0(T))=0$, and $\alpha_0(T)$ is a root of the polynomial
\[
\Phi_T(\alpha):=\det\At(T)=\det\bigl(A_\alpha(T)-d_\alpha I\bigr)\in\mathbb Q[\alpha],\qquad \deg\Phi_T\le n .
\]
In particular $\alpha_0(T)$ is a real algebraic number of degree at most $n$.
\end{lemma}

\begin{proof}
Write $\alpha_0=\alpha_0(T)$. If $\mathcal P_\alpha$ held for some $\alpha\in(\alpha_0,1]$, then $\alpha\in S(T)$ and $\alpha>\max S(T)$, which is absurd. Hence $g_T>0$ on $(\alpha_0,1]$, an interval which is nonempty because $\alpha_0<1$ (Proposition~\ref{prop:alpha1}). Since $\alpha_0\in S(T)$, we have $g_T(\alpha_0)\le0$; since $g_T$ is continuous and positive immediately to the right of $\alpha_0$, we also have $g_T(\alpha_0)=\lim_{\alpha\downarrow\alpha_0}g_T(\alpha)\ge0$. Therefore $g_T(\alpha_0)=0$, which means that $d_{\alpha_0}$ is an eigenvalue of $A_{\alpha_0}(T)$, that is, $\Phi_T(\alpha_0)=0$.

The entries of $\At$ are polynomials of degree at most $1$ in $\alpha$ with rational coefficients, see~\eqref{eq:entries}. Expanding the determinant as a sum of products of $n$ entries, $\Phi_T$ is a polynomial of degree at most $n$ with rational coefficients. It is not the zero polynomial, since at $\alpha=1$ the matrix $\At$ is diagonal and $\Phi_T(1)=\prod_v\bigl(d_v-2+\frac2n\bigr)$, where no factor vanishes because $d_v-2$ is an integer and $\frac2n$ is not an integer for $n\ge3$. Hence $\alpha_0$ is a root of a nonzero polynomial of degree at most $n$ with rational coefficients.
\end{proof}

\begin{remark}\label{rem:locate}
Lemma~\ref{lem:firstcrossing} says that $\alpha_0(T)$ is the first zero of $g_T$ met when $\alpha$ decreases from $1$. It is not necessarily the largest root of $\Phi_T$ in $[\frac12,1]$, because $\Phi_T(\alpha)=0$ whenever $d_\alpha$ is \emph{some} eigenvalue of $A_\alpha(T)$, not necessarily the midrange one. Moreover, as long as it is not known whether $S(T)$ is an interval, a sign change of $g_T$ located somewhere in $[\frac12,1]$ need not be $\alpha_0(T)$: one has to locate the first sign change from the right. Theorem~\ref{thm:threshold} below shows that one may start this search at $\alpha=\frac{n}{n+1}$, and Lemma~\ref{lem:diag} allows one to decide the validity of $\mathcal P_\alpha(T)$ at any given rational $\alpha$ in exact arithmetic.
\end{remark}

The last result of this section describes the degree structure of the trees for which $\mathcal P_\alpha$ eventually fails.

\begin{lemma}[Degree structure of trees with few leaves]\label{lem:degree}
Let $T$ be a tree with $n\ge3$ vertices and $p$ leaves. Then
\[
p=2+\sum_{v\in V_{\ge3}}(d_v-2)\ \ge\ 2+n_{\ge3}.
\]
If moreover $p\le k_n-1$, then
\[
n_{\ge3}\le p-2,\qquad n_2\ge n-2p+2,\qquad 1\le j(T)\le n_2 .
\]
In particular, more than half of the vertices have degree at most $2$, and the median degree $\lmid_1(T)$ of $T$ is exactly $2$.
\end{lemma}

\begin{proof}
Summing the degrees, $\sum_v d_v=2(n-1)$, that is, $\sum_v(d_v-2)=-2$. Splitting this sum according to the degree, the leaves contribute $-p$, the vertices of degree $2$ contribute $0$ and the hubs contribute $\sum_{v\in V_{\ge3}}(d_v-2)$; hence
\[
-p+\sum_{v\in V_{\ge3}}(d_v-2)=-2\qquad\Longrightarrow\qquad p=2+\sum_{v\in V_{\ge3}}(d_v-2)\ge2+n_{\ge3},
\]
since each hub contributes at least $1$. Consequently $n_{\ge3}\le p-2$ and $n_2=n-p-n_{\ge3}\ge n-2p+2$.
Assume now $p\le k_n-1$. Then $j(T)=k_n-p\ge1$.
For the upper bound, $j(T)\le n_2$ is equivalent to $k_n-p\le n-p-n_{\ge3}$, that is, to $k_n+n_{\ge3}\le n$. Using $n_{\ge3}\le p-2\le k_n-3$ and $2k_n\le n+1$ we get $k_n+n_{\ge3}\le 2k_n-3\le n-2<n$.
Finally $p+n_2=n-n_{\ge3}\ge n-k_n+3>k_n$, because $2k_n\le n+1<n+3$. So, in the ordered degree sequence, the first $p<k_n$ entries are equal to $1$ and the entries in positions $p+1,\dots,p+n_2$ (a range containing $k_n$) are equal to $2$. Hence the $k_n$-th smallest degree is $2$.
\end{proof}

\section{Trees with many leaves}\label{sec:exceptional}

The first main result characterizes the trees for which the conclusion of Property~\ref{property} extends to the whole interval $[0,1]$.

\begin{theorem}[Dichotomy]\label{thm:dichotomy}
Let $T$ be a tree with $n\ge3$ vertices and $p$ leaves. The following are equivalent:
\begin{enumerate}[label=(\alph*)]
\item $p\ge k_n=\lceil\frac n2\rceil$;
\item $\lmid_\alpha(T)\le d_\alpha$ for every $\alpha\in[0,1]$;
\item $\lmid_1(T)\le d_1$.
\end{enumerate}
Moreover, if $p\ge k_n$, then for every $\alpha\in[0,1]$ at least $p$ eigenvalues of $A_\alpha(T)$ are less than or equal to $\alpha$, and $\alpha\le d_\alpha$.
\end{theorem}

\begin{proof}
(a)$\Rightarrow$(b). Since $n\ge3$, no two leaves are adjacent: two adjacent leaves would form a connected component isomorphic to $K_2$, while $T$ is connected and has at least $3$ vertices. Hence the principal submatrix $A_\alpha(T)[V_1]$ has zero off-diagonal entries, and its diagonal entries are $\alpha d_v=\alpha$; that is, $A_\alpha(T)[V_1]=\alpha I_p$, whose eigenvalues are all equal to $\alpha$. By the left inequality in~\eqref{eq:cauchy},
\[
\lambda_i(A_\alpha(T))\le\lambda_i(\alpha I_p)=\alpha\qquad\text{for }1\le i\le p ,
\]
which proves the claim that at least $p$ eigenvalues are at most $\alpha$. Since $p\ge k_n$, we can take $i=k_n$ and obtain $\lmid_\alpha(T)\le\alpha$. Finally $\alpha\le\alpha\bigl(2-\frac2n\bigr)=d_\alpha$, because $2-\frac2n\ge1$ for $n\ge2$ and $\alpha\ge0$. Hence $\lmid_\alpha(T)\le d_\alpha$.

(b)$\Rightarrow$(c) is the case $\alpha=1$, and (c)$\Rightarrow$(a) is Proposition~\ref{prop:alpha1}.
\end{proof}

Theorem~\ref{thm:dichotomy} shows that trees with at least $k_n=\lceil\frac n2\rceil$ leaves are exactly the trees for which the conclusion of Property~\ref{property} extends to the whole interval $[0,1]$. We call them \emph{exceptional}. For exceptional trees $\alpha_0(T)=1$, and for all other trees $\alpha_0(T)<1$.

\begin{remark}[Exceptional trees are rare]\label{rem:rare}
Exceptional trees form an exponentially small proportion of all trees. Let $t_n$ be the number of unlabeled trees (isomorphism classes of trees) with $n$ vertices. By Otter's theorem~\cite{Otter}, $t_n=(2.9557652856\ldots)^{n+o(n)}$. On the other hand, the number of leaves of a uniformly random unlabeled tree with $n$ vertices is concentrated around $m_\ast n$, with $m_\ast=0.4381562356\ldots$~\cite{RS}, so that having at least $\lceil n/2\rceil$ leaves is a large deviation event. Its exponential rate follows from the work of Kharlamov and Orevkov~\cite{KO}, who determined the exponential growth of the number of rooted unlabeled trees half of whose vertices are leaves: the number of exceptional trees with $n$ vertices is $(2.9193800174\ldots)^{n+o(n)}$ (passing from rooted to unrooted trees changes the counts by a factor between $\frac1n$ and $1$, which does not affect the exponential rate), and therefore their proportion decays like $(0.98769\ldots)^{n+o(n)}$. Thus the phenomenon studied in this paper, the eventual failure of $\mathcal P_\alpha$, concerns almost all trees.
\end{remark}

\section{The universal threshold \texorpdfstring{$\frac{n}{n+1}$}{n/(n+1)}}\label{sec:threshold}

We now show that every non-exceptional tree loses property $\mathcal P_\alpha$ on the whole interval $[\frac n{n+1},1]$, independently of its structure.

\begin{theorem}\label{thm:threshold}
Let $T$ be a tree with $n\ge5$ vertices and $p\le k_n-1$ leaves. Then $\lmid_\alpha(T)>d_\alpha$ for every $\alpha\in\bigl[\frac{n}{n+1},1\bigr]$. Consequently
\[
\alpha_0(T)<\frac{n}{n+1}=1-\frac1{n+1}\qquad\text{and}\qquad \alpha_{\max}(n)<1-\frac{1}{n+1}.
\]
\end{theorem}

\begin{proof}
Let $W=V_2\cup V_{\ge3}$ be the set of vertices that are not leaves, so that $|W|=n-p$, and let $Y=\At[W]$.

\emph{Step 1: reduction.} We claim that it suffices to prove that $Y$ is positive definite. Indeed, if $Y\succ0$, then by~\eqref{eq:Nnpos} and~\ref{I2} of Lemma~\ref{lem:interlacing},
\[
\Nal(T)=\npos(\At)\le n-|W|=p\le k_n-1<k_n ,
\]
so $\mathcal P_\alpha$ fails (Lemma~\ref{lem:equivP}), that is, $\lmid_\alpha(T)>d_\alpha$.

\emph{Step 2: the case $\alpha=1$.} Then $\varepsilon=0$ and, by~\eqref{eq:entries}, $Y$ is diagonal with entries $d_v-2+\frac2n>0$, since every $v\in W$ has $d_v\ge2$. So $Y\succ0$.

\emph{Step 3: a sum-of-squares formula.} Let $\frac n{n+1}\le\alpha<1$, so that $\varepsilon>0$. The graph $T[W]$ is obtained from $T$ by deleting its leaves, hence it is again a tree (deleting a leaf from a tree keeps it connected and acyclic). For $v\in W$ let $d^W_v$ be the degree of $v$ in $T[W]$, so that $d^W_v\le d_v$. For $x\in\mathbb R^W$, by~\eqref{eq:entries},
\[
x^\top Yx=\sum_{v\in W}(\alpha d_v-d_\alpha)x_v^2+2\varepsilon\sum_{uv\in E(T[W])}x_ux_v .
\]
We use $2x_ux_v=(x_u+x_v)^2-x_u^2-x_v^2$ for each edge $uv$ of $T[W]$. In the resulting sum each $x_v^2$ appears once for every edge of $T[W]$ at $v$, that is, $d^W_v$ times. Therefore
\begin{equation}\label{eq:sos}
x^\top Yx=\varepsilon\sum_{uv\in E(T[W])}(x_u+x_v)^2+\sum_{v\in W}c_v\,x_v^2,\qquad c_v:=\alpha d_v-d_\alpha-\varepsilon d^W_v .
\end{equation}

\emph{Step 4: the coefficients $c_v$ are nonnegative.} Using $d^W_v\le d_v$, $d_\alpha=2\alpha-\frac{2\alpha}n$ and $\varepsilon=1-\alpha$, we obtain
\[
c_v\ge\alpha d_v-2\alpha+\frac{2\alpha}{n}-(1-\alpha)d_v=(2\alpha-1)d_v-2\alpha+\frac{2\alpha}n .
\]
If $d_v=2$, the right-hand side equals $2\alpha\bigl(1+\frac1n\bigr)-2$, which is nonnegative precisely because $\alpha\ge\frac{n}{n+1}$. If $d_v\ge3$, then, since $2\alpha-1>0$, the right-hand side is at least $3(2\alpha-1)-2\alpha=4\alpha-3>0$, because $\alpha\ge\frac n{n+1}\ge\frac56$ for $n\ge5$. Hence $c_v\ge0$ for all $v\in W$, and $x^\top Yx\ge0$ by~\eqref{eq:sos}.

\emph{Step 5: strict positivity.} We first show that some vertex $v_0\in W$ has $c_{v_0}>0$. If $T$ has a hub, any hub works by Step 4. Otherwise $T=P_n$ and $T[W]=P_{n-2}$, which has at least $3$ vertices; an end vertex $v_0$ of $T[W]$ has $d_{v_0}=2$ and $d^W_{v_0}=1$, so
\[
c_{v_0}=2\alpha-d_\alpha-\varepsilon=\frac{2\alpha}n-\varepsilon\ \ge\ \frac{2}{n+1}-\frac{1}{n+1}>0 ,
\]
where we used $\alpha\ge\frac n{n+1}$ and $\varepsilon\le\frac1{n+1}$.
Now let $x\in\mathbb R^W$ satisfy $x^\top Yx=0$. All the terms in~\eqref{eq:sos} are nonnegative, so they all vanish: $x_{v_0}=0$ and $x_u=-x_v$ for every edge $uv$ of $T[W]$. The second condition implies $|x_u|=|x_v|$ along every edge, and since $T[W]$ is connected, $|x_v|=|x_{v_0}|=0$ for every $v\in W$. Hence $x=0$, and $Y\succ0$.

\emph{Step 6: consequences.} By Steps 1--5, $\mathcal P_\alpha$ fails on $[\frac n{n+1},1]$, so $S(T)\subseteq[\frac12,\frac n{n+1})$ and $\alpha_0(T)<\frac n{n+1}$. Since there are finitely many trees with $n$ vertices, the maximum $\alpha_{\max}(n)$ is attained by some tree, and therefore it is also smaller than $\frac n{n+1}$.
\end{proof}

\begin{remark}[Dynamical meaning of $\frac n{n+1}$]\label{rem:dynamics}
Consider a run of degree-$2$ vertices processed by $\Diag$ (Lemma~\ref{lem:diag}) applied to $\At$. If $y_k$ is the value at a vertex of the run and $y_{k+1}$ the value at the next one, then~\eqref{eq:diag} gives the rational recurrence
\[
y_{k+1}=\frac{2\alpha}{n}-\frac{(1-\alpha)^2}{y_k},
\]
which is recurrence (1) of~\cite{BO}; we refer to~\cite{OT} for the systematic use of such recurrences in spectral graph theory. The Möbius map $y\mapsto c-\varepsilon^2/y$, $c=\frac{2\alpha}n$, has real fixed points (the solutions of $y^2-cy+\varepsilon^2=0$) if and only if the discriminant $c^2-4\varepsilon^2$ is nonnegative, that is, if and only if $\frac{2\alpha}{n}\ge2(1-\alpha)$, which is equivalent to
\[
\alpha\ge\frac{n}{n+1}.
\]
For $\alpha\ge\frac n{n+1}$ the map is hyperbolic or parabolic: its fixed points $y_\pm=\frac{c\pm\sqrt{c^2-4\varepsilon^2}}{2}$ are positive, $y_+$ is attracting, and every orbit started at a value larger than $y_-$ stays positive and converges monotonically to $y_+$. The values that enter a run are of this type (after a leaf one gets $c-\varepsilon^2/x_1>c\ge y_+$, and after a hub one gets $c-\varepsilon^2/a$ with $a$ a large positive value), so heuristically a run produces no nonpositive values in this regime; Theorem~\ref{thm:threshold} is the rigorous version of this observation. For $\alpha<\frac n{n+1}$ the map is elliptic: writing $\frac{c}{2\varepsilon}=\frac{\alpha}{n(1-\alpha)}=\cos\psi$ with $\psi\in(0,\frac\pi2)$, the associated linear recurrence $P_k=cP_{k-1}-\varepsilon^2P_{k-2}$ has the characteristic roots $\varepsilon e^{\pm i\psi}$, so that the orbit behaves like a rotation by the angle $\psi$, and a run of length $\ell$ produces about $(\ell+1)\psi/\pi$ sign changes. This heuristic is made precise in Lemma~\ref{lem:pendant} below. In particular, a single run of length $\lfloor n/2\rfloor$ produces its first sign change when $\psi\approx\frac{\pi}{\lfloor n/2\rfloor+1}$, that is, at $\alpha\approx\frac{n}{n+1}-O(n^{-3})$. This is the mechanism behind the sharpness of Theorem~\ref{thm:threshold}, which we prove in the next section.
\end{remark}

\section{Sharpness: the spider \texorpdfstring{$S_n^\ast$}{S\_n*}}\label{sec:spider}

In this section we show that the threshold of Theorem~\ref{thm:threshold} is asymptotically sharp. The extremal candidate is a spider with many short legs and one long leg.

\begin{definition}\label{def:spider}
For $n\ge7$ let $\ell_n=\lfloor n/2\rfloor$ and $p=k_n-1$. The spider $S_n^\ast$ consists of
\begin{itemize}
\item a centre $u$;
\item $p-1=k_n-2$ pendant vertices $v_1,\dots,v_{p-1}$ attached to $u$;
\item one pendant path $u\,w_1w_2\cdots w_{\ell_n}\,z$, with $\ell_n$ vertices $w_1,\dots,w_{\ell_n}$ of degree $2$ and a leaf $z$.
\end{itemize}
\end{definition}

\begin{figure}[ht]
\centering
\begin{tikzpicture}[scale=1.1]
\node[hub,label={[label distance=1pt]below:$u$}] (u) at (0,0) {};
\node[leaf,label=left:$v_1$] (v1) at (140:1.5) {};
\node[leaf,label=left:$v_2$] (v2) at (165:1.5) {};
\node at (192:1.35) {$\vdots$};
\node[leaf,label=left:$v_{p-1}$] (vp) at (220:1.5) {};
\draw (u)--(v1); \draw (u)--(v2); \draw (u)--(vp);
\node[dtwo,label=above:$w_1$] (w1) at (1.2,0) {};
\node[dtwo,label=above:$w_2$] (w2) at (2.4,0) {};
\node (dots) at (3.6,0) {$\cdots$};
\node[dtwo,label=above:$w_{\ell_n-1}$] (wl1) at (4.8,0) {};
\node[dtwo,label=above:$w_{\ell_n}$] (wl) at (6.0,0) {};
\node[leaf,label=above:$z$] (z) at (7.2,0) {};
\draw (u)--(w1)--(w2)--(dots)--(wl1)--(wl)--(z);
\draw[decorate,decoration={brace,amplitude=6pt,mirror}] (1.1,-0.35)--(6.1,-0.35) node[midway,below=7pt] {\small $\ell_n=\lfloor n/2\rfloor$ vertices of degree $2$};
\draw[decorate,decoration={brace,amplitude=6pt}] (-2.75,-1.1)--(-2.75,1.1) node[midway,left=8pt,align=right] {\small $p-1=$\\ \small $\lceil n/2\rceil-2$\\ \small leaves};
\end{tikzpicture}
\caption{The spider $S_n^\ast$. It has $p=\lceil n/2\rceil-1$ leaves (the $v_i$ and $z$), one hub $u$ of degree $p$, and a single run $w_1\cdots w_{\ell_n}$. For $n=7$ it is the second tree of Figure~\ref{fig:n7}.}
\label{fig:spider}
\end{figure}
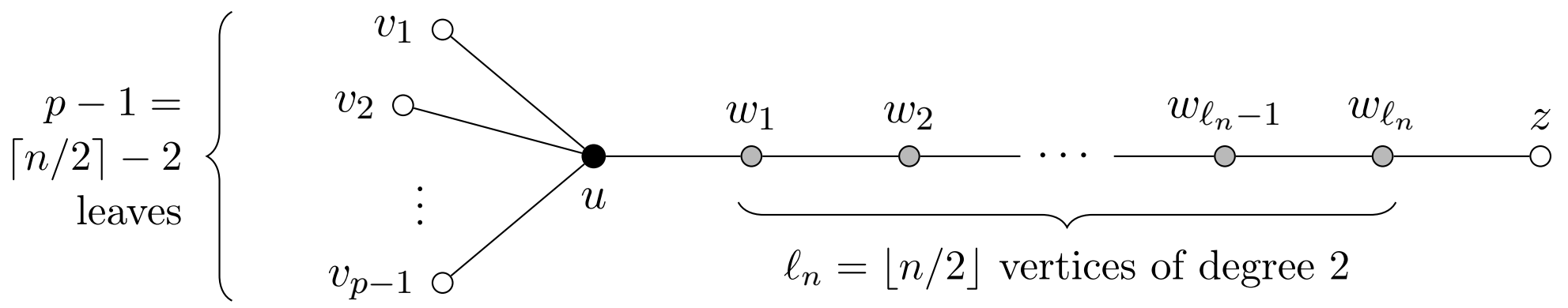

Let us check the basic parameters of $S_n^\ast$. Counting vertices, $|S_n^\ast|=1+(k_n-2)+\ell_n+1=k_n+\lfloor n/2\rfloor=n$. The leaves are $v_1,\dots,v_{p-1}$ and $z$, so $S_n^\ast$ has $p=k_n-1$ leaves, the largest number of leaves allowed for a non-exceptional tree by~\eqref{eq:halfn}, and $j(S_n^\ast)=k_n-p=1$. The centre has degree $(p-1)+1=p\ge3$, because $n\ge7$ gives $k_n\ge4$. The degree-$2$ forest $F_2$ is the single path $w_1\cdots w_{\ell_n}$, and it is as long as possible: by Lemma~\ref{lem:degree}, a tree with $p=k_n-1$ leaves and at least one hub has $n_2=n-p-n_{\ge3}\le n-k_n=\lfloor n/2\rfloor$.

\begin{theorem}\label{thm:spider}
For $n\ge7$ and $\ell_n=\lfloor n/2\rfloor$,
\[
\alpha_0(S_n^\ast)\ \ge\ 1-\frac{1}{1+n\cos(\pi/\ell_n)} .
\]
\end{theorem}

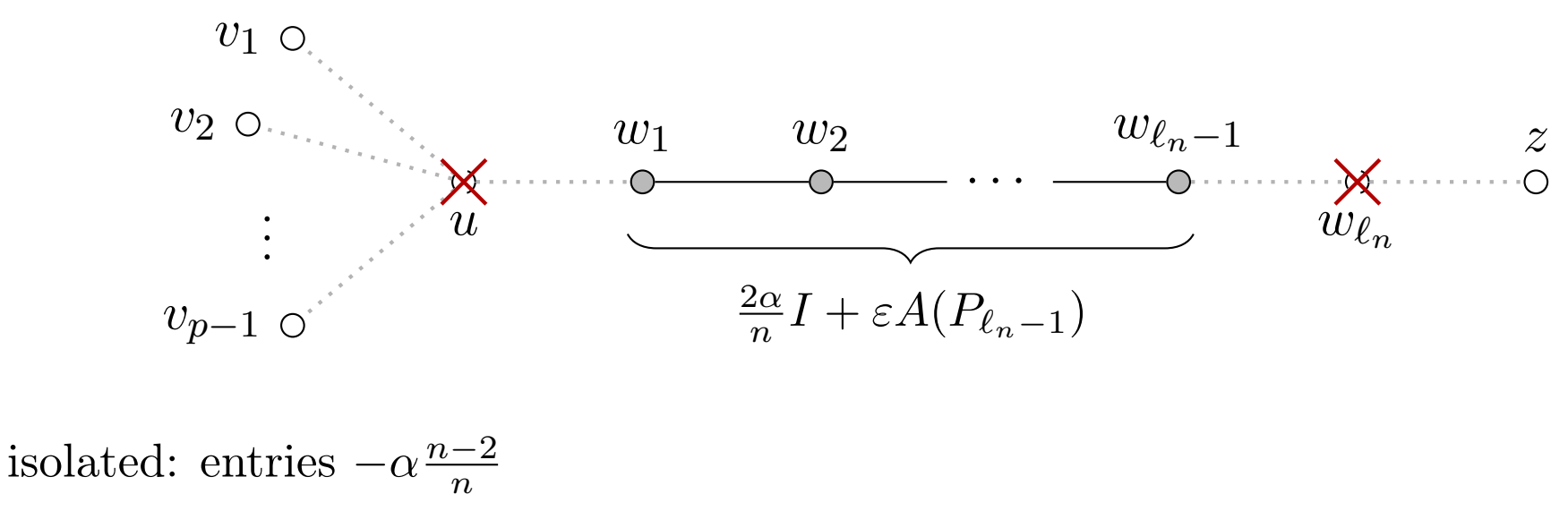
\begin{figure}[ht]
\centering
\begin{tikzpicture}[scale=1.1]
\node[gone,label=below:$u$] (u) at (0,0) {};
\node[leaf,label=left:$v_1$] (v1) at (140:1.5) {};
\node[leaf,label=left:$v_2$] (v2) at (165:1.5) {};
\node at (192:1.35) {$\vdots$};
\node[leaf,label=left:$v_{p-1}$] (vp) at (220:1.5) {};
\draw[gray!60,dotted,thick] (u)--(v1); \draw[gray!60,dotted,thick] (u)--(v2); \draw[gray!60,dotted,thick] (u)--(vp);
\node[dtwo,label=above:$w_1$] (w1) at (1.2,0) {};
\node[dtwo,label=above:$w_2$] (w2) at (2.4,0) {};
\node (dots) at (3.6,0) {$\cdots$};
\node[dtwo,label=above:$w_{\ell_n-1}$] (wl1) at (4.8,0) {};
\node[gone,label=below:$w_{\ell_n}$] (wl) at (6.0,0) {};
\node[leaf,label=above:$z$] (z) at (7.2,0) {};
\draw[gray!60,dotted,thick] (u)--(w1); \draw (w1)--(w2)--(dots)--(wl1); \draw[gray!60,dotted,thick] (wl1)--(wl)--(z);
\draw[red!70!black,thick] ($(u)+(-0.15,-0.15)$)--($(u)+(0.15,0.15)$) ($(u)+(-0.15,0.15)$)--($(u)+(0.15,-0.15)$);
\draw[red!70!black,thick] ($(wl)+(-0.15,-0.15)$)--($(wl)+(0.15,0.15)$) ($(wl)+(-0.15,0.15)$)--($(wl)+(0.15,-0.15)$);
\draw[decorate,decoration={brace,amplitude=6pt,mirror}] (1.1,-0.35)--(4.9,-0.35) node[midway,below=7pt] {\small $\frac{2\alpha}{n}I+\varepsilon A(P_{\ell_n-1})$};
\node[align=center] at (-1.4,-1.9) {\small isolated: entries $-\alpha\frac{n-2}{n}$};
\end{tikzpicture}
\caption{Proof of Theorem~\ref{thm:spider}: deleting $u$ and $w_{\ell_n}$ leaves $p$ isolated leaves, each contributing a negative eigenvalue, and the path $w_1\cdots w_{\ell_n-1}$, whose smallest eigenvalue $\frac{2\alpha}{n}-2\varepsilon\cos\frac{\pi}{\ell_n}$ becomes nonpositive as soon as $\alpha\le1-\frac{1}{1+n\cos(\pi/\ell_n)}$.}
\label{fig:spiderproof}
\end{figure}

\begin{proof}
Put $\alpha^\ast:=1-\frac1{1+n\cos(\pi/\ell_n)}=\frac{n\cos(\pi/\ell_n)}{1+n\cos(\pi/\ell_n)}$. We show that $\alpha^\ast\ge\frac12$ and that $\mathcal P_{\alpha^\ast}(S_n^\ast)$ holds; then $\alpha^\ast\in S(S_n^\ast)$ and $\alpha_0(S_n^\ast)\ge\alpha^\ast$.

\emph{Step 1: $\alpha^\ast\ge\frac12$.} Since $n\ge7$ we have $\ell_n\ge3$ and $\cos(\pi/\ell_n)\ge\cos(\pi/3)=\frac12$, so $n\cos(\pi/\ell_n)\ge\frac72>1$ and $\alpha^\ast>\frac12$.

\emph{Step 2: the condition on $\alpha$.} Let $\alpha\in(0,1)$ and $\varepsilon=1-\alpha$. The inequality
\begin{equation}\label{eq:spidercond}
\frac{2\alpha}{n}\le2\varepsilon\cos\frac{\pi}{\ell_n}
\end{equation}
is equivalent to $\alpha\le n(1-\alpha)\cos\frac{\pi}{\ell_n}$, that is, to $\alpha\bigl(1+n\cos\frac{\pi}{\ell_n}\bigr)\le n\cos\frac{\pi}{\ell_n}$, that is, to $\alpha\le\alpha^\ast$. In particular~\eqref{eq:spidercond} holds at $\alpha=\alpha^\ast$.

\emph{Step 3: a principal submatrix.} Let $U$ be the vertex set of $S_n^\ast$ with $u$ and $w_{\ell_n}$ removed (see Figure~\ref{fig:spiderproof}), and let $Y=\At[U]$, a matrix of order $n-2$. Removing $u$ and $w_{\ell_n}$ disconnects the tree into the $p$ isolated vertices $v_1,\dots,v_{p-1},z$ and the path $w_1\cdots w_{\ell_n-1}$ (here $\ell_n-1\ge2$). By~\eqref{eq:entries}, after ordering the vertices suitably,
\[
Y=\Bigl(-\alpha\,\frac{n-2}{n}\Bigr)I_p\ \oplus\ \Bigl(\frac{2\alpha}{n}I_{\ell_n-1}+\varepsilon A(P_{\ell_n-1})\Bigr).
\]
The first block has $p$ negative eigenvalues. By Lemma~\ref{lem:pathspec} with $\ell=\ell_n-1$, the eigenvalues of the second block are $\frac{2\alpha}n+2\varepsilon\cos\frac{k\pi}{\ell_n}$, $1\le k\le \ell_n-1$; the smallest one, obtained for $k=\ell_n-1$, is $\frac{2\alpha}{n}+2\varepsilon\cos\frac{(\ell_n-1)\pi}{\ell_n}=\frac{2\alpha}{n}-2\varepsilon\cos\frac{\pi}{\ell_n}$, which is nonpositive under~\eqref{eq:spidercond}.

\emph{Step 4: conclusion.} At $\alpha=\alpha^\ast$ we get $\npos(Y)\ge p+1=k_n$. By~\ref{I1} of Lemma~\ref{lem:interlacing} and~\eqref{eq:Nnpos}, $\Nal(S_n^\ast)=\npos(\At)\ge\npos(Y)\ge k_n$, which is $\mathcal P_{\alpha^\ast}(S_n^\ast)$ by Lemma~\ref{lem:equivP}.
\end{proof}

\begin{corollary}[Asymptotics of $\alpha_{\max}$]\label{cor:asymp}
For every $n\ge7$,
\[
\frac1{n+1}\ <\ 1-\alpha_{\max}(n)\ \le\ \frac{1}{1+n\cos(\pi/\lfloor n/2\rfloor)}\ =\ \frac1{n+1}+O(n^{-3}).
\]
In particular $\lim_{n\to\infty}\alpha_{\max}(n)=1$ and $n\bigl(1-\alpha_{\max}(n)\bigr)\to1$.
\end{corollary}

\begin{proof}
The left inequality is Theorem~\ref{thm:threshold}. The middle one follows from Theorem~\ref{thm:spider}, since $S_n^\ast$ has $k_n-1$ leaves and therefore $\alpha_{\max}(n)\ge\alpha_0(S_n^\ast)$. For the asymptotics, write $\ell_n=\lfloor n/2\rfloor\ge\frac{n-1}2$ and use the elementary inequality $\cos x\ge1-\frac{x^2}2$:
\[
n\cos\frac{\pi}{\ell_n}\ \ge\ n-\frac{n\pi^2}{2\ell_n^2}\ \ge\ n-\eta_n,\qquad \eta_n:=\frac{2\pi^2n}{(n-1)^2}.
\]
The sequence $\eta_n$ is decreasing, and for $n\ge7$ we have $\eta_n\le\frac{14\pi^2}{36}<4<n+1$. Hence
\[
\frac1{1+n\cos(\pi/\ell_n)}\le\frac1{n+1-\eta_n}=\frac1{n+1}+\frac{\eta_n}{(n+1)(n+1-\eta_n)} .
\]
Since $\eta_n=O(1/n)$ and $(n+1)(n+1-\eta_n)\ge(n+1)(n-3)$, the last term is $O(n^{-3})$. Finally, multiplying the two-sided estimate by $n$ we get $\frac{n}{n+1}<n(1-\alpha_{\max}(n))\le\frac{n}{n+1}+O(n^{-2})$, and both bounds tend to $1$.
\end{proof}

The bound of Theorem~\ref{thm:spider} is not strong enough to decide, for moderate values of $n$, whether $\alpha_{\max}(n)$ is larger than $1-\frac1n$. To settle this question we need a finer description of the matrices associated with pendant paths, which will also be used in Section~\ref{sec:leaves}.

\begin{lemma}[Pendant paths]\label{lem:pendant}
Let $n\ge3$, $\alpha\in(0,1)$, $\varepsilon=1-\alpha$, and let $x_1$ and $c$ be as in~\eqref{eq:x1c}. For an integer $m\ge0$ let $M_m=M_m(\alpha)$ be the symmetric tridiagonal matrix of order $m+1$ with diagonal $(x_1,c,c,\dots,c)$ and all the entries adjacent to the diagonal equal to $\varepsilon$. (If $z$ is a leaf of a tree $T$ and $z\,w_m\cdots w_1$ is a path in which $w_1,\dots,w_m$ have degree $2$, then $M_m$ is the principal submatrix of $\At(T)$ indexed by $z,w_m,\dots,w_1$, in this order, by~\eqref{eq:entries}.) Define
\[
P_{-1}=1,\qquad P_0=x_1,\qquad P_k=cP_{k-1}-\varepsilon^2P_{k-2}\quad(k\ge1).
\]
Then the following hold.
\begin{enumerate}[label=(\alph*)]
\item For $0\le k\le m$, $P_k$ is the leading principal minor of order $k+1$ of $M_m$. In particular $\det M_m=P_m$.
\item $M_m$ has at least one negative eigenvalue. If $\det M_m>0$, then $M_m$ has at least two negative eigenvalues.
\item Assume that $\zeta:=\frac{\alpha}{n\varepsilon}\in(0,1)$, write $\zeta=\cos\psi$ with $\psi\in(0,\frac\pi2)$, and let $\omega\in(0,\frac\pi2)$ be defined by $\tan\omega=\frac{\tan\psi}{n-1}$. Then $0<\omega<\psi$ and, for every $k\ge-1$,
\begin{equation}\label{eq:Pclosed}
P_k=-\rho\,\varepsilon^{k}\,\frac{\sin\bigl((k+1)\psi-\omega\bigr)}{\sin\psi},\qquad \rho:=\frac{\alpha(n-1)}{n\cos\omega}>0 .
\end{equation}
In particular, $\det M_m>0$ whenever $\pi<(m+1)\psi-\omega<2\pi$.
\item Under the assumptions of (c), if $(k+1)\psi-\omega$ is not an integer multiple of $\pi$ for $0\le k\le m$, then $M_m$ is nonsingular and it has exactly
\[
1+\Bigl\lfloor\frac{(m+1)\psi-\omega}{\pi}\Bigr\rfloor
\]
negative eigenvalues.
\end{enumerate}
\end{lemma}

\begin{proof}
(a) Let $D_k$ be the leading principal minor of order $k+1$ of $M_m$, for $0\le k\le m$, and put $D_{-1}=1$. Clearly $D_0=x_1=P_0$, and $D_1=cx_1-\varepsilon^2=cP_0-\varepsilon^2P_{-1}=P_1$. For $k\ge2$, expanding the determinant of the leading block of order $k+1$ along its last row, whose only nonzero entries are $\varepsilon$ (in column $k$) and $c$ (in column $k+1$), we obtain $D_k=cD_{k-1}-\varepsilon\cdot\varepsilon D_{k-2}$; indeed, the minor associated with the entry $\varepsilon$ has a last column whose only nonzero entry is $\varepsilon$, at the bottom. Hence $D_k$ satisfies the same recurrence and the same initial conditions as $P_k$, and $D_k=P_k$ for all $0\le k\le m$.

(b) Let $e$ be the first vector of the canonical basis. By the Rayleigh principle (the case $i=1$ of the Courant--Fischer theorem), $\lambda_1(M_m)\le e^\top M_me=x_1<0$. If $\det M_m>0$, then the product of the eigenvalues of $M_m$ is positive, so no eigenvalue vanishes and the number of negative eigenvalues is even. Since it is at least one, it is at least two.

(c) Since $\tan$ is increasing on $(0,\frac\pi2)$ and $\tan\omega=\frac{\tan\psi}{n-1}<\tan\psi$, we have $0<\omega<\psi$. Let $S_k$ denote the right-hand side of~\eqref{eq:Pclosed} and put $\theta_k=(k+1)\psi-\omega$. We check that $S_k$ satisfies the recurrence and the initial conditions of $P_k$. From the identity $\sin\theta_k+\sin\theta_{k-2}=2\sin\theta_{k-1}\cos\psi$ we get
\[
S_k=-\rho\varepsilon^k\frac{2\cos\psi\sin\theta_{k-1}-\sin\theta_{k-2}}{\sin\psi}=2\varepsilon\cos\psi\,S_{k-1}-\varepsilon^2S_{k-2}=cS_{k-1}-\varepsilon^2S_{k-2},
\]
because $2\varepsilon\cos\psi=2\varepsilon\zeta=\frac{2\alpha}{n}=c$. For the initial values we use $\rho\cos\omega=\frac{\alpha(n-1)}{n}$ and $\rho\sin\omega=\rho\cos\omega\tan\omega=\frac{\alpha(n-1)}{n}\cdot\frac{\tan\psi}{n-1}=\frac{\alpha\tan\psi}{n}$. Then
\[
S_{-1}=\frac{\rho\sin\omega}{\varepsilon\sin\psi}=\frac{\alpha}{n\varepsilon\cos\psi}=\frac{\zeta}{\zeta}=1,
\]
and
\begin{align*}
S_0&=-\rho\,\frac{\sin\psi\cos\omega-\cos\psi\sin\omega}{\sin\psi}=-\rho\cos\omega+\rho\sin\omega\,\frac{\cos\psi}{\sin\psi}\\
&=-\frac{\alpha(n-1)}{n}+\frac{\alpha}{n}=-\frac{\alpha(n-2)}{n}=x_1 .
\end{align*}
By induction, $S_k=P_k$ for every $k\ge-1$. Finally, $\rho$, $\varepsilon^m$ and $\sin\psi$ are positive, so $\det M_m=P_m>0$ if and only if $\sin\bigl((m+1)\psi-\omega\bigr)<0$, which is the case when $\pi<(m+1)\psi-\omega<2\pi$.

(d) Root the path $z\,w_m\cdots w_1$ underlying $M_m$ at $w_1$, so that the bottom-up order is $z,w_m,\dots,w_1$, and apply $\Diag(M_m,0)$. The value at $z$ is $x_1=P_0/P_{-1}$, and if the value at the $k$-th processed vertex is $P_{k-1}/P_{k-2}\ne0$, then by~\eqref{eq:diag} the value at the next one is
\[
c-\varepsilon^2\frac{P_{k-2}}{P_{k-1}}=\frac{cP_{k-1}-\varepsilon^2P_{k-2}}{P_{k-1}}=\frac{P_k}{P_{k-1}} .
\]
By (c) and the hypothesis, $P_k\neq0$ for $0\le k\le m$, so the third branch of the algorithm is never used and the values are $P_k/P_{k-1}$, $0\le k\le m$, all nonzero. By Lemma~\ref{lem:diag}, $M_m$ is nonsingular and its number of negative eigenvalues is the number of indices $k\in\{0,\dots,m\}$ such that $P_k$ and $P_{k-1}$ have opposite signs. By~\eqref{eq:Pclosed}, the sign of $P_k$ is the opposite of the sign of $\sin y_{k+1}$, where $y_i:=i\psi-\omega$ for $0\le i\le m+1$. Hence we must count the sign changes in the sequence $\sin y_0,\sin y_1,\dots,\sin y_{m+1}$, none of whose terms vanishes. For $y\notin\pi\mathbb Z$ the sign of $\sin y$ is $(-1)^{\lfloor y/\pi\rfloor}$. Since $y_i-y_{i-1}=\psi<\pi$, the integers $\lfloor y_i/\pi\rfloor$ form a nondecreasing sequence whose consecutive differences are $0$ or $1$, and a sign change occurs exactly when the difference is $1$. Therefore the number of sign changes is
\[
\Bigl\lfloor\frac{y_{m+1}}{\pi}\Bigr\rfloor-\Bigl\lfloor\frac{y_0}{\pi}\Bigr\rfloor=\Bigl\lfloor\frac{(m+1)\psi-\omega}{\pi}\Bigr\rfloor+1,
\]
because $y_0=-\omega\in(-\frac\pi2,0)$, so that $\lfloor y_0/\pi\rfloor=-1$.
\end{proof}

Lemma~\ref{lem:pendant}~(d) is the precise form of the heuristic of Remark~\ref{rem:dynamics}: in the elliptic regime, a pendant path with $m$ vertices of degree $2$ ending at a leaf behaves like a rotation by the angle $\psi$ started at the phase $-\omega$, and every passage of the phase through a multiple of $\pi$ produces a negative eigenvalue. We now use it to decide the comparison between $\alpha_{\max}(n)$ and $1-\frac1n$.

\begin{corollary}\label{cor:false}
$\alpha_{\max}(n)>1-\frac1n$ for $n=16$ and for every $n\ge18$.
\end{corollary}

\begin{proof}
Since $\alpha_{\max}(n)\ge\alpha_0(S_n^\ast)$, it suffices to prove that $\alpha_0(S_n^\ast)>1-\frac1n$. Let $\ell_n=\lfloor n/2\rfloor$. We consider two ranges of $n$.

\emph{Case 1: $n\ge22$.} By Theorem~\ref{thm:spider}, it suffices that $1+n\cos(\pi/\ell_n)>n$, that is, $\cos(\pi/\ell_n)>1-\frac1n$; indeed, then $\frac{1}{1+n\cos(\pi/\ell_n)}<\frac1n$. By the estimate in the proof of Corollary~\ref{cor:asymp}, $\cos(\pi/\ell_n)\ge1-\frac{2\pi^2}{(n-1)^2}$, and $\frac{2\pi^2}{(n-1)^2}<\frac1n$ is equivalent to $(n-1)^2>2\pi^2n$. This inequality holds for $n=22$ ($441>434.3$), and therefore for all $n\ge22$, because $(n-1)^2/n=n-2+\frac1n$ is increasing for $n\ge1$.

\emph{Case 2: $n=16$ or $18\le n\le35$.} Let $\alpha=1-\frac1n$, so that $\varepsilon=\frac1n$ and $\zeta=\frac{\alpha}{n\varepsilon}=1-\frac1n\in(0,1)$. Let $\psi_n=\arccos(1-\frac1n)$ and let $\omega_n$ be as in Lemma~\ref{lem:pendant}~(c). Since $\sin\psi_n=\frac{\sqrt{2n-1}}{n}$, we have $\tan\psi_n=\frac{\sqrt{2n-1}}{n-1}$, hence
\[
\omega_n\le\tan\omega_n=\frac{\sqrt{2n-1}}{(n-1)^2}.
\]
We first bound $\psi_n$. From $\cos x\ge1-\frac{x^2}{2}$ we get $\cos\sqrt{2/n}\ge1-\frac1n=\cos\psi_n$, and since $\cos$ is decreasing on $[0,\pi]$, $\psi_n\ge\sqrt{2/n}$. From $\cos x\le1-\frac{x^2}2+\frac{x^4}{24}$ we get $\cos\frac{2}{\sqrt n}\le1-\frac2n+\frac{2}{3n^2}<1-\frac1n=\cos\psi_n$, hence $\psi_n<\frac{2}{\sqrt n}$.

We claim that $\pi<(\ell_n+1)\psi_n-\omega_n<2\pi$. For the upper bound, $(\ell_n+1)\psi_n-\omega_n<(\ell_n+1)\frac{2}{\sqrt n}\le\frac{n+2}{2}\cdot\frac{2}{\sqrt n}=\sqrt n+\frac{2}{\sqrt n}$, which is increasing in $n$ and equal to $\sqrt{35}+\frac{2}{\sqrt{35}}<5.917+0.339<6.26<2\pi$ for $n=35$. For the lower bound, if $n=16$ then $\ell_n+1=9$ and
\[
(\ell_n+1)\psi_n-\omega_n\ge\frac{9}{2\sqrt2}-\frac{\sqrt{31}}{225}>3.1819-0.0248>3.157>\pi .
\]
If $n\ge18$, then $\ell_n+1\ge\frac{n+1}{2}$, so $(\ell_n+1)\psi_n\ge\frac{n+1}{\sqrt{2n}}$, which is increasing in $n$ and at least $\frac{19}{6}>3.1666$; on the other hand $\frac{\sqrt{2n-1}}{(n-1)^2}$ is decreasing in $n$ and at most $\frac{\sqrt{35}}{289}<0.0205$. Hence $(\ell_n+1)\psi_n-\omega_n>3.1666-0.0205>3.146>\pi$.

By Lemma~\ref{lem:pendant}~(c) with $m=\ell_n$, the determinant of $M_{\ell_n}(\alpha)$ is positive at $\alpha=1-\frac1n$. Since $\det M_{\ell_n}(\alpha)$ is a polynomial in $\alpha$ (Lemma~\ref{lem:pendant}~(a)), it remains positive on an interval $[1-\frac1n,1-\frac1n+\delta)$ for some $\delta>0$. Fix $\alpha$ in this interval with $\alpha>1-\frac1n$ and $\alpha<1$. Let $U$ be the vertex set of $S_n^\ast$ with the centre $u$ removed. Removing $u$ leaves the $p-1$ isolated leaves $v_i$ and the path $z\,w_{\ell_n}\cdots w_1$, so
\[
\At(S_n^\ast)[U]=x_1I_{p-1}\oplus M_{\ell_n}(\alpha).
\]
The first block has $p-1$ negative eigenvalues, and by Lemma~\ref{lem:pendant}~(b) the second block has at least two. Hence $\npos(\At[U])\ge p+1=k_n$, and by~\ref{I1} of Lemma~\ref{lem:interlacing}, $\mathcal P_\alpha(S_n^\ast)$ holds. Therefore $\alpha_0(S_n^\ast)\ge\alpha>1-\frac1n$.

The two cases cover $n=16$ and all $n\ge18$.
\end{proof}

\begin{remark}\label{rem:n17}
Case 2 of the proof requires $(\ell_n+1)\psi_n>\pi+\omega_n$. For $n=17$ we have $\ell_n+1=9$ and $\cos\frac{\pi}{9}=0.9396\ldots<\frac{16}{17}=0.9411\ldots$, so that $\psi_{17}<\frac\pi9$ and $9\psi_{17}<\pi$: at $\alpha=1-\frac1{17}$ the pendant path of $S_{17}^\ast$ has not yet produced its second negative eigenvalue, and the argument does not apply. This parity effect, due to $\ell_n=\lfloor n/2\rfloor$, is the reason why $n=17$ is not covered by Corollary~\ref{cor:false}.
\end{remark}

\begin{remark}[Small values of $n$]\label{rem:table}
For $5\le n\le17$ the value of $\alpha_{\max}(n)$ can be determined exactly by examining all the non-isomorphic trees with $n$ vertices and at most $k_n-1$ leaves. For each such tree, $\alpha_0(T)$ is located as described in Remark~\ref{rem:locate}: one decreases $\alpha$ from $\frac{n}{n+1}$ along a fine grid until $\mathcal P_\alpha(T)$ holds for the first time, refines by bisection inside this first bracket only, and decides the validity of $\mathcal P_\alpha(T)$ at each rational value of $\alpha$ by counting the nonpositive values produced by $\Diag$ (Lemma~\ref{lem:diag}) in exact arithmetic. The results are collected in Table~\ref{tab:alphamax}. For every $7\le n\le17$ the maximum is attained by $S_n^\ast$, and for $n\le14$ it is attained only by $S_n^\ast$; for instance, the second largest value of $\alpha_0$ for $n=14$ is $0.92673923$, against $\alpha_0(S_{14}^\ast)=0.92814853$. For $n=5,6$ we have $k_n=3$, and the only tree with $p\le k_n-1=2$ leaves is the path. The table also shows that $\alpha_{\max}(n)<1-\frac1n$ for $n\le15$ and for $n=17$, while $\alpha_{\max}(16)>1-\frac1{16}$, in agreement with Corollary~\ref{cor:false} and Remark~\ref{rem:n17}. Finally, for every $n$ in the table, $\alpha_0(S_n^\ast)$ coincides, up to $10^{-10}$, with the largest root of the function $R_n$ of Proposition~\ref{prop:spidereq} in $(\frac12,\frac{n}{n+1})$.
\end{remark}

\begin{table}[ht]
\centering
\begin{tabular}{rrcccl}
\toprule
$n$ & \#$\{T: p\le k_n-1\}$ & $\alpha_{\max}(n)$ & $1-\frac1n$ & $\frac n{n+1}$ & maximiser\\
\midrule
5 & 1 & 0.74271573 & 0.800000 & 0.833333 & $P_5$\\
6 & 1 & 0.82177615 & 0.833333 & 0.857143 & $P_6$\\
7 & 4 & 0.83001656 & 0.857143 & 0.875000 & $S_7^\ast$\\
8 & 5 & 0.86570271 & 0.875000 & 0.888889 & $S_8^\ast$\\
9 & 20 & 0.87838476 & 0.888889 & 0.900000 & $S_9^\ast$\\
10 & 31 & 0.89614291 & 0.900000 & 0.909091 & $S_{10}^\ast$\\
11 & 109 & 0.90465183 & 0.909091 & 0.916667 & $S_{11}^\ast$\\
12 & 188 & 0.91516988 & 0.916667 & 0.923077 & $S_{12}^\ast$\\
13 & 660 & 0.92117986 & 0.923077 & 0.928571 & $S_{13}^\ast$\\
14 & 1229 & 0.92814853 & 0.928571 & 0.933333 & $S_{14}^\ast$\\
15 & 4183 & 0.93261746 & 0.933333 & 0.937500 & $S_{15}^\ast$\\
16 & 8289 & \textbf{0.93758943} & 0.937500 & 0.941176 & $S_{16}^\ast$\\
17 & 27752 & 0.94104549 & 0.941176 & 0.944444 & $S_{17}^\ast$\\
\bottomrule
\end{tabular}
\caption{The value of $\alpha_{\max}(n)$ over all the non-isomorphic trees with $n$ vertices and at most $k_n-1$ leaves, together with the number of such trees and the tree attaining the maximum. In bold, the first value exceeding $1-\frac1n$ (compare with Corollary~\ref{cor:false}).}
\label{tab:alphamax}
\end{table}

The exact value of $\alpha_0(S_n^\ast)$ is a root of an explicit function, which can be written with Chebyshev polynomials.

\begin{proposition}[Characteristic equation of $S_n^\ast$]\label{prop:spidereq}
Let $n\ge7$, $p=k_n-1$, $\ell_n=\lfloor n/2\rfloor$, $\varepsilon=1-\alpha$, $\zeta=\frac{\alpha}{n\varepsilon}$, $x_1=-\alpha\frac{n-2}{n}$, let $U_k$ be the Chebyshev polynomials of the second kind, defined by $U_0=1$, $U_1(\zeta)=2\zeta$ and $U_k=2\zeta U_{k-1}-U_{k-2}$, extended by $U_{-1}=0$ and $U_{-2}=-1$, and put
\[
Q_k(\alpha)=x_1U_k(\zeta)-\varepsilon U_{k-1}(\zeta),\qquad
h(\alpha)=\alpha\Bigl(p-2+\frac2n\Bigr)+\frac{(p-1)\,n\,\varepsilon^2}{\alpha(n-2)} .
\]
Then, for $\alpha\in(0,1)$,
\[
\det\At(S_n^\ast)=x_1^{\,p-1}\,\varepsilon^{\ell_n}\,R_n(\alpha),\qquad R_n(\alpha):=h(\alpha)Q_{\ell_n}(\alpha)-\varepsilon\,Q_{\ell_n-1}(\alpha),
\]
and $\alpha_0(S_n^\ast)$ is a root of $R_n$ in $(\frac12,\frac n{n+1})$.
\end{proposition}

\begin{figure}[ht]
\centering
\begin{tikzpicture}[scale=1.15]
\node[hub,label={[label distance=1pt]above:$u$}] (u) at (0,0) {};
\node (ulab) at (0.9,-1.25) {\scriptsize $a(u)=h-\varepsilon^2\frac{P_{\ell_n-1}}{P_{\ell_n}}$};
\draw[gray,->,shorten >=2pt] (ulab.north west)--(u);
\node[leaf,label=left:{\scriptsize $x_1$}] (v1) at (140:1.5) {};
\node[leaf,label=left:{\scriptsize $x_1$}] (v2) at (165:1.5) {};
\node at (192:1.35) {$\vdots$};
\node[leaf,label=left:{\scriptsize $x_1$}] (vp) at (220:1.5) {};
\draw (u)--(v1); \draw (u)--(v2); \draw (u)--(vp);
\node[dtwo,label=above:{\scriptsize $\frac{P_{\ell_n}}{P_{\ell_n-1}}$}] (w1) at (1.3,0) {};
\node[dtwo,label=above:{\scriptsize $\frac{P_{\ell_n-1}}{P_{\ell_n-2}}$}] (w2) at (2.6,0) {};
\node (dots) at (3.8,0) {$\cdots$};
\node[dtwo,label=above:{\scriptsize $\frac{P_2}{P_1}$}] (wl1) at (5.0,0) {};
\node[dtwo,label=above:{\scriptsize $\frac{P_1}{P_0}$}] (wl) at (6.2,0) {};
\node[leaf,label=above:{\scriptsize $x_1=\frac{P_0}{P_{-1}}$}] (z) at (7.4,0) {};
\draw (u)--(w1)--(w2)--(dots)--(wl1)--(wl)--(z);
\node[below] at (1.3,-0.1) {\scriptsize $w_1$}; \node[below] at (2.6,-0.1) {\scriptsize $w_2$};
\node[below] at (5.0,-0.1) {\scriptsize $w_{\ell_n-1}$}; \node[below] at (6.2,-0.1) {\scriptsize $w_{\ell_n}$}; \node[below] at (7.4,-0.1) {\scriptsize $z$};
\draw[->,gray] (7.4,-0.75)--(1.3,-0.75) node[midway,below] {\scriptsize processing order};
\end{tikzpicture}
\caption{The values of $\Diag$ applied to $\At(S_n^\ast)$ and rooted at $u$, where $P_{-1}=1$, $P_0=x_1$ and $P_k=\frac{2\alpha}{n}P_{k-1}-\varepsilon^2P_{k-2}$ (compare with Lemma~\ref{lem:pendant}).}
\label{fig:spiderdiag}
\end{figure}
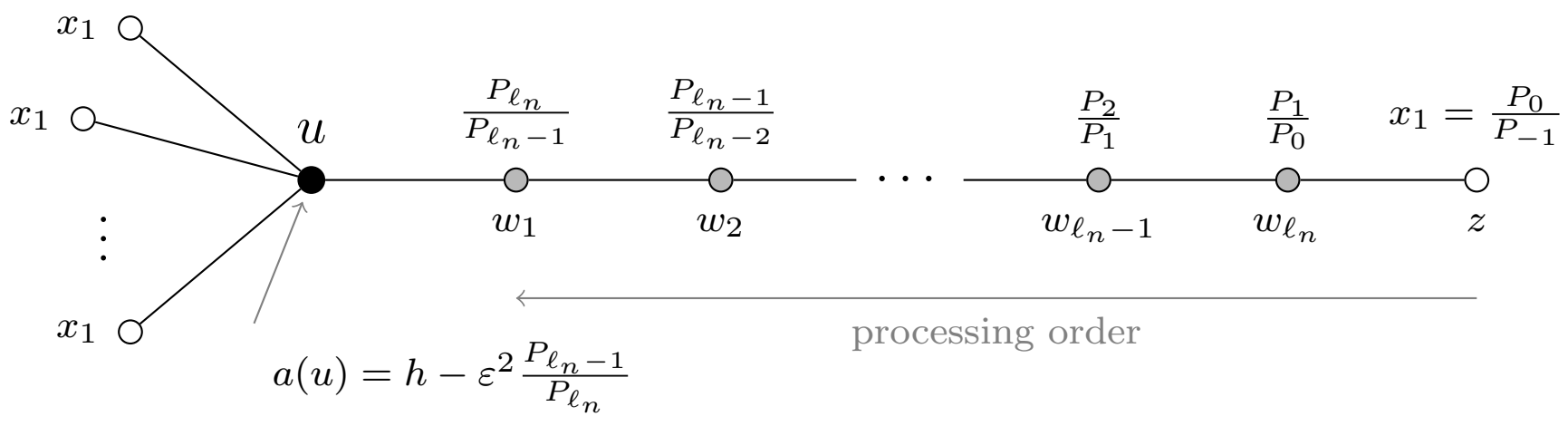

\begin{proof}
\emph{Step 1: the values of $\Diag$.} Root $S_n^\ast$ at $u$ and apply the algorithm $\Diag$ to $A_\alpha(S_n^\ast)$ with $x=-d_\alpha$. First suppose that all the values below are nonzero, so that Lemma~\ref{lem:diag} applies in the form~\eqref{eq:diag}. Each pendant vertex $v_i$ has no children, so $a(v_i)=(\At)_{v_iv_i}=\alpha-d_\alpha=x_1$; the same holds for $z$. Let $P_k$ be as in Lemma~\ref{lem:pendant}. We claim that $a(w_{\ell_n-k+1})=P_k/P_{k-1}$ for $k=1,\dots,\ell_n$. For $k=1$, the vertex $w_{\ell_n}$ has the single child $z$, so
\[
a(w_{\ell_n})=\frac{2\alpha}n-\frac{\varepsilon^2}{x_1}=\frac{(2\alpha/n)P_0-\varepsilon^2P_{-1}}{P_0}=\frac{P_1}{P_0}.
\]
If $a(w_{\ell_n-k+2})=P_{k-1}/P_{k-2}$, then, since $w_{\ell_n-k+2}$ is the only child of $w_{\ell_n-k+1}$,
\[
a(w_{\ell_n-k+1})=\frac{2\alpha}n-\varepsilon^2\frac{P_{k-2}}{P_{k-1}}=\frac{(2\alpha/n)P_{k-1}-\varepsilon^2P_{k-2}}{P_{k-1}}=\frac{P_k}{P_{k-1}} .
\]
Finally, $u$ has the $p-1$ children $v_i$ and the child $w_1$, and $(\At)_{uu}=\alpha p-d_\alpha$, so
\[
a(u)=\alpha p-d_\alpha-(p-1)\frac{\varepsilon^2}{x_1}-\varepsilon^2\frac{P_{\ell_n-1}}{P_{\ell_n}}=h(\alpha)-\varepsilon^2\frac{P_{\ell_n-1}}{P_{\ell_n}},
\]
because $\alpha p-d_\alpha=\alpha\bigl(p-2+\frac2n\bigr)$ and $-\frac{\varepsilon^2}{x_1}=\frac{n\varepsilon^2}{\alpha(n-2)}$.

\emph{Step 2: the determinant.} By Lemma~\ref{lem:diag}, $\det\At$ is the product of all the values. The product along $z,w_{\ell_n},\dots,w_1$ telescopes: $\frac{P_0}{P_{-1}}\cdot\frac{P_1}{P_0}\cdots\frac{P_{\ell_n}}{P_{\ell_n-1}}=P_{\ell_n}$. Therefore
\[
\det\At=x_1^{p-1}\,P_{\ell_n}\Bigl(h-\varepsilon^2\frac{P_{\ell_n-1}}{P_{\ell_n}}\Bigr)=x_1^{p-1}\bigl(hP_{\ell_n}-\varepsilon^2P_{\ell_n-1}\bigr).
\]

\emph{Step 3: Chebyshev form.} Put $Q_k=P_k/\varepsilon^k$ for $k\ge-1$. Dividing the recurrence for $P_k$ by $\varepsilon^k$ gives $Q_k=\frac{2\alpha}{n\varepsilon}Q_{k-1}-Q_{k-2}=2\zeta Q_{k-1}-Q_{k-2}$, with $Q_{-1}=\varepsilon$ and $Q_0=x_1$. The sequence $x_1U_k(\zeta)-\varepsilon U_{k-1}(\zeta)$ satisfies the same recurrence (because $U_k$ does, also for $k=0$ with the conventions $U_{-1}=0$ and $U_{-2}=-1$) and has the same initial values: at $k=-1$ it equals $x_1\cdot0-\varepsilon\cdot(-1)=\varepsilon$, and at $k=0$ it equals $x_1\cdot1-\varepsilon\cdot0=x_1$. Hence $Q_k=x_1U_k(\zeta)-\varepsilon U_{k-1}(\zeta)$ for all $k\ge-1$, and
\[
hP_{\ell_n}-\varepsilon^2P_{\ell_n-1}=\varepsilon^{\ell_n}\bigl(hQ_{\ell_n}-\varepsilon Q_{\ell_n-1}\bigr)=\varepsilon^{\ell_n}R_n(\alpha).
\]

\emph{Step 4: removing the nonvanishing assumption.} Both $\det\At$ and $x_1^{p-1}\varepsilon^{\ell_n}R_n(\alpha)$ are rational functions of $\alpha$ on $(0,1)$ (the second one has the factor $\frac1\alpha$ in $h$), and they coincide at every $\alpha\in(0,1)$ for which no $P_k$ with $0\le k\le \ell_n$ vanishes. Each $P_k$ is a polynomial in $\alpha$ (because $x_1$, $\frac{2\alpha}n$ and $\varepsilon^2$ are), and it is not the zero polynomial: at $\alpha=1$ we have $\varepsilon=0$ and $P_k(1)=(\frac2n)^kx_1(1)=-(\frac2n)^k\frac{n-2}{n}\ne0$. Hence each $P_k$ vanishes only at finitely many points. Two rational functions that agree on an infinite set are equal, so the identity holds for all $\alpha\in(0,1)$. The statement about $\alpha_0(S_n^\ast)$ follows from Lemma~\ref{lem:firstcrossing}, because $x_1\neq0$ and $\varepsilon\neq0$ for $\alpha\in(\frac12,1)$, because $\alpha_0(S_n^\ast)>\frac12$ by Theorem~\ref{thm:spider}, and because $\alpha_0(S_n^\ast)<\frac n{n+1}$ by Theorem~\ref{thm:threshold}.
\end{proof}

Since $U_k(\cos\psi)=\frac{\sin((k+1)\psi)}{\sin\psi}$, in the elliptic regime $\zeta=\cos\psi\in(0,1)$ the functions $Q_k$ are trigonometric, and the closed form~\eqref{eq:Pclosed} is precisely $Q_k=x_1U_k(\zeta)-\varepsilon U_{k-1}(\zeta)$ rewritten as a single sine. In this form, $R_n=0$ is the condition that the value at the hub, $a(u)=h-\varepsilon Q_{\ell_n-1}/Q_{\ell_n}$, vanishes; since $h(\alpha)\approx p-2$ is large, this forces $Q_{\ell_n}$ to be close to zero, that is, $(\ell_n+1)\psi-\omega$ close to $\pi$. This observation suggests the second-order asymptotics $\frac{1}{1-\alpha_0(S_n^\ast)}=n+1-\frac{2\pi^2}{n}+O(n^{-2})$, whose inequality ``$\ge$'' follows from Theorem~\ref{thm:spider} (see the proof of Corollary~\ref{cor:asymp}).

\section{First-order theory: the degree-two forest}\label{sec:firstorder}

We now explain, for a general tree, which structural features decide the position of the crossing near $\alpha=1$. Order the vertices as $V_1$, $V_2$, $V_{\ge3}$. By~\eqref{eq:entries}, $\At$ has the block form
\begin{equation}\label{eq:blocks}
\At=\begin{pmatrix}
-\alpha\frac{n-2}{n}I_p & \varepsilon E_{12} & \varepsilon E_{13}\\[1mm]
\varepsilon E_{12}^\top & \frac{2\alpha}{n}I+\varepsilon A(F_2) & \varepsilon E_{23}\\[1mm]
\varepsilon E_{13}^\top & \varepsilon E_{23}^\top & \operatorname{diag}\bigl(\alpha(d_v-2+\frac2n)\bigr)+\varepsilon A(T[V_{\ge3}])
\end{pmatrix},
\end{equation}
where $E_{ab}$ is the $0/1$ matrix recording the edges between the classes $V_a$ and $V_b$ (with $V_3:=V_{\ge3}$). The first diagonal block is diagonal because no two leaves are adjacent (for $n\ge3$). When $\varepsilon=O(1/n)$, the first diagonal block is close to $-I$, the third one is, roughly speaking, at least $I$, and the middle one is of order $1/n$. The couplings between the blocks are of order $\varepsilon$, so eliminating the first and the third blocks (by Schur complements) perturbs the middle block only by terms of order $\varepsilon^2$. The first block contributes $p$ negative eigenvalues and the third one $n_{\ge3}$ positive ones. Therefore $\mathcal P_\alpha$ requires $j(T)=k_n-p$ nonpositive eigenvalues coming from the middle block, which, to first order, happens if and only if
\[
\frac{2\alpha}{n}+\varepsilon\,\theta_{j}(F_2)\le0 .
\]
The following theorem turns this heuristic into rigorous two-sided estimates.

\begin{theorem}[Two-sided first-order estimate]\label{thm:firstorder}
Let $T$ be a tree with $n\ge7$ vertices, $p\le k_n-1$ leaves and $j=j(T)$.
\begin{enumerate}[label=(\alph*)]
\item \emph{(Lower bound.)} If $j\le\mu(\Fo)$, then $\theta_j(\Fo)<0$ and
\[
\alpha_0(T)\ \ge\ \max\Bigl\{\tfrac12,\ 1-\frac{2}{2-n\,\theta_j(\Fo)}\Bigr\}.
\]
\item \emph{(Failure criterion.)} Let $0<\varepsilon\le\frac17$. If
\[
\theta_j(F_2)\ >\ \frac{12\,\varepsilon}{1-\varepsilon}-\frac{2(1-\varepsilon)}{n\,\varepsilon},
\]
then $\mathcal P_{1-\varepsilon}(T)$ fails.
\item \emph{(Slow regime.)} If $j>\mu(F_2)$, equivalently $\theta_j(F_2)\ge0$, then
\[
\alpha_0(T)\ \le\ 1-\frac{1}{1+\sqrt{6n}} .
\]
\end{enumerate}
\end{theorem}

\begin{proof}
(a) \emph{Step 1: the submatrix.} Let $U=V_1\cup V_2^\circ$ and $Y=\At[U]$. A neighbour of a leaf is either a hub or a vertex of degree $2$ adjacent to a leaf (it cannot be a leaf, since $n\ge3$); neither of them belongs to $U$. Hence no vertex of $V_1$ is adjacent to another vertex of $U$, and by~\eqref{eq:entries}
\[
Y=\Bigl(-\alpha\frac{n-2}n\Bigr)I_p\ \oplus\ \Bigl(\frac{2\alpha}nI+\varepsilon A(\Fo)\Bigr).
\]
\emph{Step 2: counting.} The first block gives $p$ negative eigenvalues; the eigenvalues of the second one are $\frac{2\alpha}n+\varepsilon\theta_i(\Fo)$. By~\eqref{eq:matching} applied to the union of paths $\Fo$, the hypothesis $j\le\mu(\Fo)$ gives $\theta_j(\Fo)<0$. If moreover
\begin{equation}\label{eq:condA}
\frac{2\alpha}{n}+\varepsilon\,\theta_j(\Fo)\le0,
\end{equation}
then $\theta_1(\Fo),\dots,\theta_j(\Fo)$ all satisfy the same inequality (they are smaller than or equal to $\theta_j(\Fo)$ and $\varepsilon>0$), so $\npos(Y)\ge p+j=k_n$, and by~\ref{I1} of Lemma~\ref{lem:interlacing} the property $\mathcal P_\alpha$ holds.

\emph{Step 3: solving for $\alpha$.} Write $\theta=\theta_j(\Fo)<0$ and $\alpha=1-\varepsilon$. Then~\eqref{eq:condA} reads $2(1-\varepsilon)\le n\varepsilon|\theta|$, that is, $\varepsilon(2+n|\theta|)\ge2$, that is, $\varepsilon\ge\frac2{2-n\theta}$. So $\mathcal P_\alpha$ holds at $\alpha^\ast=1-\frac2{2-n\theta}$. If $\alpha^\ast\ge\frac12$, then $\alpha^\ast\in S(T)$ and $\alpha_0(T)\ge\alpha^\ast$; in any case $\alpha_0(T)\ge\frac12$ by Property~\ref{property}.

(b) \emph{Step 1: the submatrix.} Let $W=V_2\cup V_{\ge3}$ and $Y=\At[W]$. By~\ref{I2} of Lemma~\ref{lem:interlacing}, $\Nal(T)\le p+\npos(Y)$. Since $p+j=k_n$, it suffices to show that $\npos(Y)\le j-1$.

\emph{Step 2: a quadratic form estimate.} Write $x\in\mathbb R^W$ as $x=(y,\xi)$, with $y$ indexed by the hubs and $\xi$ by $V_2$. Then
\[
x^\top Yx=\sum_{v\in V_{\ge3}}(\alpha d_v-d_\alpha)y_v^2+2\varepsilon\!\!\sum_{\substack{uv\in E\\ u,v\in V_{\ge3}}}\!\!y_uy_v+2\varepsilon\!\!\sum_{\substack{vw\in E\\ v\in V_{\ge3},\,w\in V_2}}\!\!y_v\xi_w+\xi^\top\Bigl(\frac{2\alpha}nI+\varepsilon A(F_2)\Bigr)\xi .
\]
Fix $\delta=\frac\alpha6$; the assumption $\varepsilon\le\frac17$ gives $\alpha\ge\frac67$ and $\delta=\frac{1-\varepsilon}{6}\ge\frac17\ge\varepsilon$. We use two elementary inequalities:
\[
2\varepsilon y_uy_v\ge-\varepsilon(y_u^2+y_v^2),\qquad 2\varepsilon y_v\xi_w\ge-\delta y_v^2-\frac{\varepsilon^2}{\delta}\xi_w^2,
\]
the first one being $\varepsilon(y_u+y_v)^2\ge0$ and the second one being $\bigl(\sqrt\delta\,y_v+\frac{\varepsilon}{\sqrt\delta}\xi_w\bigr)^2\ge0$. In this way a hub $v$ loses at most $\varepsilon$ for each hub neighbour and $\delta$ for each degree-$2$ neighbour, hence at most $\delta d_v$ in total (recall $\varepsilon\le\delta$). A vertex $w$ of degree $2$ has at most two hub neighbours, so it loses at most $2\varepsilon^2/\delta=12\varepsilon^2/\alpha$. Therefore
\[
x^\top Yx\ \ge\ \sum_{v\in V_{\ge3}}\bigl(\alpha d_v-d_\alpha-\delta d_v\bigr)y_v^2+\xi^\top\Bigl(\frac{2\alpha}nI+\varepsilon A(F_2)-\frac{12\varepsilon^2}{\alpha}I\Bigr)\xi .
\]
\emph{Step 3: hubs are harmless.} Since $d_\alpha\le2\alpha$ and $d_v\ge3$, we have $\alpha d_v-d_\alpha-\delta d_v\ge(\alpha-\delta)d_v-2\alpha\ge3\cdot\frac{5\alpha}6-2\alpha=\frac\alpha2>0$. Hence, in the Loewner order,
\[
Y\ \succeq\ \frac\alpha2I_{n_{\ge3}}\ \oplus\ \Bigl(\frac{2\alpha}nI+\varepsilon A(F_2)-\frac{12\varepsilon^2}{\alpha}I\Bigr)=:Z .
\]
\emph{Step 4: counting.} By~\ref{W1} of Lemma~\ref{lem:weyl}, $\npos(Y)\le\npos(Z)$. The first block of $Z$ is positive definite, and the eigenvalues of the second one are $\frac{2\alpha}n+\varepsilon\theta_i(F_2)-\frac{12\varepsilon^2}\alpha$. Dividing by $\varepsilon>0$ and using $\alpha=1-\varepsilon$, such an eigenvalue is positive if and only if $\theta_i(F_2)>\frac{12\varepsilon}{1-\varepsilon}-\frac{2(1-\varepsilon)}{n\varepsilon}$. Under the hypothesis this holds for $i=j$, hence for all $i\ge j$ (the $\theta_i$ are ordered), and therefore $\npos(Z)\le j-1$. Here $j\le n_2$ by Lemma~\ref{lem:degree}, so $\theta_j(F_2)$ is defined.

(c) By~\eqref{eq:matching}, $j>\mu(F_2)$ is equivalent to $\theta_j(F_2)\ge0$. The hypothesis of (b) then holds as soon as $\frac{12\varepsilon}{1-\varepsilon}<\frac{2(1-\varepsilon)}{n\varepsilon}$, that is, $6n\varepsilon^2<(1-\varepsilon)^2$, that is, $\frac{\varepsilon}{1-\varepsilon}<\frac1{\sqrt{6n}}$, that is, $\varepsilon<\varepsilon_n:=\frac1{1+\sqrt{6n}}$. Note that $\varepsilon_n\le\frac17$ for $n\ge6$. Thus $\mathcal P_\alpha$ fails for every $\alpha\in(1-\varepsilon_n,1)$, and also at $\alpha=1$ (Proposition~\ref{prop:alpha1}). Since $\alpha_0(T)\in S(T)$, we conclude that $\alpha_0(T)\le1-\varepsilon_n$.
\end{proof}

\begin{remark}\label{rem:firstorder}
(i) For $T=S_n^\ast$ one has $j=1$ and $\Fo=P_{\ell_n-1}$ (the vertex $w_{\ell_n}$ is adjacent to the leaf $z$), with $\theta_1(P_{\ell_n-1})=-2\cos\frac{\pi}{\ell_n}$ by Lemma~\ref{lem:pathspec}. Hence Theorem~\ref{thm:firstorder}(a) gives $\alpha_0(S_n^\ast)\ge1-\frac{2}{2+2n\cos(\pi/\ell_n)}=1-\frac{1}{1+n\cos(\pi/\ell_n)}$, which is exactly Theorem~\ref{thm:spider}.

(ii) Every adjacency eigenvalue of a path is larger than $-2$ (Lemma~\ref{lem:pathspec}), so the first-order prediction $1-\frac2{2-n\theta_j}$ is always smaller than $1-\frac{2}{2+2n}=\frac{n}{n+1}$, in agreement with Theorem~\ref{thm:threshold}.

(iii) The two regimes of Theorem~\ref{thm:firstorder} are separated by a purely combinatorial condition: the comparison of $j(T)=k_n-p$ with the matching number of the forest of degree-$2$ vertices.

(iv) The constant $12$ in (b), and hence the constant $6$ in (c), can be improved at the price of the condition $\alpha\ge\frac34$; see Lemma~\ref{lem:loewner} and Remark~\ref{rem:improvedslow} in Section~\ref{sec:leaves}.
\end{remark}

\begin{conjecture}[First-order law]\label{conj:firstorder}
For every $\delta>0$, uniformly over the trees $T$ with $p\le k_n-1$ and $\theta_{j(T)}(F_2)\le-\delta$, one has
\[
1-\alpha_0(T)=\frac{2}{2-n\,\theta_{j(T)}(F_2)}+o\!\left(\frac1n\right).
\]
\end{conjecture}

The natural route to Conjecture~\ref{conj:firstorder} is to replace the Loewner estimate in the proof of (b) by the exact Schur complement of~\eqref{eq:blocks} onto the middle block. When $\lambda$ is an eigenvalue of $\At$ close to $0$, this Schur complement is the middle block corrected by $-\varepsilon^2E_{23}(H-\lambda I)^{-1}E_{23}^\top$, where $H$ is the hub block, and by $-\varepsilon^2E_{12}^\top(-\alpha\frac{n-2}{n}I-\lambda I)^{-1}E_{12}$, which is approximately $+\varepsilon^2/|x_1|$ at the vertices of the runs adjacent to leaves. These corrections act only at the ends of the runs, so the shift of the $j$-th eigenvalue is $O(\varepsilon^2)$ times the squared mass of the corresponding eigenvector at the ends of its run. This mass can be controlled with the explicit eigenvectors of Lemma~\ref{lem:pathspec}: for an eigenvector of $P_\ell$ normalized to have norm one, the squared entries at the ends are $\frac{2}{\ell+1}\sin^2\frac{k\pi}{\ell+1}$.

\begin{remark}[The slow regime, heuristically]
When $\theta_j(F_2)=0$, the $j$-th eigenvalue lives in the kernel of $A(F_2)$, which, by Lemma~\ref{lem:pathspec}, is spanned by one alternating vector $(1,0,-1,0,\dots)$ on each run of odd length. On this kernel the middle block of~\eqref{eq:blocks} acts as $\frac{2\alpha}{n}I$, and the decision is taken at second order by an effective matrix $\varepsilon^2K$: a run end adjacent to a leaf contributes $+\frac{2}{\ell+1}$, a run end adjacent to a hub of degree $d$ contributes $-\frac{2}{(\ell+1)(d-2)}$, and runs sharing a hub interact. This predicts $1-\alpha_0(T)\asymp\sqrt{2/(n|\tau|)}$ for an appropriate eigenvalue $\tau$ of $K$, which is consistent with the $n^{-1/2}$ upper bound of Theorem~\ref{thm:firstorder}(c).
\end{remark}

\section{The dependence on the number of leaves}\label{sec:leaves}

By Theorem~\ref{thm:dichotomy}, the number of leaves alone decides whether the midrange eigenvalue ever overtakes the average: this happens if and only if $p\le k_n-1$. In this section we show that the number of leaves also gives a quantitative estimate for the \emph{moment} at which this happens. More precisely, for every tree $T$ with $n$ vertices and $p\le k_n-1$ leaves we give an explicit number $\alpha_n(p)$, depending only on $n$ and $p$, such that $\lmid_\alpha(T)>d_\alpha$ for every $\alpha\in(\alpha_n(p),1]$ (Theorem~\ref{thm:leaves}), and we show that this estimate is sharp, up to an additive constant in $\frac{1}{1-\alpha}$, for a wide range of values of $p$ (Theorem~\ref{thm:leaves-sharp} and Corollary~\ref{cor:twosided}).

Throughout this section we use the notation
\begin{equation}\label{eq:Dkappa}
j=j_n(p)=k_n-p,\qquad D_n(p)=\min\bigl\{n-1,\ n+k_n-2p-1\bigr\},\qquad \kappa_n(p)=\cos\frac{\pi\,j_n(p)}{D_n(p)} .
\end{equation}
Note that $D_n(p)=n+k_n-2p-1$ if $j\le p$, that is, $2p\ge k_n$, or equivalently $p\ge\lceil k_n/2\rceil$; and $D_n(p)=n-1$ if $j\ge p$, that is, $p\le\lfloor k_n/2\rfloor$. When $k_n$ is even and $p=\frac{k_n}2$ both formulas give the same value.

\subsection{A refined quadratic form estimate}

The first ingredient is a version of Theorem~\ref{thm:firstorder}(b) with a better constant and without the condition $\varepsilon\le\frac17$.

\begin{lemma}\label{lem:loewner}
Let $T$ be a tree with $n\ge3$ vertices, let $W=V_2\cup V_{\ge3}$, let $\frac34\le\alpha<1$ and $\varepsilon=1-\alpha$. Then
\[
\npos\bigl(\At[W]\bigr)\ \le\ \#\Bigl\{i:\ \theta_i(F_2)\le\frac{6\varepsilon}{\alpha}-\frac{2\alpha}{n\varepsilon}\Bigr\}.
\]
Consequently, $\Nal(T)\le p+\#\{i:\theta_i(F_2)\le\frac{6\varepsilon}{\alpha}-\frac{2\alpha}{n\varepsilon}\}$, and if $p\le k_n-1$ and
\[
\theta_{j(T)}(F_2)\ >\ \frac{6\varepsilon}{\alpha}-\frac{2\alpha}{n\varepsilon},
\]
then $\mathcal P_\alpha(T)$ fails.
\end{lemma}

\begin{proof}
Write $Y=\At[W]$ and $x\in\mathbb R^W$ as $x=(y,\xi)$, with $y$ indexed by the hubs and $\xi$ by $V_2$. The edges of the tree $T[W]$ are of three kinds: edges between two vertices of degree $2$ (these are the edges of $F_2$), edges between a hub and a vertex of degree $2$, and edges between two hubs. Hence, by~\eqref{eq:entries},
\[
x^\top Yx=\sum_{v\in V_{\ge3}}(\alpha d_v-d_\alpha)y_v^2+2\varepsilon\!\!\sum_{\substack{vv'\in E\\ v,v'\in V_{\ge3}}}\!\!y_vy_{v'}+2\varepsilon\!\!\sum_{\substack{vw\in E\\ v\in V_{\ge3},\,w\in V_2}}\!\!y_v\xi_w+\xi^\top\Bigl(\frac{2\alpha}nI+\varepsilon A(F_2)\Bigr)\xi .
\]
For an edge $vv'$ between two hubs we use $2\varepsilon y_vy_{v'}\ge-\varepsilon y_v^2-\varepsilon y_{v'}^2$, which is $\varepsilon(y_v+y_{v'})^2\ge0$. For an edge $vw$ between a hub $v$ and a vertex $w$ of degree $2$ we use
\[
2\varepsilon y_v\xi_w\ \ge\ -\frac{\alpha}{3}\,y_v^2-\frac{3\varepsilon^2}{\alpha}\,\xi_w^2,
\]
which is $\bigl(\sqrt{\alpha/3}\,y_v+\varepsilon\sqrt{3/\alpha}\,\xi_w\bigr)^2\ge0$. Since $\alpha\ge\frac34$, we have $\varepsilon\le\frac14\le\frac\alpha3$. Hence every hub $v$ loses at most $\frac\alpha3$ for each of its neighbours in $W$, and therefore at most $\frac\alpha3d_v$ in total, while every vertex $w$ of degree $2$ loses $\frac{3\varepsilon^2}{\alpha}$ for each of its hub neighbours, of which there are at most two, hence at most $\frac{6\varepsilon^2}{\alpha}$ in total. Bounding the loss of every vertex of degree $2$ by $\frac{6\varepsilon^2}{\alpha}$ (which can only decrease the right-hand side) we obtain
\[
x^\top Yx\ \ge\ \sum_{v\in V_{\ge3}}\Bigl(\frac{2\alpha}{3}d_v-d_\alpha\Bigr)y_v^2+\xi^\top\Bigl(\frac{2\alpha}nI+\varepsilon A(F_2)-\frac{6\varepsilon^2}{\alpha}I\Bigr)\xi .
\]
For a hub, $d_v\ge3$ and $d_\alpha=2\alpha-\frac{2\alpha}{n}$, so $\frac{2\alpha}{3}d_v-d_\alpha\ge2\alpha-2\alpha+\frac{2\alpha}{n}=\frac{2\alpha}{n}>0$. Therefore, in the Loewner order,
\[
Y\ \succeq\ \frac{2\alpha}{n}I_{n_{\ge3}}\ \oplus\ \Bigl(\frac{2\alpha}nI+\varepsilon A(F_2)-\frac{6\varepsilon^2}{\alpha}I\Bigr)=:Z .
\]
By~\ref{W1} of Lemma~\ref{lem:weyl}, $\npos(Y)\le\npos(Z)$. The first block of $Z$ is positive definite, and the eigenvalues of the second block are $\frac{2\alpha}{n}+\varepsilon\theta_i(F_2)-\frac{6\varepsilon^2}{\alpha}$; dividing by $\varepsilon>0$, such an eigenvalue is nonpositive if and only if $\theta_i(F_2)\le\frac{6\varepsilon}{\alpha}-\frac{2\alpha}{n\varepsilon}$. This proves the first inequality.

The second one follows from $\Nal(T)=\npos(\At)\le(n-|W|)+\npos(\At[W])=p+\npos(\At[W])$, by~\ref{I2} of Lemma~\ref{lem:interlacing}. Finally, if $p\le k_n-1$, then $1\le j(T)\le n_2$ by Lemma~\ref{lem:degree}, so $\theta_{j(T)}(F_2)$ is defined. If $\theta_{j(T)}(F_2)>\frac{6\varepsilon}{\alpha}-\frac{2\alpha}{n\varepsilon}$, then, since the $\theta_i$ are ordered, at most $j(T)-1$ of them satisfy the opposite inequality, and $\Nal(T)\le p+j(T)-1=k_n-1$; so $\mathcal P_\alpha(T)$ fails.
\end{proof}

\begin{remark}\label{rem:improvedslow}
Lemma~\ref{lem:loewner} improves the slow regime of Theorem~\ref{thm:firstorder}(c). Indeed, if $\theta_{j(T)}(F_2)\ge0$, the failure condition of Lemma~\ref{lem:loewner} holds as soon as $\frac{6\varepsilon}{\alpha}<\frac{2\alpha}{n\varepsilon}$, that is, $3n\varepsilon^2<(1-\varepsilon)^2$, that is, $\varepsilon<\frac{1}{1+\sqrt{3n}}$; and $\frac{1}{1+\sqrt{3n}}\le\frac14$ for $n\ge3$, so that $\alpha\ge\frac34$ in this range. Arguing as in the proof of Theorem~\ref{thm:firstorder}(c) we obtain $\alpha_0(T)\le1-\frac{1}{1+\sqrt{3n}}$ for every tree with $p\le k_n-1$ and $j(T)>\mu(F_2)$.
\end{remark}

\subsection{Small eigenvalues of a union of paths}

The second ingredient is a lower bound for the $j$-th smallest eigenvalue of a union of paths, in terms of its number of vertices and of its number of components.

\begin{lemma}\label{lem:pathcount}
Let $F$ be a disjoint union of $r\ge1$ paths with $N$ vertices in total, and let $1\le j\le N$. Then
\[
\theta_j(F)\ \ge\ -2\cos\frac{\pi j}{N+\min\{j,r\}} .
\]
Equality holds, for instance, when $F$ consists of $r$ copies of $P_\ell$ and $j=sr$ for some integer $1\le s\le\ell$.
\end{lemma}

\begin{proof}
By Lemma~\ref{lem:pathspec}, every eigenvalue of a path lies in the open interval $(-2,2)$, so we can write $\theta_j(F)=-2\cos x$ for a unique $x\in(0,\pi)$. For a component $P_\ell$ of $F$, let $m_\ell$ be the number of its eigenvalues that are less than or equal to $-2\cos x$. Since $-2\cos x=2\cos(\pi-x)$ and $\cos$ is decreasing on $[0,\pi]$, an eigenvalue $2\cos\frac{k\pi}{\ell+1}$ of $P_\ell$ satisfies this inequality if and only if $\frac{k\pi}{\ell+1}\ge\pi-x$, that is, $k'\le\frac{(\ell+1)x}{\pi}$, where $k'=\ell+1-k$. As $k$ runs over $\{1,\dots,\ell\}$, $k'$ also runs over $\{1,\dots,\ell\}$, and $\frac{(\ell+1)x}{\pi}<\ell+1$; therefore
\[
m_\ell=\Bigl\lfloor\frac{(\ell+1)x}{\pi}\Bigr\rfloor .
\]
This gives two estimates. First, $m_\ell\le\frac{(\ell+1)x}{\pi}$. Second, $m_\ell(\pi-x)\le\ell x$: this is trivial if $m_\ell=0$, and if $m_\ell\ge1$ it follows from $m_\ell\pi\le(\ell+1)x=\ell x+x\le\ell x+m_\ell x$.

Since $\theta_1(F)\le\dots\le\theta_j(F)=-2\cos x$, at least $j$ eigenvalues of $F$ are less than or equal to $-2\cos x$, so $j\le\sum m_\ell$, the sum running over the $r$ components of $F$. Summing the first estimate over the components we get $j\le\frac{(N+r)x}{\pi}$, that is, $x\ge\frac{\pi j}{N+r}$. Summing the second one we get $j(\pi-x)\le Nx$, that is, $x\ge\frac{\pi j}{N+j}$. Hence $x\ge x_0:=\frac{\pi j}{N+\min\{j,r\}}$, and $x_0\in(0,\pi)$. Since $-2\cos$ is increasing on $[0,\pi]$, we conclude that $\theta_j(F)=-2\cos x\ge-2\cos x_0$.

For the equality case, let $F=rP_\ell$ and $j=sr$ with $1\le s\le\ell$. The $s$ smallest eigenvalues of $P_\ell$ are $-2\cos\frac{k\pi}{\ell+1}$, $1\le k\le s$, and each of them has multiplicity $r$ in $F$; hence $\theta_{sr}(F)=-2\cos\frac{s\pi}{\ell+1}$. On the other hand $N=r\ell$ and $\min\{j,r\}=r$, so $\frac{\pi j}{N+\min\{j,r\}}=\frac{\pi sr}{r(\ell+1)}=\frac{s\pi}{\ell+1}$.
\end{proof}

The third ingredient is a simple counting argument relating the runs to the number of leaves.

\begin{lemma}\label{lem:n2r}
Let $T$ be a tree with $n\ge3$ vertices, $p$ leaves, $n_2$ vertices of degree $2$ and $r$ runs.
\begin{enumerate}[label=(\alph*)]
\item $n_2+r\le n-1$.
\item If $T$ has a hub, then $n_2\le n-p-1$.
\item If $n\ge7$ and $p\le k_n-1$, then $D_T:=n_2+\min\{j(T),r\}$ satisfies $D_T\le D_n(p)$, and $D_n(p)>2j(T)$.
\end{enumerate}
\end{lemma}

\begin{proof}
(a) Consider the set of edges of $T$ having at least one endpoint of degree $2$. If both endpoints of such an edge have degree $2$, they are adjacent in $F_2$ and belong to the same run. Hence each such edge can be assigned to the unique run containing its endpoints of degree $2$, and different runs receive disjoint sets of edges. A run $w_1w_2\cdots w_\ell$ receives its $\ell-1$ internal edges and the edges joining its vertices to vertices outside the run; since the sum of the degrees of $w_1,\dots,w_\ell$ in $T$ is $2\ell=2(\ell-1)+2$, there are exactly two edges of the second type. So each run of length $\ell$ receives exactly $\ell+1$ edges. Summing over the $r$ runs, $n_2+r\le|E(T)|=n-1$.

(b) This follows from $n_2=n-p-n_{\ge3}$ and $n_{\ge3}\ge1$.

(c) If $T$ has a hub, then by (a) and (b)
\[
D_T\le n_2+r\le n-1\qquad\text{and}\qquad D_T\le n_2+j(T)\le n-p-1+k_n-p=n+k_n-2p-1,
\]
so $D_T\le D_n(p)$. If $T$ has no hub, then all its degrees are at most $2$ and $T$ is the path $P_n$; then $p=2$, $n_2=n-2$, $r=1$ and $j(T)=k_n-2\ge2$ (because $k_n\ge4$ for $n\ge7$), so $D_T=n-1$, while $D_n(2)=\min\{n-1,n+k_n-5\}=n-1$ because $k_n\ge4$. For the last claim, $2j(T)=2k_n-2p\le n+1-2p<n-1$ because $p\ge2$, and $2k_n-2p<n+k_n-2p-1$ because $k_n<n-1$.
\end{proof}

\subsection{\texorpdfstring{An upper bound for $\alpha_0$ in terms of the number of leaves}{An upper bound for alpha0 in terms of the number of leaves}}

We are now ready to answer the question of when the midrange eigenvalue overtakes the average, in terms of the number of leaves.

\begin{theorem}\label{thm:leaves}
Let $T$ be a tree with $n\ge7$ vertices and $p\le k_n-1$ leaves, and let $\kappa=\kappa_n(p)$ be as in~\eqref{eq:Dkappa}. Then $0<\kappa<1$, and the function
\[
f_n(\alpha)=\frac{\alpha}{n(1-\alpha)}-\frac{3(1-\alpha)}{\alpha}
\]
is a strictly increasing bijection from $(0,1)$ onto $\mathbb R$. Let $\alpha_n(p)$ be the unique solution of $f_n(\alpha)=\kappa$. Then $\alpha_n(p)\in(\frac34,1)$ and
\[
\lmid_\alpha(T)>d_\alpha\qquad\text{for every }\alpha\in\bigl(\alpha_n(p),1\bigr].
\]
In particular $\alpha_0(T)\le\alpha_n(p)$ and $\alpha_{\max}(n,p)\le\alpha_n(p)$. Moreover,
\begin{equation}\label{eq:explicitleaves}
\frac{1}{1-\alpha_0(T)}\ \le\ \frac{1}{1-\alpha_n(p)}\ \le\ 1+n\,\kappa_n(p)+\min\Bigl\{\frac{3}{\kappa_n(p)},\ \sqrt{3n}\Bigr\}.
\end{equation}
\end{theorem}

\begin{proof}
\emph{Step 1: the number $\kappa$.} By Lemma~\ref{lem:n2r}(c), $0<\frac{\pi j}{D_n(p)}<\frac\pi2$, so $0<\kappa<1$.

\emph{Step 2: the function $f_n$.} On $(0,1)$ the function $\alpha\mapsto\frac{\alpha}{n(1-\alpha)}$ is continuous and strictly increasing, and so is $\alpha\mapsto-\frac{3(1-\alpha)}{\alpha}=3-\frac3\alpha$. Hence $f_n$ is continuous and strictly increasing; moreover $f_n(\alpha)\to-\infty$ as $\alpha\to0^+$ and $f_n(\alpha)\to+\infty$ as $\alpha\to1^-$. By the intermediate value theorem $f_n$ is a bijection onto $\mathbb R$, and $\alpha_n(p)$ is well defined. Since $f_n(\frac34)=\frac{3}{n}-1<0<\kappa$, we have $\alpha_n(p)>\frac34$.

\emph{Step 3: a lower bound for $\theta_j(F_2)$.} Let $j=j(T)$ and $r$ be the number of runs of $T$. By Lemma~\ref{lem:degree}, $1\le j\le n_2$, so Lemma~\ref{lem:pathcount} applies to $F=F_2$, with $N=n_2$, and gives $\theta_j(F_2)\ge-2\cos\frac{\pi j}{D_T}$, where $D_T=n_2+\min\{j,r\}$. By Lemma~\ref{lem:n2r}(c), $D_T\le D_n(p)$, so $\frac{\pi j}{D_T}\ge\frac{\pi j}{D_n(p)}$; both angles lie in $(0,\pi)$, where $-2\cos$ is increasing, hence
\begin{equation}\label{eq:thetaj-kappa}
\theta_j(F_2)\ \ge\ -2\cos\frac{\pi j}{D_n(p)}=-2\kappa .
\end{equation}

\emph{Step 4: failure of $\mathcal P_\alpha$.} Let $\alpha\in(\alpha_n(p),1)$ and $\varepsilon=1-\alpha$. Then $\alpha>\frac34$ and $f_n(\alpha)>f_n(\alpha_n(p))=\kappa$, that is,
\[
\frac{\alpha}{n\varepsilon}-\frac{3\varepsilon}{\alpha}>\kappa
\qquad\Longleftrightarrow\qquad
\frac{6\varepsilon}{\alpha}-\frac{2\alpha}{n\varepsilon}<-2\kappa .
\]
Combining with~\eqref{eq:thetaj-kappa}, $\theta_j(F_2)>\frac{6\varepsilon}{\alpha}-\frac{2\alpha}{n\varepsilon}$, and Lemma~\ref{lem:loewner} shows that $\mathcal P_\alpha(T)$ fails, that is, $\lmid_\alpha(T)>d_\alpha$. At $\alpha=1$ the property fails by Proposition~\ref{prop:alpha1}. Hence $S(T)\subseteq[\frac12,\alpha_n(p)]$ and $\alpha_0(T)\le\alpha_n(p)$. Since this holds for every tree with $n$ vertices and $p$ leaves, $\alpha_{\max}(n,p)\le\alpha_n(p)$.

\emph{Step 5: the explicit bound.} Let $w>0$ and $\hat\varepsilon=\frac{1}{1+n\kappa+w}\in(0,1)$, so that $\frac{1-\hat\varepsilon}{\hat\varepsilon}=n\kappa+w$. Then
\[
f_n(1-\hat\varepsilon)=\frac{1-\hat\varepsilon}{n\hat\varepsilon}-\frac{3\hat\varepsilon}{1-\hat\varepsilon}=\kappa+\frac{w}{n}-\frac{3}{n\kappa+w}.
\]
This is at least $\kappa$ provided $w(n\kappa+w)\ge3n$. The choice $w=\frac3\kappa$ satisfies this condition, because $w\,n\kappa=3n$, and the choice $w=\sqrt{3n}$ satisfies it, because $w^2=3n$. For both choices $f_n(1-\hat\varepsilon)\ge\kappa=f_n(\alpha_n(p))$, and since $f_n$ is increasing, $\alpha_n(p)\le1-\hat\varepsilon$, that is, $\frac{1}{1-\alpha_n(p)}\le\frac1{\hat\varepsilon}=1+n\kappa+w$. Taking the better of the two choices gives~\eqref{eq:explicitleaves}, the first inequality of which follows from $\alpha_0(T)\le\alpha_n(p)$.
\end{proof}

Let us comment on the meaning of Theorem~\ref{thm:leaves}. The angle $\frac{\pi j}{D_n(p)}$ measures the \emph{leaf deficit} $j=\lceil n/2\rceil-p$, that is, the number of small eigenvalues that the vertices of degree $2$ must provide, relative to the room $D_n(p)$ available for the runs. The fewer leaves a tree has, the larger the deficit, the smaller $\kappa_n(p)$, and the smaller the value $\alpha_n(p)$ beyond which the midrange eigenvalue is guaranteed to stay above the average.

\begin{remark}[The three regimes]\label{rem:regimes}
(i) \emph{Almost exceptional trees.} If $p=k_n-1$, then $j=1$, $D_n(p)=n-k_n+1=\lfloor n/2\rfloor+1$ and $\kappa_n(p)=\cos\frac{\pi}{\lfloor n/2\rfloor+1}=1-O(n^{-2})$. In this case~\eqref{eq:explicitleaves} gives $\frac{1}{1-\alpha_0(T)}\le n+O(1)$, which has the correct order of magnitude but is weaker than Theorem~\ref{thm:threshold}, which gives $\frac{1}{1-\alpha_0(T)}<n+1$.

(ii) \emph{Linear number of leaves.} If $p/n$ stays in a compact subset of $(0,\frac12)$, then $\kappa_n(p)$ stays bounded away from $0$, and~\eqref{eq:explicitleaves} gives $1-\alpha_0(T)\ge\frac{1}{n\kappa_n(p)}\bigl(1-O(\frac1n)\bigr)$. This is the regime in which the bound is sharp (Corollary~\ref{cor:twosided}).

(iii) \emph{Few leaves.} If $p$ is fixed, or more generally $p=o(\sqrt n)$, then for $n$ large enough we have $p\le\lfloor k_n/2\rfloor$, and a direct computation from~\eqref{eq:Dkappa} gives $D_n(p)=n-1$ and $\kappa_n(p)=\sin\frac{\pi(p-1)}{n-1}$ for odd $n$, $\kappa_n(p)=\sin\frac{\pi(p-1/2)}{n-1}$ for even $n$, so that $n\kappa_n(p)\le\pi p\,\frac{n}{n-1}$. Then~\eqref{eq:explicitleaves} gives $\frac{1}{1-\alpha_0(T)}\le\sqrt{3n}\,(1+o(1))$: the crossing happens at distance at least of order $n^{-1/2}$ from $1$. This agrees with the slow regime of Theorem~\ref{thm:firstorder}(c) and Remark~\ref{rem:improvedslow}.
\end{remark}

As a first consequence, the leaf density of a sequence of trees gives an asymptotic lower bound for $n(1-\alpha_0)$. For $0<\beta<\frac12$ we define
\begin{equation}\label{eq:Lambda}
\Lambda(\beta)=\begin{cases}\dfrac{2}{1-2\beta}, & 0<\beta\le\frac14,\\[3mm] \dfrac{3-4\beta}{1-2\beta}=2+\dfrac{1}{1-2\beta}, & \frac14\le\beta<\frac12 .\end{cases}
\end{equation}
Both expressions are equal to $4$ at $\beta=\frac14$, so $\Lambda$ is continuous and strictly increasing, with $\Lambda(\beta)\to2$ as $\beta\to0$ and $\Lambda(\beta)\to\infty$ as $\beta\to\frac12$. In particular $\cos\frac{\pi}{\Lambda(\beta)}>0$.

\begin{corollary}\label{cor:density}
Let $\beta\in(0,\frac12)$ and let $(T_n)$ be a sequence of trees, where $T_n$ has $n$ vertices and $p_n\le k_n-1$ leaves, with $p_n/n\to\beta$. Then
\[
\liminf_{n\to\infty}\ n\bigl(1-\alpha_0(T_n)\bigr)\ \ge\ \frac{1}{\cos\bigl(\pi/\Lambda(\beta)\bigr)} .
\]
\end{corollary}

\begin{proof}
Since $k_n/n\to\frac12$ and $p_n/n\to\beta$, we have $j_n(p_n)/n\to\frac12-\beta$ and
\[
\frac{D_n(p_n)}{n}=\min\Bigl\{1-\frac1n,\ 1+\frac{k_n}{n}-\frac{2p_n}{n}-\frac1n\Bigr\}\ \longrightarrow\ \min\Bigl\{1,\ \frac32-2\beta\Bigr\},
\]
which equals $1$ if $\beta\le\frac14$ and $\frac32-2\beta$ if $\beta\ge\frac14$. Hence $\frac{\pi j_n(p_n)}{D_n(p_n)}\to\frac{\pi(\frac12-\beta)}{\min\{1,\frac32-2\beta\}}=\frac{\pi}{\Lambda(\beta)}$, and $\kappa_n(p_n)\to\kappa_\infty:=\cos\frac{\pi}{\Lambda(\beta)}>0$. By~\eqref{eq:explicitleaves},
\[
n\bigl(1-\alpha_0(T_n)\bigr)\ \ge\ \frac{n}{1+n\kappa_n(p_n)+3/\kappa_n(p_n)}\ \longrightarrow\ \frac{1}{\kappa_\infty}. \qedhere
\]
\end{proof}

\begin{remark}[Typical trees]\label{rem:typicaldensity}
The number of leaves of a uniformly random unlabeled tree with $n$ vertices satisfies $p/n\to m_\ast=0.4381562356\ldots$ in probability~\cite{RS}, and $m_\ast\in(\frac14,\frac12)$. Since the bound~\eqref{eq:explicitleaves} depends only on $n$ and $p$, and $\Lambda$ is continuous, Corollary~\ref{cor:density} suggests that for most trees
\[
n\bigl(1-\alpha_0(T)\bigr)\ \gtrsim\ \frac{1}{\cos\bigl(\pi/\Lambda(m_\ast)\bigr)}=1.0505\ldots,\qquad \Lambda(m_\ast)=\frac{3-4m_\ast}{1-2m_\ast}=10.0848\ldots,
\]
that is, the crossing of a typical tree happens at distance at least about $\frac{1.05}{n}$ from $1$. By Corollary~\ref{cor:profile}, the latest crossing among the trees with leaf density $m_\ast$ happens at distance between about $\frac{1.0505}{n}$ and $\frac{1}{n\cos(\pi/10)}=\frac{1.0514\ldots}{n}$ from $1$ (see the marked point in Figure~\ref{fig:profile}), since $\lceil\Lambda(m_\ast)\rceil-1=10$. We expect typical trees to cross considerably earlier.
\end{remark}

\subsection{Sharpness: subdivided trees}

We now construct, for a wide range of values of $p$, trees with $p$ leaves for which the bound of Theorem~\ref{thm:leaves} is attained up to an additive constant. The idea, suggested by the equality case of Lemma~\ref{lem:pathcount}, is to distribute the vertices of degree $2$ among exactly $j$ runs of almost equal lengths, and to arrange that each run provides exactly one small eigenvalue.

\begin{definition}\label{def:Tnp}
Let $n\ge7$ and let $p$ be an integer with $3\le p\le k_n-1$ and $j:=k_n-p\le2p-3$. Put
\[
h=\max\{1,\ j-p+1\},
\]
so that $1\le h\le p-2$. Let $R$ be the following tree with $p$ leaves, $h$ hubs and no vertex of degree $2$:
\begin{itemize}
\item if $h=1$, $R$ is the star $K_{1,p}$;
\item if $h\ge2$, $R$ is the caterpillar with spine $u_1u_2\cdots u_h$, two leaves attached to $u_1$, one leaf attached to each $u_i$ with $2\le i\le h-1$, and $p-h$ leaves attached to $u_h$ (note $p-h\ge2$).
\end{itemize}
The tree $R$ has $p+h$ vertices and $p+h-1$ edges. Let $E^\ast$ be a set of $j$ edges of $R$: if $h=1$, then $j\le p$ and $E^\ast$ consists of $j$ of the $p$ edges of the star; if $h\ge2$, then $p+h-1=j$ and $E^\ast=E(R)$. Let $N^\ast=n-p-h$, $\ell=\lfloor N^\ast/j\rfloor$ and $s=N^\ast-j\ell\in\{0,\dots,j-1\}$. The tree $T_{n,p}$ is obtained from $R$ by subdividing $s$ of the edges of $E^\ast$ with $\ell+1$ new vertices each, and the remaining $j-s$ edges of $E^\ast$ with $\ell$ new vertices each.
\end{definition}

Let us check that $T_{n,p}$ is well defined and record its parameters. If $h\ge2$, then $h=j-p+1\le p-2$ because $j\le2p-3$; and $h\ge2$ means $j\ge p+1>p$, while $h=1$ means $j\le p$. In the caterpillar, $u_1$ has degree $3$, $u_2,\dots,u_{h-1}$ have degree $3$ and $u_h$ has degree $1+(p-h)\ge3$, and the number of leaves is $2+(h-2)+(p-h)=p$. Subdividing an edge does not change the degrees of the existing vertices and creates vertices of degree $2$; hence $T_{n,p}$ is a tree with $p+h+N^\ast=n$ vertices, exactly $p$ leaves, exactly $h$ hubs, and $n_2=N^\ast$ vertices of degree $2$, organized in $j$ runs of lengths $\ell$ or $\ell+1$. Moreover
\[
N^\ast+j=D_n(p).
\]
Indeed, if $h=1$ then $j\le p$ and $N^\ast+j=n-p-1+k_n-p=n+k_n-2p-1=D_n(p)$; if $h=j-p+1$ then $j>p$ and $N^\ast+j=n-p-(j-p+1)+j=n-1=D_n(p)$. Consequently
\[
\ell+1=\Bigl\lfloor\frac{N^\ast}{j}\Bigr\rfloor+1=\Bigl\lfloor\frac{D_n(p)}{j}\Bigr\rfloor=:q .
\]
Finally, $q\ge3$, that is, $D_n(p)\ge3j$: if $j\le p$, this reads $n+k_n-2p-1\ge3k_n-3p$, that is, $p\ge2k_n-n+1$, which holds because $2k_n-n+1\le2<p$; if $j>p$, then $3j\le2j+2p-3=2k_n-3\le n-2<n-1$. So every run of $T_{n,p}$ has at least $2$ vertices. Note also that, for $j=1$, the tree $T_{n,k_n-1}$ is the spider $S_n^\ast$ of Definition~\ref{def:spider}.

\begin{figure}[ht]
\centering
\begin{tikzpicture}[scale=0.9]
\begin{scope}
\node[hub,label=below:$u$] (u) at (0,0) {};
\foreach \ang in {150,180,210} {\node[leaf] at (\ang:1.1) {}; \draw (u)--(\ang:1.1);}
\foreach \ang/\len in {60/3,20/3,-20/3,-60/4} {
  \draw (u)--(\ang:0.8*\len+0.8);
  \foreach \k in {1,...,\len} {\node[dtwo] at (\ang:0.8*\k) {};}
  \node[leaf] at (\ang:0.8*\len+0.8) {};
}
\node[hub] at (0,0) {};
\node at (0.6,-4.2) {\small $h=1$};
\end{scope}
\begin{scope}[shift={(8.2,0)}]
\coordinate (h1) at (0,0); \coordinate (h2) at (2.4,0);
\draw (h1)--(h2);
\foreach \x in {0.6,1.2,1.8} {\node[dtwo] at (\x,0) {};}
\foreach \ang/\len in {135/3,225/2} {
  \draw (h1)--++(\ang:0.6*\len+0.6);
  \foreach \k in {1,...,\len} {\node[dtwo] at ($(h1)+(\ang:0.6*\k)$) {};}
  \node[leaf] at ($(h1)+(\ang:0.6*\len+0.6)$) {};
}
\foreach \ang/\len in {45/3,0/2,-45/2} {
  \draw (h2)--++(\ang:0.6*\len+0.6);
  \foreach \k in {1,...,\len} {\node[dtwo] at ($(h2)+(\ang:0.6*\k)$) {};}
  \node[leaf] at ($(h2)+(\ang:0.6*\len+0.6)$) {};
}
\node[hub] at (h1) {}; \node[hub] at (h2) {};
\node at (1.2,-4.2) {\small $h\ge2$};
\end{scope}
\end{tikzpicture}
\caption{The trees $T_{n,p}$ of Definition~\ref{def:Tnp}. Left: a case with $j\le p$ (so $h=1$): $n=21$, $p=7$, $j=4$, $D_{21}(7)=17$, and four legs with $3,3,3,4$ vertices of degree $2$. Right: a case with $j>p$ (so $h=j-p+1\ge2$): $n=22$, $p=5$, $j=6$, $h=2$, $D_{22}(5)=21$, and every edge of the caterpillar $R$ subdivided, with runs of lengths $2$ or $3$. In both cases there are exactly $j$ runs, of lengths $\lfloor D_n(p)/j\rfloor-1$ or $\lfloor D_n(p)/j\rfloor$.}
\label{fig:Tnp}
\end{figure}
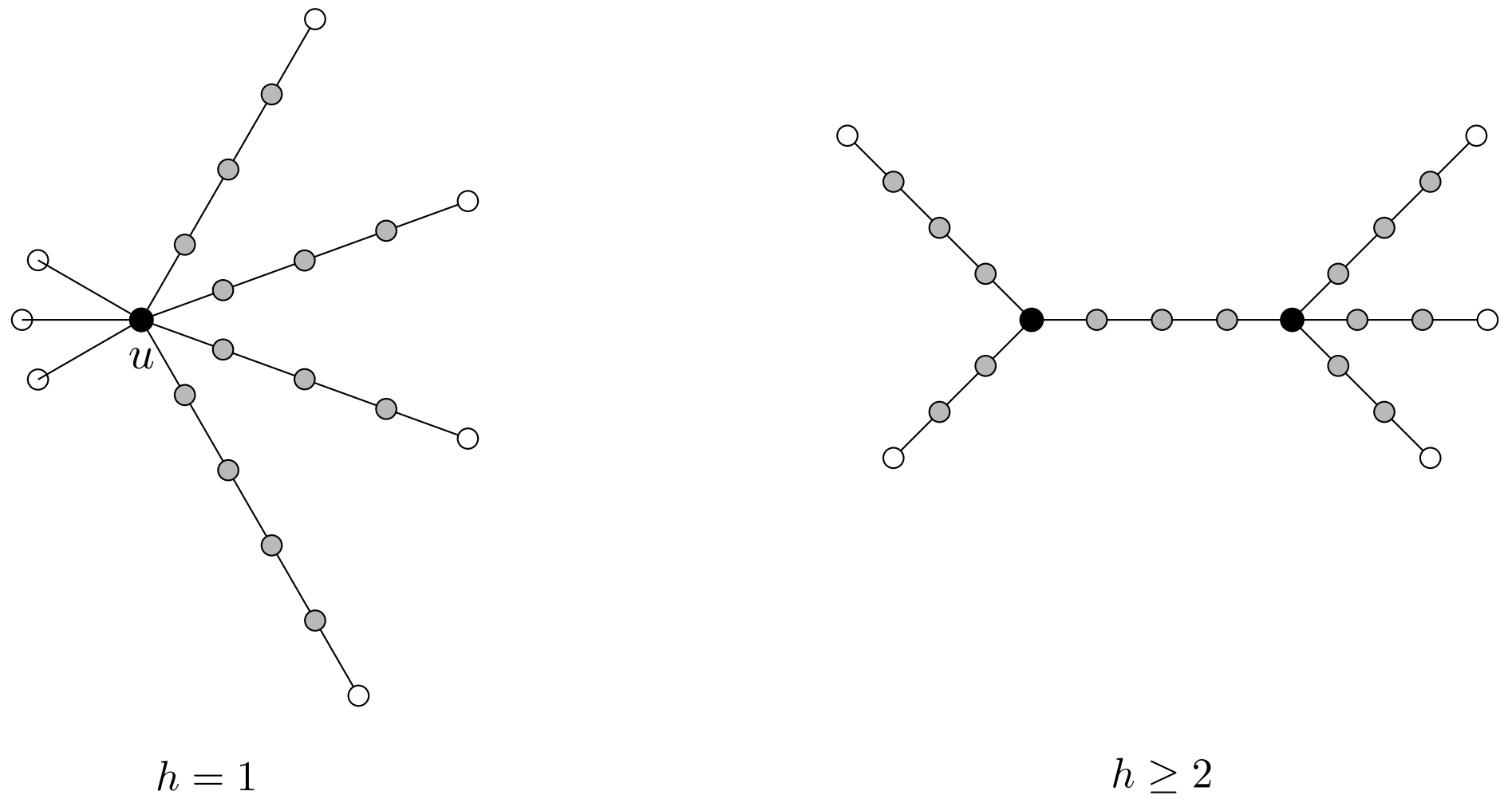

\begin{theorem}\label{thm:leaves-sharp}
Let $n\ge7$ and let $p$ be an integer with $3\le p\le k_n-1$ and $j=k_n-p\le2p-3$. Let $q=\lfloor D_n(p)/j\rfloor$ (so that $q\ge3$) and
\[
\psi^\ast=\frac{\pi+\frac{3}{n-1}}{q},\qquad \alpha^\ast=1-\frac{1}{1+n\cos\psi^\ast}.
\]
Then $\mathcal P_{\alpha^\ast}(T_{n,p})$ holds. Consequently
\[
\alpha_{\max}(n,p)\ \ge\ \alpha_0(T_{n,p})\ \ge\ \alpha^\ast,\qquad\text{that is,}\qquad \frac{1}{1-\alpha_{\max}(n,p)}\ \ge\ 1+n\cos\frac{\pi+\frac{3}{n-1}}{\bigl\lfloor D_n(p)/j\bigr\rfloor}.
\]
\end{theorem}

\begin{proof}
\emph{Step 1: the angle.} Since $q\ge3$ and $\frac{3}{n-1}\le\frac12$ for $n\ge7$, we have $0<\psi^\ast\le\frac{\pi+1/2}{3}<1.22$. Using the elementary inequality $\cos x\ge1-\frac{x^2}{2}+\frac{x^4}{24}-\frac{x^6}{720}$, valid for every real $x$, and the fact that $\cos$ is decreasing on $[0,\pi]$, we get $\cos\psi^\ast\ge\cos1.22>0.343>\frac{1}{\sqrt{10}}$. Hence $\psi^\ast\in(0,\frac\pi2)$ and $\tan^2\psi^\ast=\frac{1}{\cos^2\psi^\ast}-1<9$, that is, $\tan\psi^\ast<3$.

\emph{Step 2: the parameter $\alpha^\ast$.} Since $n\cos\psi^\ast>\frac{7}{\sqrt{10}}>1$, we have $\alpha^\ast=\frac{n\cos\psi^\ast}{1+n\cos\psi^\ast}\in(\frac12,1)$. Put $\varepsilon^\ast=1-\alpha^\ast=\frac{1}{1+n\cos\psi^\ast}$. Then $\frac{\alpha^\ast}{n\varepsilon^\ast}=\cos\psi^\ast\in(0,1)$, so we are in the situation of Lemma~\ref{lem:pendant}(c) with $\psi=\psi^\ast$, and the corresponding angle $\omega^\ast$ satisfies
\[
\omega^\ast<\tan\omega^\ast=\frac{\tan\psi^\ast}{n-1}<\frac{3}{n-1}.
\]

\emph{Step 3: the pieces of $T_{n,p}$.} Let $H$ be the set of hubs of $T_{n,p}$ and $U=V(T_{n,p})\setminus H$. The connected components of $T_{n,p}-H$ correspond to the edges of $R$, as follows. An edge of $R$ joining a hub to a leaf, not subdivided, gives an isolated leaf; subdivided with $m$ vertices, it gives a pendant path $z\,w_m\cdots w_1$ ending at a leaf $z$, which we call a \emph{leg}. An edge of $R$ joining two hubs, subdivided with $m$ vertices, gives a path $P_m$ of vertices of degree $2$, which we call an \emph{internal run}. (Edges of $R$ joining two hubs are never left unsubdivided: there are none if $h=1$, and all edges are subdivided if $h\ge2$.) Since there are no edges between different components, $\At(T_{n,p})[U]$ is the direct sum of the following blocks: $(x_1)$ for each isolated leaf, $M_m(\alpha)$ (as in Lemma~\ref{lem:pendant}) for each leg with $m$ vertices of degree $2$, and $cI_m+\varepsilon A(P_m)$ for each internal run with $m$ vertices.

\emph{Step 4: counting at $\alpha=\alpha^\ast$.} Every run has length $m\in\{q-1,q\}$, so $m+1\ge q$ and $m+1\le q+1$.
\begin{itemize}
\item An isolated leaf contributes one negative eigenvalue, $x_1<0$.
\item For a leg with $m$ vertices of degree $2$,
\[
(m+1)\psi^\ast-\omega^\ast\ \ge\ q\psi^\ast-\omega^\ast=\pi+\frac{3}{n-1}-\omega^\ast>\pi,
\]
and
\[
(m+1)\psi^\ast-\omega^\ast<(q+1)\psi^\ast=\Bigl(1+\frac1q\Bigr)\Bigl(\pi+\frac{3}{n-1}\Bigr)\le\frac43\Bigl(\pi+\frac12\Bigr)<4.86<2\pi .
\]
By Lemma~\ref{lem:pendant}(c), $\det M_m(\alpha^\ast)>0$, and by Lemma~\ref{lem:pendant}(b), $M_m(\alpha^\ast)$ has at least two negative eigenvalues.
\item For an internal run with $m$ vertices, the smallest eigenvalue of $cI_m+\varepsilon^\ast A(P_m)$ is, by Lemma~\ref{lem:pathspec},
\[
c-2\varepsilon^\ast\cos\frac{\pi}{m+1}=2\varepsilon^\ast\Bigl(\cos\psi^\ast-\cos\frac{\pi}{m+1}\Bigr)<0,
\]
because $c=\frac{2\alpha^\ast}{n}=2\varepsilon^\ast\cos\psi^\ast$ and $\frac{\pi}{m+1}<\psi^\ast$ (recall $(m+1)\psi^\ast\ge q\psi^\ast>\pi$), while $\cos$ is decreasing on $[0,\pi]$. So the internal run contributes at least one negative eigenvalue.
\end{itemize}
Every leaf of $R$ lies on exactly one edge of $R$, which is either left unsubdivided (isolated leaf) or subdivided (leg); and exactly $j$ edges of $R$ are subdivided, each giving a leg or an internal run. Adding up,
\[
\npos\bigl(\At[U]\bigr)\ \ge\ \#\{\text{isolated leaves}\}+2\,\#\{\text{legs}\}+\#\{\text{internal runs}\}=p+j=k_n .
\]
By~\ref{I1} of Lemma~\ref{lem:interlacing}, $\Nal(T_{n,p})\ge k_n$, that is, $\mathcal P_{\alpha^\ast}(T_{n,p})$ holds. Since $\alpha^\ast>\frac12$, we get $\alpha^\ast\in S(T_{n,p})$, so $\alpha_0(T_{n,p})\ge\alpha^\ast$, and $\alpha_{\max}(n,p)\ge\alpha_0(T_{n,p})$ because $T_{n,p}$ has $n$ vertices and $p$ leaves.
\end{proof}

\begin{remark}
For $j=1$ we have $T_{n,k_n-1}=S_n^\ast$ and $q=\lfloor n/2\rfloor+1$, so Theorem~\ref{thm:leaves-sharp} gives $\alpha_0(S_n^\ast)\ge1-\frac{1}{1+n\cos\frac{\pi+3/(n-1)}{\ell_n+1}}$, $\ell_n=\lfloor n/2\rfloor$. This slightly improves Theorem~\ref{thm:spider}, because $\frac{\pi+3/(n-1)}{\ell_n+1}<\frac{\pi}{\ell_n}$ (equivalently, $\frac{3\ell_n}{n-1}<\pi$, which holds since $\frac{3\ell_n}{n-1}\le\frac{3n}{2(n-1)}\le\frac74$). Both lower bounds for $\frac{1}{1-\alpha_0(S_n^\ast)}$ are of the form $n+1-\frac{2\pi^2}{n}+O(n^{-2})$.
\end{remark}

Combining Theorems~\ref{thm:leaves} and~\ref{thm:leaves-sharp} we obtain a two-sided estimate for the latest crossing among trees with a prescribed number of leaves.

\begin{corollary}[Two-sided estimate]\label{cor:twosided}
Let $n\ge7$ and let $p$ be an integer with $3\le p\le k_n-1$ and $j=k_n-p\le2p-3$, and let $\kappa=\kappa_n(p)$. Then
\[
1+n\cos\frac{\pi+\frac{3}{n-1}}{\lfloor D_n(p)/j\rfloor}\ \le\ \frac{1}{1-\alpha_{\max}(n,p)}\ \le\ 1+n\kappa+\frac{3}{\kappa}.
\]
If moreover $j$ divides $D_n(p)$, then
\[
n\kappa-\frac16\ \le\ \frac{1}{1-\alpha_{\max}(n,p)}\ \le\ n\kappa+7 ,
\]
so that $1-\alpha_{\max}(n,p)=\dfrac{1}{n\,\kappa_n(p)}\Bigl(1+O\bigl(\tfrac1n\bigr)\Bigr)$ with an absolute implied constant.
\end{corollary}

\begin{proof}
The lower bound is Theorem~\ref{thm:leaves-sharp} and the upper bound is~\eqref{eq:explicitleaves} applied to a tree attaining $\alpha_{\max}(n,p)$. Assume now that $j$ divides $D_n(p)$ and put $q=D_n(p)/j\ge3$. Then $\kappa=\cos\frac{\pi}{q}\ge\cos\frac\pi3=\frac12$, so $\frac3\kappa\le6$ and the upper bound is at most $n\kappa+7$. For the lower bound, since $|\cos a-\cos b|\le|a-b|$,
\[
\cos\frac{\pi+\frac{3}{n-1}}{q}\ \ge\ \cos\frac{\pi}{q}-\frac{3}{(n-1)q}\ \ge\ \kappa-\frac{1}{n-1},
\]
and therefore $1+n\cos\frac{\pi+3/(n-1)}{q}\ge1+n\kappa-\frac{n}{n-1}\ge n\kappa+1-\frac76=n\kappa-\frac16$, because $\frac{n}{n-1}\le\frac76$ for $n\ge7$. The last statement follows by taking reciprocals, since $n\kappa\ge\frac n2$.
\end{proof}

When the leaf density converges, Corollary~\ref{cor:twosided} gives the first-order behaviour of $\alpha_{\max}(n,p)$.

\begin{corollary}[Asymptotic profile]\label{cor:profile}
Let $\beta\in(\frac16,\frac12)$, let $\Lambda(\beta)$ be as in~\eqref{eq:Lambda}, and let $(p_n)$ be a sequence of integers with $2\le p_n\le k_n-1$ and $p_n/n\to\beta$. Then
\begin{align*}
\frac{1}{\cos\bigl(\pi/\Lambda(\beta)\bigr)}\ &\le\ \liminf_{n\to\infty}n\bigl(1-\alpha_{\max}(n,p_n)\bigr)\\
&\le\ \limsup_{n\to\infty}n\bigl(1-\alpha_{\max}(n,p_n)\bigr)\ \le\ \frac{1}{\cos\bigl(\pi/(\lceil\Lambda(\beta)\rceil-1)\bigr)} .
\end{align*}
Moreover, along every subsequence for which $k_n-p_n$ divides $D_n(p_n)$,
\[
n\bigl(1-\alpha_{\max}(n,p_n)\bigr)\ \longrightarrow\ \frac{1}{\cos(\pi/\Lambda(\beta))}.
\]
\end{corollary}

\begin{proof}
For the lower bound, let $T_n$ be a tree with $n$ vertices and $p_n$ leaves such that $\alpha_0(T_n)=\alpha_{\max}(n,p_n)$, and apply Corollary~\ref{cor:density}. For the upper bound, let $j_n=k_n-p_n$ and $D_n=D_n(p_n)$. Since $\beta>\frac16$, we have $j_n/n\to\frac12-\beta<2\beta$, so $j_n\le2p_n-3$ for large $n$; also $p_n\ge3$ and $p_n\le k_n-1$ for large $n$. Hence Theorem~\ref{thm:leaves-sharp} applies when $n$ is large. As in the proof of Corollary~\ref{cor:density}, $D_n/j_n\to\Lambda(\beta)$, and $\Lambda(\beta)>\Lambda(\frac16)=3$. Since $\Lambda(\beta)>\lceil\Lambda(\beta)\rceil-1$, for all large $n$ we have $D_n/j_n>\lceil\Lambda(\beta)\rceil-1$, and therefore $q_n:=\lfloor D_n/j_n\rfloor\ge\lceil\Lambda(\beta)\rceil-1\ge3$. By Theorem~\ref{thm:leaves-sharp},
\begin{align*}
n\bigl(1-\alpha_{\max}(n,p_n)\bigr)&\le\frac{n}{1+n\cos\frac{\pi+3/(n-1)}{q_n}}\le\frac{n}{1+n\cos\frac{\pi+3/(n-1)}{\lceil\Lambda(\beta)\rceil-1}}\\
&\longrightarrow\ \frac{1}{\cos\frac{\pi}{\lceil\Lambda(\beta)\rceil-1}} .
\end{align*}
Finally, along a subsequence where $j_n$ divides $D_n$, Corollary~\ref{cor:twosided} gives $\frac{n}{n\kappa_n(p_n)+7}\le n(1-\alpha_{\max}(n,p_n))\le\frac{n}{n\kappa_n(p_n)-1/6}$, and both bounds converge to $\frac{1}{\cos(\pi/\Lambda(\beta))}$, because $\kappa_n(p_n)\to\cos\frac{\pi}{\Lambda(\beta)}$.
\end{proof}

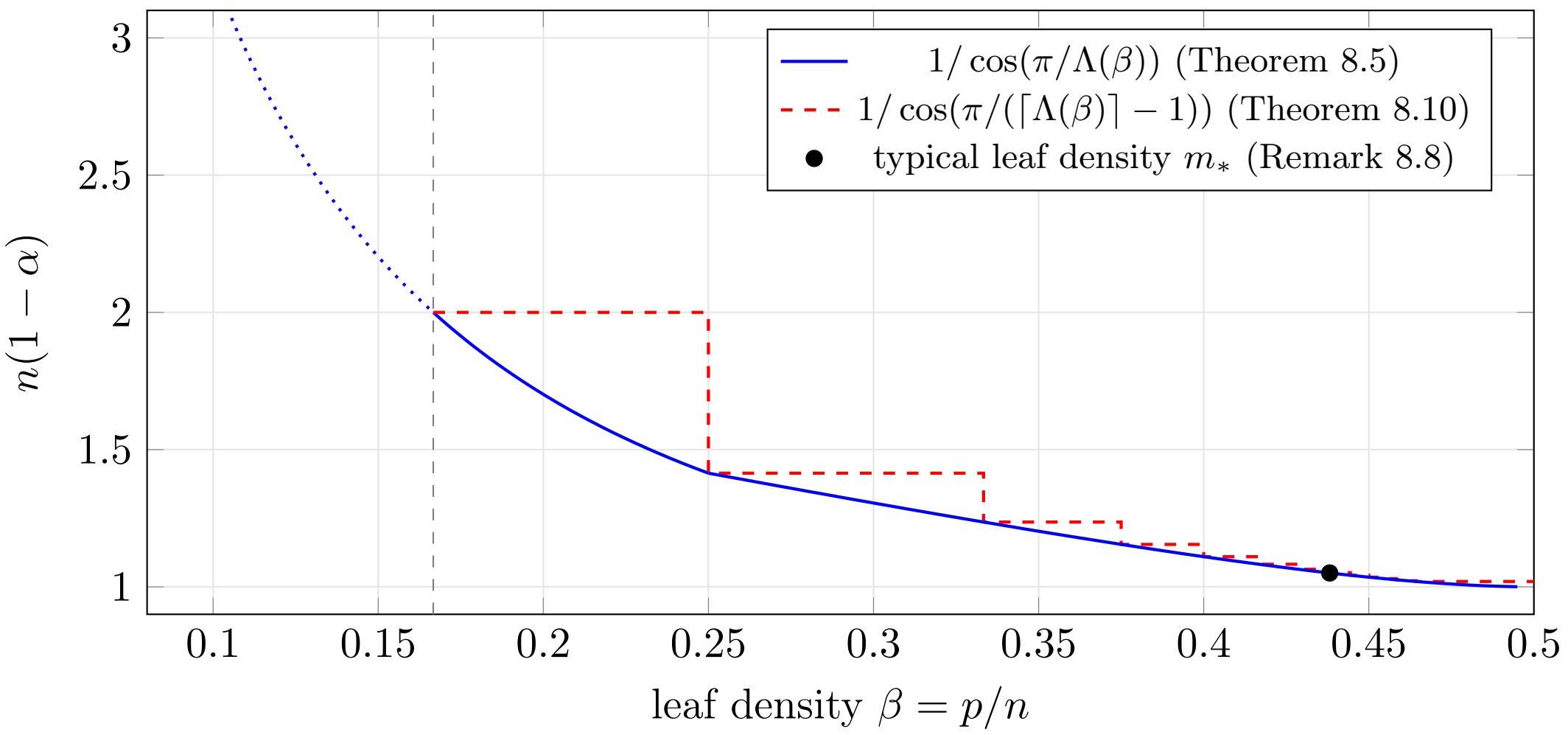
\begin{figure}[ht]
\centering
\begin{tikzpicture}
\begin{axis}[width=0.9\textwidth,height=7cm,xlabel={leaf density $\beta=p/n$},ylabel={$n(1-\alpha)$},xmin=0.08,xmax=0.5,ymin=0.9,ymax=3.1,grid=both,grid style={gray!20},legend pos=north east,legend style={font=\footnotesize}]
\addplot[domain=0.1667:0.25,samples=60,thick,blue] {1/cos(deg(pi*(1-2*x)/2))};
\addlegendentry{$1/\cos(\pi/\Lambda(\beta))$ (Theorem~\ref{thm:leaves})}
\addplot[const plot,thick,red,dashed] coordinates {(0.16667,2) (0.25000,1.41421) (0.33333,1.23607) (0.37500,1.15470) (0.40000,1.10992) (0.41667,1.08239) (0.42857,1.06418) (0.43750,1.05146) (0.44444,1.04222) (0.45000,1.03528) (0.45455,1.02993) (0.45833,1.02572) (0.46154,1.02234) (0.46429,1.01959) (0.5,1.01959)};
\addlegendentry{$1/\cos(\pi/(\lceil\Lambda(\beta)\rceil-1))$ (Theorem~\ref{thm:leaves-sharp})}
\addplot[domain=0.25:0.495,samples=80,thick,blue,forget plot] {1/cos(deg(pi*(1-2*x)/(3-4*x)))};
\addplot[domain=0.08:0.1667,samples=40,thick,blue,dotted,forget plot] {1/cos(deg(pi*(1-2*x)/2))};
\addplot[only marks,mark=*,mark size=2pt,black] coordinates {(0.43816,1.05056)};
\addlegendentry{typical leaf density $m_\ast$ (Remark~\ref{rem:typicaldensity})}
\draw[gray,dashed] (axis cs:0.16667,0.9)--(axis cs:0.16667,3.1);
\end{axis}
\end{tikzpicture}
\caption{The asymptotic profile of Corollary~\ref{cor:profile}. For a leaf density $\beta\in(\frac16,\frac12)$, the limit points of $n(1-\alpha_{\max}(n,p))$, $p/n\to\beta$, lie between the solid curve and the dashed step function, and they lie on the solid curve along divisible subsequences. For $\beta\le\frac16$ (dotted) only the lower bound of Corollary~\ref{cor:density} is available.}
\label{fig:profile}
\end{figure}

The function $\beta\mapsto1/\cos(\pi/\Lambda(\beta))$ is continuous and strictly decreasing on $(\frac16,\frac12)$; it equals $2$ at $\beta=\frac16$, $\sqrt2$ at $\beta=\frac14$, and tends to $1$ as $\beta\to\frac12$, in agreement with $n(1-\alpha_{\max}(n))\to1$ (Corollary~\ref{cor:asymp}). The upper bound in Corollary~\ref{cor:profile} is a step function which coincides with the lower one at the points where $\Lambda(\beta)$ is an integer, namely $\beta_q=\frac{q-3}{2q-4}$ for $q\ge4$ (with $\beta_4=\frac14$, $\beta_5=\frac13$, $\beta_6=\frac38$, \dots). Both curves are displayed in Figure~\ref{fig:profile}. In words: for the extremal trees, the crossing happens at
\[
\alpha\ \approx\ 1-\frac{1}{n\cos\bigl(\pi/\Lambda(p/n)\bigr)},
\]
and the smaller the leaf density, the earlier the midrange eigenvalue overtakes the average.

\begin{remark}[Trees with few leaves]\label{rem:fewleaves}
The hypothesis $j\le2p-3$ in Theorem~\ref{thm:leaves-sharp}, that is, $3p\ge k_n+3$, or equivalently $p\ge\lceil\frac{k_n+3}{3}\rceil$ (asymptotically $\beta>\frac16$), guarantees that a tree with $p$ leaves can have $j$ runs, each of which has to provide only one small eigenvalue. Indeed, by Lemma~\ref{lem:n2r}(a) and Lemma~\ref{lem:degree}, a tree with $p$ leaves has at most $r\le n-1-n_2=p+n_{\ge3}-1\le2p-3$ runs. When $j>2p-3$, some runs must provide several small eigenvalues. Lemma~\ref{lem:pendant}(d) shows that a leg with $m$ vertices of degree $2$ provides $1+\lfloor((m+1)\psi-\omega)/\pi\rfloor$ negative eigenvalues (when no $P_k$ vanishes), and Lemma~\ref{lem:pathspec} shows that an internal run with $m$ vertices provides $\lfloor(m+1)\psi/\pi\rfloor$ nonpositive ones, so trees with runs of suitable lengths can be built in the same way. However, the integrality of the run lengths then produces rounding effects of the same order as the main term, and we leave the precise first-order behaviour of $\alpha_{\max}(n,p)$ for $p\le\lceil\frac{k_n+3}{3}\rceil-1$ as an open problem. In this range Theorem~\ref{thm:leaves} still applies, and Remark~\ref{rem:regimes}(iii) shows that for $p=o(\sqrt n)$ the crossing happens at distance at least of order $n^{-1/2}$ from $1$.
\end{remark}

\section{Transfer to other matrix families}\label{sec:families}

In the same way that $A_\alpha$ interpolates between the adjacency matrix and the signless Laplacian, other one-parameter deformations of the classical graph matrices are frequently studied. Following~\cite[Section~10]{BO}, we consider two of them. Samanta, Deepshikha and Das~\cite{SDD} introduced the family
\[
B_\beta(G)=\beta A(G)+(1-\beta)L(G),\qquad 0\le\beta\le1,
\]
for which $B_0=L$, $B_{1/2}=\frac12D$, $B_{2/3}=\frac13Q$ and $B_1=A$. Morbidi~\cite{Morbidi} introduced, in the context of consensus protocols, the \emph{deformed Laplacian}
\[
M_G(s)=I-sA(G)+s^2\bigl(D(G)-I\bigr),\qquad s\in\mathbb R,
\]
for which $M_G(0)=I$, $M_G(1)=L$ and $M_G(-1)=Q$; see also~\cite{GHN,DOT}. As observed in~\cite[Proposition~10.1]{BO}, both families are affine reparametrisations of $A_\alpha$ on trees, and this allows us to transfer all the results of this paper to them.

\begin{remark}[The matrices $B_\beta$ and $M_T(s)$]\label{rem:families}
Let $T$ be a tree with $n$ vertices and $p$ leaves. For a real symmetric matrix $X$ of order $n$, let $\lmid(X)=\lambda_{k_n}(X)$ be its midrange eigenvalue and $d(X)=\frac1n\operatorname{tr}X$ the average of its eigenvalues; as in Lemma~\ref{lem:equivP}, $\lmid(X)\le d(X)$ if and only if at least $\lceil n/2\rceil$ eigenvalues of $X$ are less than or equal to $d(X)$, and we say that $X$ has property $\mathcal P$ when this holds.

\emph{(a) Affine relations.} Since $L=D-A$, we have $B_\beta=(2\beta-1)A+(1-\beta)D$. For $\alpha\in[0,1]$ put $\beta_\alpha=\frac{1}{1+\alpha}$; then $2\beta_\alpha-1=\frac{1-\alpha}{1+\alpha}$ and $1-\beta_\alpha=\frac{\alpha}{1+\alpha}$, so that
\begin{equation}\label{eq:Bbeta}
B_{\beta_\alpha}(T)=\frac{1}{1+\alpha}\,A_\alpha(T).
\end{equation}
For $t>0$ put $\alpha_t=\frac{t}{1+t}\in(0,1)$; then $t(1+t)A_{\alpha_t}=tA+t^2D$, and therefore
\begin{equation}\label{eq:Mdeformed}
M_T(-t)=I+tA+t^2(D-I)=(1-t^2)I+t(1+t)\,A_{\alpha_t}(T).
\end{equation}
Moreover, $M_T(s)$ and $M_T(-s)$ are similar. Indeed, a tree is bipartite; if $J$ is the diagonal matrix with entries $1$ on one class of a bipartition and $-1$ on the other, then $JAJ=-A$, while $J$ commutes with $D$ and $I$, so $JM_T(s)J=M_T(-s)$ and $J^{-1}=J$ (see also~\cite{DOT}).

\emph{(b) Transfer of the property.} Let $X=aI+bA_\alpha(T)$ with $a\in\mathbb R$ and $b>0$. Every eigenvector of $A_\alpha(T)$ with eigenvalue $\lambda$ is an eigenvector of $X$ with eigenvalue $a+b\lambda$, and $\lambda\mapsto a+b\lambda$ is increasing; hence $\lambda_i(X)=a+b\lambda_i(A_\alpha(T))$ for every $i$, and $d(X)=a+b\,d_\alpha$ by linearity of the trace. Consequently $\lmid(X)-d(X)=b\bigl(\lmid_\alpha(T)-d_\alpha\bigr)$, and $X$ has property $\mathcal P$ if and only if $\mathcal P_\alpha(T)$ holds. Similar matrices have the same eigenvalues and the same trace. By~\eqref{eq:Bbeta} and~\eqref{eq:Mdeformed} we conclude:
\begin{itemize}
\item for $\beta\in[\frac12,1]$, $B_\beta(T)$ has property $\mathcal P$ if and only if $\mathcal P_{\alpha}(T)$ holds with $\alpha=\frac1\beta-1$;
\item for $s\neq0$, $M_T(s)$ has property $\mathcal P$ if and only if $\mathcal P_{\alpha}(T)$ holds with $\alpha=\frac{|s|}{1+|s|}$; and $M_T(0)=I$ trivially has property $\mathcal P$.
\end{itemize}
The map $\beta\mapsto\frac1\beta-1$ is a decreasing bijection from $[\frac12,1]$ onto $[0,1]$, and $t\mapsto\frac{t}{1+t}$ is an increasing bijection from $(0,\infty)$ onto $(0,1)$. The range $0\le\alpha\le\frac12$ of Property~\ref{property} corresponds to $\frac23\le\beta\le1$ and to $|s|\le1$, which is~\cite[Proposition~10.1]{BO}, while the range $\frac12<\alpha<1$ studied here corresponds to $\frac12<\beta<\frac23$ and to $|s|>1$ (and $\alpha=1$ to $\beta=\frac12$). Note also that $\frac{1}{1-\alpha_t}=1+t$, and that $\frac{1}{1-\alpha}=\frac{\beta}{2\beta-1}$ when $\alpha=\frac1\beta-1$.

\emph{(c) Translation of the results.} Combining (b) with the results of this paper we obtain the following statements.
\begin{itemize}
\item \emph{Dichotomy.} $B_\beta(T)$ has property $\mathcal P$ for every $\beta\in[\frac12,1]$ if and only if $p\ge k_n$, and $M_T(s)$ has property $\mathcal P$ for every $s\in\mathbb R$ if and only if $p\ge k_n$ (Theorems~\ref{thm:dichotomy} and~\ref{thm:threshold}).
\item \emph{Threshold.} If $p\le k_n-1$, then $B_\beta(T)$ fails property $\mathcal P$ for every $\beta\in[\frac12,\frac{n+1}{2n+1}]$, and $M_T(s)$ fails it for every $|s|\ge n$ (Theorem~\ref{thm:threshold}; note that $\alpha=\frac{n}{n+1}$ corresponds to $\beta=\frac{n+1}{2n+1}$ and to $|s|=n$).
\item \emph{The deformation parameter of the crossing.} Let $p\le k_n-1$ and put
\[
t_0(T)=\frac{\alpha_0(T)}{1-\alpha_0(T)}=\frac{1}{1-\alpha_0(T)}-1,\qquad \beta_0(T)=\frac{1}{1+\alpha_0(T)} .
\]
Since $S(T)\subseteq[\frac12,\frac{n}{n+1})$ is compact and the reparametrisations are monotone, $t_0(T)$ is the largest $t\ge1$ for which $M_T(t)$ has property $\mathcal P$, and the property fails for all $|s|>t_0(T)$; similarly, $\beta_0(T)$ is the smallest $\beta\in[\frac12,\frac23]$ for which $B_\beta(T)$ has property $\mathcal P$, and the property fails for all $\beta\in[\frac12,\beta_0(T))$.
\end{itemize}
In terms of the deformed Laplacian the main results of the paper take a particularly transparent form, since they are all estimates for $\frac{1}{1-\alpha_0}=1+t_0$. Writing $t_{\max}(n)$ and $t_{\max}(n,p)$ for the maxima of $t_0(T)$ over the trees with $n$ vertices and at most $k_n-1$ leaves, respectively with exactly $p$ leaves, we obtain: $t_0(T)<n$ for every non-exceptional tree (Theorem~\ref{thm:threshold}); $n-\frac{2\pi^2n}{(n-1)^2}\le n\cos\frac{\pi}{\lfloor n/2\rfloor}\le t_0(S_n^\ast)\le t_{\max}(n)<n$ (Theorem~\ref{thm:spider} and Corollary~\ref{cor:asymp}); $t_{\max}(n)>n-1$ for $n=16$ and for all $n\ge18$ (Corollary~\ref{cor:false}); $t_0(T)\le n\kappa_n(p)+\min\{3/\kappa_n(p),\sqrt{3n}\}$ for every tree with $p\le k_n-1$ leaves (Theorem~\ref{thm:leaves}); and $n\kappa_n(p)-\frac76\le t_{\max}(n,p)\le n\kappa_n(p)+6$ under the hypotheses of Corollary~\ref{cor:twosided} when $k_n-p$ divides $D_n(p)$. In particular, in the situation of Corollary~\ref{cor:profile},
\[
\cos\frac{\pi}{\lceil\Lambda(\beta)\rceil-1}\ \le\ \liminf_{n\to\infty}\frac{t_{\max}(n,p_n)}{n}\ \le\ \limsup_{n\to\infty}\frac{t_{\max}(n,p_n)}{n}\ \le\ \cos\frac{\pi}{\Lambda(\beta)},
\]
so that the largest deformation parameter $|s|$ for which some tree with $p$ leaves still has property $\mathcal P$ is approximately $n\cos\frac{\pi}{\Lambda(p/n)}$. The corresponding statements for $B_\beta$ are obtained from the identity $\frac{1}{1-\alpha_0(T)}=\frac{\beta_0(T)}{2\beta_0(T)-1}$; for instance, Theorem~\ref{thm:threshold} gives $\beta_0(T)>\frac{n+1}{2n+1}$ for every non-exceptional tree. Finally, as already noted in~\cite[Section~10]{BO}, the counterexamples to $\mathcal P_\alpha$ constructed there transfer to $B_\beta$ and to $M_T(s)$ in the same way.
\end{remark}

\section*{Data availability}
No datasets were generated or analysed in this study. The values reported in Figure~\ref{fig:Nalpha} are roots of the explicit polynomial $\Phi_T$ of Lemma~\ref{lem:firstcrossing}, computed to the precision shown. The numerical values reported in Example~\ref{ex:n7} and in Table~\ref{tab:alphamax} were obtained by exhaustive computation over all trees with the indicated number of vertices, following the procedure described in Remark~\ref{rem:table}.

\section*{Declaration of competing interest}
The authors declare no competing interests.

\section*{Declaration of AI usage}
AI tools were used to assist with outlining, checking arguments, searching for examples, proofreading and revising the manuscript. The authors reviewed, edited and verified all content and take full responsibility for the final text.

\end{document}